# LQG online learning

**Giorgio Gnecco** GIORGIO.GNECCO@IMTLUCCA.IT
*DYSCO Research Unit - IMT School for Advanced Studies*
*Piazza S. Francesco, 19 - 55110 Lucca, Italy*

**Alberto Bemporad** ALBERTO.BEMPORAD@IMTLUCCA.IT
*DYSCO Research Unit - IMT School for Advanced Studies*
*Piazza S. Francesco, 19 - 55110 Lucca, Italy*

**Marco Gori** MARCO@DIISM.UNISI.IT
*DIISM Department - University of Siena*
*Via Roma, 56 - 53100 Siena, Italy*

**Marcello Sanguineti** MARCELLO.SANGUNETI@UNIGE.IT
*DIBRIS Department - University of Genoa*
*Via Opera Pia, 13 - 16145 Genova, Italy*

## Abstract

Optimal control theory and machine learning techniques are combined to formulate and solve in closed form an optimal control formulation of online learning from supervised examples with regularization of the updates. The connections with the classical Linear Quadratic Gaussian (LQG) optimal control problem, of which the proposed learning paradigm is a non-trivial variation as it involves random matrices, are investigated. The obtained optimal solutions are compared with the Kalman-filter estimate of the parameter vector to be learned. It is shown that the proposed algorithm is less sensitive to outliers with respect to the Kalman estimate (thanks to the presence of the regularization term), thus providing smoother estimates with respect to time. The basic formulation of the proposed online-learning framework refers to a discrete-time setting with a finite learning horizon and a linear model. Various extensions are investigated, including the infinite learning horizon and, via the so-called "kernel trick", the case of nonlinear models.

## 1. Introduction

In recent years, the combination of techniques from the fields of optimization/optimal control and machine learning has led to a successful interaction between the two disciplines. The cross-fertilization between these two fields shows itself in both directions.

### 1.1 Application of machine-learning techniques to optimization/optimal control

Sparsity-inducing regularization techniques from machine learning have been exploited to find suboptimal solutions to an initially unregularized optimization problem, having at the same time a sufficiently small number of nonzero arguments. For instance, the Least

Absolute Shrinkage and Selection Operator (LASSO) [49] was applied in [28] to consensus problems, and in [23] to Model Predictive Control (MPC).

Applications of machine-learning techniques to control can be found, e.g., in [48], and in the series of papers [21, 22, 29], where Least Squares Support Vector Machines (LS-SVMs) and one-hidden-layer perceptron neural networks, respectively, were applied to find suboptimal solutions to optimal control problems. In [36], spectral graph theory methods - already exploited successfully in machine-learning problems [5] - were applied to the control of multi-agent dynamical systems.

Least Squares Support Vector Machines and spectral graph theory have been also applied, respectively, to system identification [32] and control of epidemics [12].

### 1.2 Application of optimization/optimal-control techniques to machine learning

This is the direction followed in the present work: we develop and approach that exploits for machine learning techniques from optimization and optimal control.

Specifically, we propose and solve in closed form an optimal-control formulation of online learning with supervised examples and regularization of the updates. In the online framework, the examples become available one by one as time passes and the training of the learning machine is performed continuously. Online learning problems have been investigated, e.g., in [33, 38, 43, 44, 51, 52], but without using an approach based on optimal control theory. As suggested by the preliminary results that we obtained in [24], such an approach can provide a strong theoretical foundation to the choice of a specific online learning algorithm, by selecting the parameter updates as the outputs of a sequence of control laws that solve a suitable optimal control problem modeling online learning itself[1]. A distinguishing feature of our study is that we derive online learning algorithms as closed-form optimal solutions to suitable online learning problems. In contrast, typically, works in the literature propose a certain algorithm and then investigate its properties, but do not analyze the optimality of such an algorithm with respect to a suitable online learning problem. An exception is [8], but it refers to a deterministic optimization problem and, differently from our approach, it does not contain any regularization of the updates.

In a nutshell, our contributions are the following:

- we make the machine-learning community aware of a point of view that till now might have been overlooked;

- by exploiting such viewpoint, we develop a novel machine-learning paradigm, for which we provide closed-form solutions;

- we make connections between our results and other machine-learning algorithms.

### 1.3 The adopted learning model

The learning model that we adopt can be considered a nontrivial variation (due to the presence of suitable random matrices) of the Linear Quadratic (LQ) and Linear Quadratic Gaussian (LQG) optimal control problems, which we briefly summarize in the following. The LQ problem [7] consists in the minimization - with respect to a set of control laws,

1. The results from [24] correspond, in the present work, to a subset of the results contained in Section 3.

one for each decision stage - of a convex quadratic cost related to the control of a linear dynamical system, which is decomposed into the summation of several convex quadratic per-stage costs, associated with *a-priori given* cost matrices. At each stage, a control law is applied. It is a function of an information vector, which collects all the information available to the controller up to that stage. More precisely, the information vector is formed by the sequence of controls applied to the dynamical system up to the previous stage, and by the sequence of measures of the state of the dynamical system itself, acquired up to the current stage. A peculiarity of the LQ problem is that such measures are linearly related to the state, again through suitable *a-priori given* measurement matrices. The measures may be corrupted by additive noise, with given covariance matrices. When all the noise vectors are Gaussian, one obtains the LQG problem, for which closed-form optimal control laws in feedback form are known[2]. They are computed by solving recursively suitable Riccati equations and applying the Kalman filter [47] to estimate the current state of the dynamical system.

The main difference between the LQ/LQG problems and the proposed formulation of online learning with supervised examples is the following. In our approach both the cost and measurement matrices are *random*, being associated with the input examples, which become available as time goes on. It is worth mentioning that randomness of some matrices in the context of the LQ optimal control problem was considered also in [7, Section 4.1], but in a way not directly applicable to the online learning problem investigated in this paper (see Remark [7, Section 4.1] for further details). First we consider a linear relationship between the input examples and their labels, possibly corrupted by additive noise, and collect into the state vector both the current estimate of the parameter vector modeling the input-output relationship, and the parameter vector itself, which is unknown. Then we relax the linearity assumption and address a more general nonlinear context. The goal consists in finding an optimal online estimate of the parameter vector, on the basis of the information associated with the incoming examples, modeled in the simplest case as independent and identically distributed random vectors[3].

Each decision stage corresponds to the presentation of one example to the learning machine, whereas the convex per-stage cost penalizes quadratically the difference between the observed output and the one predicted by the learning machine, by using the current estimate of the parameter vector. Causality in the updates (i.e., their independence on future examples, which is important for a truly "online" framework) is preserved by constraining the updates to depend only on the "history" of the observations and updates up to the decision time, likewise in the LQ/LQG problems.

At each stage, the error on future examples is taken into account through the conditional expectation of the summation of the associated per-stage costs, conditioned on the current information vector. The link between the examples used for training and the future examples is only in their common generation model. In order to reduce the influence of outliers on the online estimate of the parameter vector, its smoothness with respect to time is also enforced through the presence of a suitable regularization term in the functional

2. Specifically, as functions of an estimate of the current state of the dynamical system.
3. This framework is also extended in the paper to other probability models for the generation of the examples (see Section 8).

| Update | Control |
|---|---|
| Updating function | Control function |
| Problem OLL (On-Line Learning) | Optimal control problem |
| Learning horizon | Optimization horizon |
| Learning functional | Cost functional |
| Average learning functional | Average cost functional |

Table 1: Some correspondences between the machine-learning terminology and the optimal-control one.

to be optimized, weighted by a positive regularization parameter[4]. The optimal solution is obtained by applying Dynamic Programming (DP) and requires the solution of suitable Riccati equations. The above-mentioned difference between the classical LQ/LQG problems and the proposed online learning framework (i.e., the random nature of the matrices) determines two different forms for such equations, for the backward and forward phases of DP, respectively. When the optimization horizon is infinite, it is necessary to take into account the random nature of the matrices to perform a convergence analysis of the online estimate of the unknown parameter vector.

Table 1 provides the correspondence between the notations used for optimal control, and the ones used for the proposed online learning framework.

### 1.4 Relationships with other machine-learning techniques

The approaches to online learning most closely related to this work are Kalman filtering [47] (see also [7, Appendix E]) and its kernel version [33,38], in which, however, no penalization is made directly on the control (updating) variables. Indeed, one of our contributions consists in developing a theoretical framework in which such a penalization is taken into account and in providing in most cases closed-form results. Interestingly, the obtained solutions can be interpreted as smoothed versions (with respect to time) of the solution obtained applying the Kalman filter only. Most importantly, we show, both theoretically and numerically, that our solutions are less sensitive to outliers than the Kalman-filter estimates. This is very useful, e.g., if one wants to give more importance to a whole set of most recently presented examples than to the current example, allowing to obtain estimates that change more smoothly with respect to time (smoothness of an estimate is a desirable property, e.g., in applications to online system identification and control, in which one has also to control the system just identified).

The updating formula that provides the solution to the proposed learning paradigm is similar to the one of other online estimates obtained through various machine learning techniques, such as stochastic gradient descent. However, there is a substantial difference: we derive it as the optimal solution of an optimization problem modeling online learning,

4. We shall present a comparison with the sequence of Kalman-filter estimates of the unknown parameter vector that shows the larger smoothness and less sensitivity to outliers of the sequences of estimates obtained solving the proposed optimal control formulations of online learning (see Section 5).

and this allows us to prove various interesting properties. We believe that this approach could be fruitfully applied also to other machine learning techniques used in online learning.

A number of extensions is also described with some detail at the end of the paper, providing hints for further research in several directions, and showing the generality of the basic theoretical framework studied in the paper.

### 1.5 Organization of the paper

Section 2 is a non-technical overview of the main results derived in the paper, written to allow readers who are not familiar with optimal control, but work in the field of machine learning, to appreciate the nature of our approach and its contributions. At the same time, it provides a summary of the main results of the paper. Section 3 introduces and solves the proposed model of online learning as an LQ optimal control problem with random matrices and finite optimization horizon, and provides closed-form expressions for the optimal solutions in the LQG case. Section 4 extends the analysis to the infinite-horizon framework. Section 5 investigates convergence properties of the algorithm, whereas Section 6 compares the proposed online approach with average regret minimization and the Kalman-filter estimates, both theoretically and numerically. Section 7 extends the analysis to nonlinear models (kernel methods). Other extensions are described in Section 8. Section 9 is a conclusive discussion. To improve the readability, most technical proofs are contained in the Appendix.

## 2. Overview of the main results

In the following, we summarize the main results with links to the parts of the paper in which they are presented, providing a guidance to the reading of the paper.

- *We derive closed-form optimal solutions* for the proposed optimal control formulation of online learning, and for some of its extensions. They are expressed in terms of two Riccati equations (see Section 3), associated, respectively, with the backward phase of DP (to determine the gain matrix of the optimal controller) and with the determination of the gain matrix in the Kalman filter (in the case of Gaussian random vectors). Differently from the LQG problem, the two Riccati equations have different natures: one involves expectations of random matrices (so, it may be called an "average Riccati equation"), whereas the other involves realizations of random matrices (so, it may be called a "stochastic Riccati equation"). As a consequence, a specific study - detailed in the paper - is needed to study properties of their solutions, which confirms that the proposed problem is not a trivial application of the LQG problem to an online learning framework.

- *We analyse both theoretically and numerically the role of the regularization parameter* (see Subsection 3.3).

- In the *infinite-horizon case, we investigate the existence of a stationary (and linear) optimal updating function* (see Section 4), *stability issues, and the convergence to* 0 *of the mean-square error between the parameter vector and its online estimate* when the number of examples goes to infinity (see Section 5). In this context, another non-trivial difference with respect to the classical LQG problem arises: when computing certain expectations conditioned on the current information vector, one has to take into account that the information vector at a generic stage has additional components deriving from the knowledge of the sequence of output measurement matrices up to the stage itself (as these are random matrices, associated with the input examples). As a consequence, the Kalman gain matrix, which is shown in the paper to be embedded in the optimal solution, is not only stage-dependent but also a random matrix (although it becomes deterministic when conditioned on the input examples already presented to the learning machine). This motivates the investigation of issues such as its convergence in probability and the convergence of its expectation when the number of examples goes to infinity.

· *We discuss the connection of the proposed online learning framework with average regret minimization.* We prove that the sequence of our estimates minimizes the average regret functional (see Subsection 6.1).

· *We investigate the connections between our solution with the Kalman-filter estimate and stochastic gradient descent* (Remark 3 and Subsection 6.2). We prove that our solution can be interpreted as a smoothed Kalman-filter estimate, with time-varying gain matrix, and we show that it outperforms the latter in terms of its larger smoothness (Subsection 6.3; see also Section 8 *e)*) and its smaller sensitivity to outliers (Subsection 6.4).

· *We discuss cases in which the proposed solutions can be computed efficiently* (see, e.g., the comments presented after Proposition 2, Remark 10, and the numerical results reported in Subsection 6.3).

· *We address the case of nonlinear input-output relationships, modeled using the "kernel trick" of kernel machines* (see Section 7). As is well-known, the kernel trick is based on a preliminary (in general nonlinear) mapping of the input space to a larger-dimensional feature space, to which the original linear model is applied in a second step. The "kernel trick", which consists in the computation of certain inner products in the auxiliary feature space through a suitable function called "kernel", can be applied in our context since we show that the optimal solution can be expressed in terms of inner products in the feature space that can be computed through a kernel.

· *We describe various other possible extensions* (see Section 8), such as the case of a time-varying parameter vector to be learned, the introduction of a discount factor, the inclusion of additional regularization terms, a continuous-time extension framework, and a possible extension of the problem formulation through techniques from robust estimation and control.

Table 2 collects some acronyms of frequent use in the paper.

| | |
|---|---|
| **Problem $\text{OLL}_\gamma^N$** | On-Line Learning Problem over finite horizon $N$ and with regularization parameter $\gamma$ |
| **Problem $\text{OLL}_\gamma^\infty$** | Online Learning Problem over infinite horizon and with regularization parameter $\gamma$ |
| **OLL estimate** | Estimate obtained solving Problem $\text{OLL}_\gamma^\infty$ or $\text{OLL}_\gamma^N$ |
| **LQ** | Linear Quadratic |
| **LQG** | Linear Quadratic Gaussian |
| **LQR** | Linear Quadratic Regulator |
| **ARE** | Average Riccati Equation |
| **SRE** | Stochastic Riccati Equation |
| **KF** | Kalman Filter |
| **MSE** | Mean-Square Error |

Table 2: Acronyms of frequent use.

## 3. The basic case: discrete-time, finite horizon, and linear model

For simplicity, we consider first a discrete-time setting with a finite learning horizon and a linear model. Then, we shall address the extensions to an infinite learning horizon and nonlinear models.

### 3.1 Problem formulation

**Assumption 1 (Linear data generation model)** *At each time $k = 0, 1, \ldots$, a learning machine can observe the supervised pair $(x_k, y_k)$, where $x_k \in \mathbb{R}^d$ is a column vector and $y_k \in \mathbb{R}$. The output $y_k$ is generated from the input $x_k$ according to the following linear model:*

$$y_k = w' x_k + \varepsilon_k \,, \tag{1}$$

*where $\varepsilon_k \in \mathbb{R}$ is a measurement noise, and $w \in \mathbb{R}^d$ is a random vector, unknown to the learning machine, and to be estimated by the learning machine itself by using the sequence of examples $(x_k, y_k)$ as they become available.*

**Assumption 2 (Random variables)** *The random variables $w$, $\{x_k\}$, $\{\varepsilon_k\}$ are mutually independent*[5] *and (only for simplicity of notation and without any loss of generality) have mean 0. The random variables $\varepsilon_k$ have the same variance $\sigma_\varepsilon^2$, and each $x_k$ has finite covariance matrix $\mathop{\mathbb{E}}_{x_k}\{x_k x_k'\}$.*

**Assumption 3 (Learning machine)** *Starting from the initial estimate $\hat{w}_0 := 0$ of $w$, at each time $k + 1 = 1, 2, \ldots$, the learning machine builds an estimate $\hat{w}_k$ of $w$, generated*

5. As another extension, one could consider the case in which the inputs $x_k$ are generated by the learning machine as the states of another controlled dynamical system. This, together with the optimization of a suitable learning functional similar to (5), would model the problem of *online active learning*, as the learning machine would have an influence even on the choice of the sequence of input examples (see item *n)* in Section 8).

*according to*

$$\hat{w}_{k+1} = \hat{w}_k + u_k \,, \qquad (2)$$

*where $u_k$ is the* update *of the estimate of $w$ at the time $k$ (to be optimized according to a suitable optimality criterion, defined later on).*

**Remark 1** It is important to observe that the model (1) is time-invariant[6], in the sense that the same $w$ is used to generate every $y_k$, starting from every $x_k$ and every $\varepsilon_k$. So, once a realization of the random vector $w$ has been generated, this can be interpreted as a *fixed* vector, to be estimated by the learning machine using the online supervised examples.

To analyze the time-evolution of the estimate, one has to consider the following dynamical system (see [2] for a similar approach), with state vector $(w'_k, \hat{w}'_k)'$ and initial conditions $w_0 := w$ and $\hat{w}_0 := 0$:

$$\begin{cases} w_{k+1} = w_k \,, \\ \hat{w}_{k+1} = \hat{w}_k + u_k \,, \end{cases} \qquad (3)$$

together with the measures

$$y_k = C_k w_k + \varepsilon_k \,, \qquad (4)$$

where $C_k := x'_k$.

**Assumption 4 (Available information and updating functions)** *The* update $u_k$ *at the time $k$ is chosen as a function $u_k(I_k)$, called* updating function*, of the* information vector *$I_k$ at the same time, which collects the "history" up to the time $k$, and is defined as*

$$I_k := \{(x_j, y_j) \text{ for } j = 0, \dots, k, \text{ and } u_j \text{ for } j = 0, \dots, k-1\}$$

*for $k = 1, 2, \dots$, and*

$$I_0 := \{(x_0, y_0)\} \,.$$

Hence, the update $u_k$ depends only on the sequence of examples $(x_j, y_j)$ observed up to the current stage $k$ and on the sequence of previous updates $u_j$ (or equivalently, since $\hat{w}_0 = 0$, on the sequence of previous updates of the estimate of $w$).

In our model, the updating functions $u_k$ are chosen in order to minimize a *learning functional* over a finite learning horizon, defined as follows,

**Definition 1 (Learning functional over horizon $N$)** *Let $N$ be a positive integer, $\gamma > 0$, $Q_k := x_k x'_k$, and*

$$\begin{aligned} & J_\gamma^N \left( \{u_k(I_k)\}_{k=0}^{N-1} \right) \\ := & \mathop{\mathbb{E}}_{w, \{x_k\}_{k=0}^N, \{\varepsilon_k\}_{k=0}^{N-1}} \left\{ \sum_{k=0}^{N-1} \left[ \left((\hat{w}_k - w_k)' x_k\right)^2 + \gamma u'_k u_k \right] + \left((\hat{w}_N - w_N)' x_N\right)^2 \right\} \\ = & \mathop{\mathbb{E}}_{w, \{x_k\}_{k=0}^N, \{\varepsilon_k\}_{k=0}^{N-1}} \left\{ \sum_{k=0}^{N-1} \left[ (\hat{w}_k - w_k)' Q_k (\hat{w}_k - w_k) + \gamma u'_k u_k \right] + (\hat{w}_N - w_N)' Q_N (\hat{w}_N - w_N) \right\} . \end{aligned} \qquad (5)$$

6. An extension to the case of a (slowly) time-varying parameter vector will be discussed in item *e)* of Section 8.

We state the following On-Line Learning Problem (in the paper, the symbol "$\circ$" is used to denote optimality).

**Problem $\text{OLL}_\gamma^N$ (On-Line Learning over a finite horizon)** *Given the finite learning horizon $N$, the examples $(x_k, y_k)$ generated at each time instant $k = 0, 1, \ldots, N$ according to the model defined by Assumptions 1 and 2, and the learning machine defined by Assumption 3, find the finite sequence $u_0^\circ(I_0), \ldots, u_{N-1}^\circ(I_{N-1})$ of* optimal updating functions *with the structure defined by Assumption 4, that minimizes the learning functional (5).*

Problem $\text{OLL}_\gamma^N$ can be considered a *parameter identification problem* or an *optimal estimation problem*, as the final goal consists in estimating the parameter vector $w$ relating input examples to their outputs, given the current subsequence of examples and the adopted optimality criterion. It can also be considered an *optimal control problem*, interpreting the updating function $u_k$ as a control function for the dynamical system (3). Although this last interpretation may seem less natural than the first two, it is motivated by the fact that Problem $\text{OLL}_\gamma^N$ and its variations presented later in the paper can be investigated using optimal control techniques, as it is done in the following.

For every $k = 0, 1, \ldots, N$, we shall call $\hat{w}_k$ *online estimate* (*OLL estimate*, for short).

**Remark 2** The term $((\hat{w}_k - w_k)' x_k)^2$ in the learning functional (5) penalizes a large deviation of the learning machine estimate $\hat{w}_k' x_k$ of the label $y_k$ from its best estimate (in a mean-square sense) $w_k' x_k = w' x_k$ which would have been obtained if $w$ were known, whereas the term $u_k' u_k$ penalizes a large square of the norm of the update $u_k$ of the estimate of $w$, and $\gamma$ is a regularization term, which measures the trade-off between the two terms.

**Remark 3** The OLL estimates correspond to the limit case $\gamma = 0$ in the formulation of Problem $\text{OLL}_\gamma^N$. Indeed, in such a case, each term

$$\mathop{\mathbb{E}}_{w, \{x_t\}_{t=0}^{k}, \{\varepsilon_t\}_{t=0}^{k-1}} \left\{ (\hat{w}_k - w_k)' Q_k (\hat{w}_k - w_k) \right\}$$

in (5) is minimized when $\hat{w}_k$ is the conditional expectation of $w_k$ given $I_{k-1}$, i.e., when it is the Kalman-filter estimate of $w_k$ at time $k-1$ (see, e.g., [7, Proposition E.1])[7]. It is worth mentioning that, since the parameter vector to be learned is constant and the data generation model is described by equation (1), the specific Kalman-estimation problem is equivalent to recursive least squares (see [41, Section 12.A] for a proof of this equivalence).

In Subsections 6.3 and 6.4, we discuss some relationships of the proposed learning framework with the classical Kalman filter [47]. As it will be shown by Proposition 8 and by the numerical results in Figure 1, the presence of the regularization term in the learning functional (5) can make the resulting sequence of optimal estimates of $w$ smoother with respect to the time index, and less sensitive to outliers, than the sequence of estimates obtained by using the classical Kalman filter, under Gaussian assumptions on the random variables $w$ and $\varepsilon_k$.

7. Note that [7, Proposition E.1] is formulated in terms of the square of the Euclidean norm of the error vector, which is $\hat{w}_k - w_k$ in our case. However, the proposition can be still applied if one moves from the square of the Euclidean norm to the square of the (semi)norm induced by $Q_k$, or to its expectation (as in the present case).

**Remark 4** The constraint that each update $u_k$ (hence also each updating function $u_k^\circ$) depends only on the sequence of examples $(x_j, y_j)$ observed up to the current stage $k$ and on the sequence of previous updates $u_j$, implies that no future examples are taken into account to update the current estimate of $w$. Hence, the proposed solution is actually a model of *online learning*. Instead, *batch learning* corresponds to the case where one assumes that all the sequence $\{(x_j, y_j), j = 0, \ldots, N\}$ of examples is known to the learning machine, starting from the time $k = 0$.

**Remark 5** An alternative definition of the learning functional can be obtained by replacing the term $((\hat{w}_k - w_k)' x_k)^2$ in (5) by $(\hat{w}_k' x_k - y_k)^2$, i.e., by the square of the difference between the label estimated by the learning machine before measuring $y_k$ (but knowing $x_k$), and the label $y_k$ generated by the model (1) at the time $k$ (note that, differently from the term $w_k' x_k$, they are both observable at the time $k$). However, by taking expectations and recalling that $\varepsilon_k$ has mean 0 and is mutually independent from $x_k$, $w_k$, and $\hat{w}_k$, one obtains

$$\begin{aligned}
& J_{\gamma,y}^N \left( \{u_k(I_k)\}_{k=0}^{N-1} \right) \\
:= & \mathop{\mathbb{E}}_{w, \{x_k\}_{k=0}^N, \{\varepsilon_k\}_{k=0}^N} \left\{ \sum_{k=0}^{N-1} \left[ \left(\hat{w}_k' x_k - y_k\right)^2 + \gamma u_k' u_k \right] + \left(\hat{w}_N' x_N - y_N\right)^2 \right\} \\
= & \mathop{\mathbb{E}}_{w, \{x_k\}_{k=0}^N, \{\varepsilon_k\}_{k=0}^N} \left\{ \sum_{k=0}^{N-1} \left[ (\hat{w}_k - w_k)' Q_k (\hat{w}_k - w_k) + \gamma u_k' u_k \right] + (\hat{w}_N - w_N)' Q_N (\hat{w}_N - w_N) \right\} \\
& + (N+1)\sigma_\varepsilon^2 \,. \qquad (6)
\end{aligned}$$

Hence, since the last term in (6) is a constant, the learning functionals (5) and (6) have the same sequence of optimal updating functions. It is worth noting that, in both formulas (5) and (6), in order to generate the estimates $\hat{w}_k$, one uses only the probability distribution of $w_k$ conditioned on the already available observations.

The statement of Problem $\mathrm{OLL}_\gamma^N$ can be simplified by defining the *learning error*

$$e_k := \hat{w}_k - w_k \,,$$

which evolves according to

$$e_{k+1} = e_k + u_k \,, \qquad (7)$$

where

$$e_0 := -w_0 = -w \,.$$

Of course, $e_k \simeq 0$ means $\hat{w}_k \simeq w_k = w$. Moreover, since both $\hat{w}_k$ and $x_k$ are known at the time $k$, one can replace the measures $y_k$ by

$$\tilde{y}_k := \hat{w}_k' x_k - y_k \,,$$

hence obtaining the measurement equation

$$\tilde{y}_k = C_k e_k + \tilde{\varepsilon}_k \,, \tag{8}$$

where $\tilde{\varepsilon}_k := -\varepsilon_k$, and has the same variance $\sigma_\varepsilon^2$ as $\varepsilon_k$. In this case, the "history" of the learning machine, measures, and past updates up to the time $k$ is collected in the new information vector $\tilde{I}_k$, defined as

$$\tilde{I}_k := \{(x_j, \tilde{y}_j)\,\text{for}\, j = 0, \ldots, k, \,\text{and}\, u_j \,\text{for}\, j = 0, \ldots, k-1\}$$

for $k = 1, 2, \ldots$, and

$$\tilde{I}_0 := \{(x_0, \tilde{y}_0)\} \,.$$

There is a one-to-one correspondence between the information vectors $I_k$ and $\tilde{I}_k$. So, the optimization of the learning functional (5) assuming that the dynamical system evolves according to equation (3), the sequence of measures is provided by equation (4), and the updating functions $u_k$ have the form $u_k(I_k)$, is equivalent to the optimization of the following learning functional:

$$\begin{aligned}
&\tilde{J}_\gamma^N \left( \left\{ u_k(\tilde{I}_k) \right\}_{k=0}^{N-1} \right) \\
:= &\mathop{\mathbb{E}}_{e_0, \{x_k\}_{k=0}^N, \{\tilde{\varepsilon}_k\}_{k=0}^{N-1}} \left\{ \sum_{k=0}^{N-1} \left[ (e_k' x_k)^2 + \gamma u_k' u_k \right] + (e_N' x_N)^2 \right\} \\
= &\mathop{\mathbb{E}}_{e_0, \{x_k\}_{k=0}^N, \{\tilde{\varepsilon}_k\}_{k=0}^{N-1}} \left\{ \sum_{k=0}^{N-1} \left[ e_k' Q_k e_k + \gamma u_k' u_k \right] + e_N' Q_N e_N \right\} \,,
\end{aligned} \tag{9}$$

assuming that the error vector evolves according to equation (7), the sequence of measures is provided by equation (8), and the update $u_k$ is now a function $u_k(\tilde{I}_k)$ of the information vector $\tilde{I}_k$. Such a problem is a non-trivial variation of the classical LQ problem [7, Section 5.2]. Whereas in the latter the matrices $C_k$ and $Q_k$ are deterministic, in the proposed formulation of online learning they are random, since they depend on the input examples $x_k$. Another difference is that, for $j = 0, \ldots, k$, the information vector $\tilde{I}_k$ includes the realizations of the inputs $x_j$, hence also of the matrices $C_j$ and $Q_j$.

### 3.2 Solution of the finite-horizon online learning problem

To solve Problem $\text{OLL}_\gamma^N$, we make an extensive use of the concept of *cost-to-go function* from the theory of dynamic programming [7, Chapter 1]. In our context, the cost-to-go function $\tilde{J}_k^\circ$ at the time stage $k = 0, \ldots, N-1$ is defined as

$$\tilde{J}_k^\circ(\tilde{I}_k) := \inf_{\{u_j(\tilde{I}_j)\}_{j=k}^{N-1}} \mathop{\mathbb{E}}_{e_k, \{x_j\}_{j=k+1}^N, \{\tilde{\varepsilon}_j\}_{j=k+1}^{N-1}} \left\{ \sum_{j=k}^{N-1} [e_j' Q_j e_j + \gamma u_j' u_j] + e_N' Q_N e_N \middle| \tilde{I}_k \right\} \,, \tag{10}$$

whereas

$$\tilde{J}_N^\circ(\tilde{I}_N) = \mathop{\mathbb{E}}_{e_N} \left\{ e_N' Q_N e_N \middle| \tilde{I}_N \right\} \,. \tag{11}$$

Finally, the optimal value of the learning functional (9) is

$$\tilde{J}_0^\circ = \mathop{\mathbb{E}}_{\tilde{I}_0} \left\{ \tilde{J}_0^\circ(\tilde{I}_0) \right\} .$$

Under mild regularity conditions (see the next Remark 6), the cost-to-go functions can be determined recursively by solving the *Bellman Equations*

$$\tilde{J}_k^\circ(\tilde{I}_k) \quad = \quad \inf_{u_k \in \mathbb{R}^d} \mathop{\mathbb{E}}_{e_k, \tilde{I}_{k+1}} \left\{ e_k' Q_k e_k + \gamma u_k' u_k + \tilde{J}_{k+1}^\circ(\tilde{I}_{k+1}) \big| \tilde{I}_k, u_k \right\} \tag{12}$$

for $k = N-1, \ldots, 0$.

**Remark 6** The regularity conditions mentioned above are satisfied in the case - studied in the paper - where the random vectors $w$ and $\varepsilon_k$ are Gaussian. Indeed, in such a context the optimal updating functions that will be provided by (13) are linear with respect to the information vector [7, Section 1.5], [9].

Equations (11) and (12) are similar to those for the cost-to-go functions in the LQ problem (see, e.g., [7, Section 5.2]), with the difference that in the present context the matrices $Q_k$ and $C_k$ are random. Moreover, both matrices become known to the learning machine at the time $k$, as they can be derived from the information vector $\tilde{I}_k$. In the following, we use sometimes the superscript "$\circ$" not only for the optimal updating functions, but also to denote vectors (e.g, $\hat{w}_k$ and $e_k$) evaluated when the sequence of optimal updating functions (13) is applied.

**Proposition 1 (Optimal updating functions and Average Riccati Equation (ARE))** *Let Assumptions 1, 2, 3, and 4 be satisfied. Then, the updating functions that solve Problem* $OLL_\gamma^N$ *are given, for* $k = N-1, \ldots, 0$, *by*

$$u_k^\circ(\tilde{I}_k) = L_k \mathop{\mathbb{E}}_{e_k^\circ} \left\{ e_k^\circ \big| \tilde{I}_k \right\} , \tag{13}$$

*where*

$$L_k := -(\overline{K}_{k+1} + \gamma I)^{-1} \overline{K}_{k+1} , \tag{14}$$

*and the matrices*

$$K_k := \overline{K}_{k+1} - \overline{K}_{k+1} (\overline{K}_{k+1} + \gamma I)^{-1} \overline{K}_{k+1} + Q_k , \tag{15}$$

$$F_k := \overline{K}_{k+1} (\overline{K}_{k+1} + \gamma I)^{-1} \overline{K}_{k+1} , \tag{16}$$

*and*

$$\overline{K}_k := \mathop{\mathbb{E}}_{K_k} \{K_k\} = \overline{K}_{k+1} - \overline{K}_{k+1} (\overline{K}_{k+1} + \gamma I)^{-1} \overline{K}_{k+1} + \mathop{\mathbb{E}}_{Q_k} \{Q_k\} \tag{17}$$

*are symmetric positive-semidefinite. The recursions above are initialized by*

$$K_N := Q_N \tag{18}$$

*and*

$$\overline{K}_N := \mathop{\mathbb{E}}_{K_N} \{K_N\}. \tag{19}$$

Equation (17) can be called an "*Average Riccati Equation*" (*ARE*, for short), since it contains the expectation term $\mathop{\mathbb{E}}_{Q_k}\{Q_k\}$. In practice, it can be solved likewise the classical deterministic Riccati equation of the Linear Quadratic Regulator (LQR) subproblem [7, Section 5.2], simply by replacing $Q_k$ (which is deterministic in the LQ problem) by $\mathop{\mathbb{E}}_{Q_k}\{Q_k\}$. It is worth remarking that solving the ARE (17) does not require the knowledge of future input examples, and that all the matrices $L_k$ in (14) have spectral radius[8] $|\lambda|_{\max}(L_k)$ strictly smaller than 1. Finally, the matrices $F_k$ are reported in formula (16) because they are used to express $\tilde{J}_k^\circ(\tilde{I}_k)$ (see formula (106) in the Appendix). They are also used in the infinite-horizon version of Problem $\text{OLL}_\gamma^N$ (Problem $\text{OLL}_\gamma^\infty$), to reduce one part of the proof of Proposition 4 in Section 4 to the finite-horizon case.

Due to (13), in order to generate the optimal update $u_k^\circ(\tilde{I}_k)$ at the time $k$ one has to compute $\mathop{\mathbb{E}}_{e_k^\circ}\left\{e_k^\circ\middle|\tilde{I}_k\right\}$. Let us now make the following additional assumption.

**Assumption 5 (Gaussian random variables)** *The random variables $w$ and $\varepsilon_k$ are Gaussian.*

The next proposition shows that, when the additional Assumption 5 is satisfied, the optimal estimate $\hat{w}_k^\circ$ of the proposed framework tracks the (usually time-varying) Kalman-filter estimate. Indeed, inspection of its proof shows that

$$\hat{e}_k^{\circ,\dagger} := \mathop{\mathbb{E}}_{e_k^\circ}\left\{e_k^\circ\middle|\tilde{I}_k\right\}$$

is the *Kalman-Filter* (*KF estimate*, for short) of the error vector $e_k^\circ$ at the time $k$, based on the information vector $\tilde{I}_k$, thus getting a Kalman-filter recursion scheme.

In the following, we denote by

$$\hat{w}_k^\dagger := \hat{w}_k^\circ - \hat{e}_k^{\circ,\dagger}$$

the KF estimate of $w$ at the time $k$, based on the information vector $I_k$ (or equivalently, on the corresponding information vector $\tilde{I}_k$). Moreover, let

$$\Sigma_k := \mathop{\mathbb{E}}_{e_k}\{(e_k - \mathop{\mathbb{E}}_{e_k}\{e_k|\tilde{I}_k\})(e_k - \mathop{\mathbb{E}}_{e_k}\{e_k|\tilde{I}_k\})'|\tilde{I}_k\} \tag{20}$$

be the (conditional) covariance matrix[9] of $e_k$, conditioned on the information vector $\tilde{I}_k$, and

$$\Sigma_{-1} := \mathop{\mathbb{E}}_{e_0}\left\{\left(e_0 - \mathop{\mathbb{E}}_{e_0}\{e_0\}\right)\left(e_0 - \mathop{\mathbb{E}}_{e_0}\{e_0\}\right)'\right\} = \Sigma_w \tag{21}$$

the (unconditional) covariance matrix[10] of $e_0$, which is equal to the (unconditional) covariance matrix of $w$.

8. For a square matrix $M$, we denote by $|\lambda|_{\max}(M)$ its spectral radius.
9. Here, the superscript "$\circ$" is omitted, to highlight that the expression (20) (and other expressions presented later, such as (25)), holds also when $e_k$ is not evaluated in correspondence of the sequence of optimal updating functions (13).
10. Likewise in [7, Appendix E.4], one could use the symbol $\Sigma_{k|k}$ to denote the (conditional) covariance matrix $\Sigma_k$, to distinguish it from the (conditional) covariance matrix of $e_{k+1}$, conditioned on the information vector $\tilde{I}_k$, and denoted by $\Sigma_{k+1|k}$. However, in the specific case they are equal, so they are both denoted by $\Sigma_k$.

**Proposition 2 (Optimal online estimate and Stochastic Riccati Equation (SRE))** *Let Assumptions 1, 2, 3, 4, and 5 be satisfied. Then*

$$\hat{w}^\circ_{k+1} = \hat{w}^\circ_k + L_k \left( \hat{w}^\circ_k - \mathop{\mathbb{E}}_{w}\{w|I_k\} \right) = \hat{w}^\circ_k + L_k(\hat{w}^\circ_k - \hat{w}^\dagger_k) = \hat{w}^\circ_k + L_k(\hat{e}^\circ_k - \hat{e}^{\circ,\dagger}_k) \,, \tag{22}$$

*where, for* $k = -1, 0, \ldots$

$$\hat{w}^\dagger_{k+1} = \hat{w}^\dagger_k + H_{k+1}(y_{k+1} - C_{k+1}\hat{w}^\dagger_k) \,, \tag{23}$$

$$H_{k+1} := \Sigma_{k+1} C'_{k+1} (\sigma^2_\varepsilon)^{-1} \,, \tag{24}$$

*and, for* $k = 0, 1, \ldots,$

$$\Sigma_k = \Sigma_{k-1} - \Sigma_{k-1} C'_k (C_k \Sigma_{k-1} C'_k + \sigma^2_\varepsilon)^{-1} C_k \Sigma_{k-1} \,, \tag{25}$$

*with the initializations*

$$\hat{w}^\dagger_{-1} = 0 \,, \tag{26}$$

$$\hat{w}^\circ_{-1} = 0 \,, \tag{27}$$

*and*

$$L_{-1} = -\left( \bar{K}_0 + \gamma I \right)^{-1} \bar{K}_0 \,. \tag{28}$$

Equation (25) has the form of the Riccati equation of the well-known *Kalman Filter* (*KF*, for short). Due to the stochastic nature of $C_k$, it can be called a "*Stochastic Riccati Equation*" (*SRE*, for short). From a computational point of view, solving (25) is easy even in a high-dimensional setting, i.e., when the dimension $d$ of the input space is large. Indeed, $C_k \Sigma_{k-1} C'_k + \sigma^2_\varepsilon$ (which needs to be inverted in (25)) is a scalar. Similarly, in formula (24) one has to invert the scalar $\sigma^2_\varepsilon$. In other applications of the Kalman filter, instead, one has to invert matrices.

**Remark 7** It is worth mentioning that also [7, Section 4.1] investigates an LQ optimal control problem with random matrices. In that case, however, there is perfect information on the state, and the randomness is limited to the system dynamics. For that problem, a suitable average Riccati equation is obtained therein, but no stochastic Riccati equation. Hence, that formulation, though inspiring for the present work, cannot be applied directly to our online-learning framework.

**Remark 8** Equations (13) and (23) show that the classical *separation principle of control and estimation* holds also for Problem $\mathrm{OLL}^N_\gamma$. More precisely, it is reduced to *two subproblems*, which can be solved independently: the *determination of the matrices* $L_k$ (solution of the *LQR subproblem*) and the *determination of the Kalman gain matrices* $H_k$ (solution of the *Kalman-filter estimation subproblem*). One might wonder why in Problem $\mathrm{OLL}^N_\gamma$ one gets, instead of the classical Riccati Equation, two different kinds of equations for the two subproblems, i.e., the ARE (17) and the SRE (25), in spite of the well-known duality between the LQR subproblem and the Kalman-filter estimation problem [45, Section 11.3].

The reason is that, when moving from the LQR subproblem to the Kalman-filter estimation subproblem, the roles of the matrices

$$A_k := I, B_k := I, Q_k, R_k := \gamma I$$

in the primal problem (i.e., the LQR subproblem) are played, respectively, by the following matrices of the dual problem (i.e., the Kalman-filter estimation problem):

$$A_k^{\text{dual}} := A_k' = I, B_k^{\text{dual}} := C_k', Q_k^{\text{dual}} := 0, R_k^{\text{dual}} := \sigma_\varepsilon^2\,,$$

where $Q_k^{\text{dual}}$ is the covariance matrix of the system noise (a kind of noise that is not present in the model (7)), hence it is an all 0's matrix. Now, in the primal problem, the matrix $Q_k$ is stochastic, whereas in the dual problem, the matrix $Q_k^{\text{dual}}$ is deterministic. Similarly, in the primal problem, the matrix $B_k$ is deterministic, whereas in the dual problem, the matrix $B_k^{\text{dual}}$ is stochastic. This lack of symmetry is the reason why the two Riccati equations (17) and (25) have different forms.

The next proposition states some properties of the solution to the SRE. For two symmetric square matrices $S_1$ and $S_2$ of the same dimension, $S_1 \preceq S_2$ means that $S_2 - S_1$ is symmetric and positive-semidefinite. When it is evident from the context, we use the symbol 0 to denote a matrix whose elements are all equal to 0.

**Proposition 3 (Properties of the solution to the SRE)** *Let Assumptions 1, 2, 3, 4, and 5 be satisfied. Then*

*(i)*

$$0 \preceq \Sigma_{k+1} \preceq \Sigma_k \tag{29}$$

*(i.e., the sequence is "non-negative" and monotonic "nonincreasing" in a generalized sense, according to $\preceq$), for all the realizations of the random matrices $\Sigma_{k+1}$ and $\Sigma_k$.*

*(ii) For all the realizations of these random matrices,*

$$0 \leq \text{Tr}\{\Sigma_{k+1}\} \leq \text{Tr}\{\Sigma_k\} \tag{30}$$

*and*

$$0 \leq \text{Tr}\{\Sigma_{k+1}^2\} \leq \text{Tr}\{\Sigma_k^2\}\,. \tag{31}$$

*(iii) There exists a symmetric and positive-semidefinite matrix $\overline{\Sigma}$ such that*

$$\lim_{k\to+\infty} \mathop{\mathbb{E}}_{\Sigma_k}\{\Sigma_k\} = \overline{\Sigma}\,. \tag{32}$$

*(iv) If*

$$\mathop{\mathbb{E}}_{Q_k}\{Q_k\} = \overline{Q} \tag{33}$$

*for all $k$ (e.g., if all the input examples $x_k$ have a common probability distribution with bounded support, and the same positive-definite covariance matrix $\overline{Q}$), then with a-priori probability 1 one has*

$$\lim_{k\to+\infty} \mathop{\mathbb{E}}_{\Sigma_k}\{\Sigma_k\} = \overline{\Sigma} = 0\,. \tag{34}$$

*When (34) holds, then*

$$\lim_{k\to+\infty} \mathrm{Tr}\left\{\mathop{\mathbb{E}}_{\Sigma_k}\{\Sigma_k\}\right\} = \mathrm{Tr}\left\{\overline{\Sigma}\right\} = 0\,. \tag{35}$$

*(v) For every* $k = -1, 0, 1, 2, \ldots$, *and all the realizations of the random matrices,*

$$\mathrm{Tr}\{F_{k+1}\Sigma_{k+1}\} \le \mathrm{Tr}\{F_{k+1}\Sigma_k\} \le \ldots \le \mathrm{Tr}\{F_{k+1}\Sigma_{-1}\}\,. \tag{36}$$

An intuitive explanation of the second bound in (30) is the following: when the time index moves from $k$ to $k+1$, the new information acquired at the time $k+1$ cannot deteriorate, on the average, the quality of the KF estimate, which is in accordante with its optimality properties [7, Appendix E]. Equations (29), (30), and (36) will be used, together with (34) and (35), in the convergence analysis of the proposed method for $k \to +\infty$ (see Section 4).

**Remark 9** An important assumption that is needed in the proof of Proposition 3 (iv) is that the common covariance matrix $\overline{Q}$ of the input examples is positive-definite. When this is not the case, this means that, with probability 1, all the input examples belong to a finite-dimensional subspace $\mathcal{S}$ of $\mathbb{R}^d$, hence, with probability 1, it is not possible to extract from the input-output pairs $(x_k, y_k)$ any information about the component of $w$ that it is orthogonal to that subspace, unless such a component is correlated with the projection of $w$ on $\mathcal{S}$. However, one still has the convergence of both the KF estimate and the OLL estimate of $w$ to the projection of $w$ on $\mathcal{S}$, as it can be shown by setting the problem directly on $\mathcal{S}$. Morover, the possible absence of information in the data about the component of $w$ that it is orthogonal to $\mathcal{S}$ has no negative consequences on the estimation process, in the sense that, in order to compute $w'x$ for a possibly unseen input $x$, one needs, with probability 1, to know only the component of $w$ that belongs to the subspace $\mathcal{S}$.

### 3.3 Role of the regularization parameter

Let us investigate the behavior of the optimal updating functions provided by (13) and (14) for the two limit cases $\gamma \simeq 0$ and $\gamma \to +\infty$, and for intermediate values of $\gamma$.

**The case** $\gamma \simeq 0$. The penalization of the update $u_k$ in the learning functional (9) becomes negligible, and one obtains $L_k \simeq -I$ from (14), and

$$u_k^\circ \simeq -\mathop{\mathbb{E}}_{e_k^\circ}\left\{e_k^\circ | \tilde{I}_k\right\} \tag{37}$$

from (13). Hence, one gets (from (7) and (37))

$$e_{k+1}^\circ \simeq e_k^\circ - \mathop{\mathbb{E}}_{e_k^\circ}\left\{e_k^\circ | \tilde{I}_k\right\}\,.$$

Equivalently, in terms of the unknown vector $w$ and its optimal estimates $\hat{w}_k^\circ$, $\hat{w}_{k+1}^\circ$ at the times $k$ and $k+1$, respectively, one obtains

$$(\hat{w}_{k+1}^\circ - w) \simeq (\hat{w}_k^\circ - w) - \left(\hat{w}_k^\circ - \mathop{\mathbb{E}}_{w}\{w | I_k\}\right)\,,$$

hence

$$\hat{w}^{\circ}_{k+1} \simeq \mathbb{E}_{w}\left\{w \middle| I_k\right\} ,$$

which is just the KF estimate of $w$ at the time $k$, based on the information vector $I_k$.

**The case** $\gamma \to +\infty$ The penalization of the update $u_k$ in the learning functional (9) becomes larger and larger. Indeed, for $\gamma$ large enough, one obtains $L_k \simeq 0$ from (14), and

$$u^{\circ}_k \simeq 0$$

from (13). Hence, one gets

$$e^{\circ}_{k+1} \simeq e^{\circ}_k$$

and

$$\hat{w}^{\circ}_{k+1} \simeq \hat{w}^{\circ}_k \simeq \ldots \simeq \hat{w}^{\circ}_0 = 0 .$$

**Intermediate values of** $\gamma$. In this case the estimate $\hat{w}^{\circ}_k$ enjoyes convergence properties similar to the ones of the KF estimate $\hat{w}^{\dagger}_k$, as illustrated numerically in Figure 1. Moreover, $\hat{w}^{\circ}_k$ is a smoothed version of the estimate $\hat{w}^{\dagger}_k$. The sequence of estimates $\hat{w}^{\circ}_k$ is smoother and less sensitive to outliers than the sequence of estimates $\hat{w}^{\dagger}_k$, as a large change in the estimate when moving from $\hat{w}^{\circ}_k$ to $\hat{w}^{\circ}_{k+1}$ is penalized by the presence of the term $\gamma u'_k u_k$ in the cost functional (9). This can be seen also by formula (22), as (14) implies that all the eigenvalues of the symmetric matrix $L_k$ are inside the unit circle. A deeper investigation of these two issues (convergence and smoothness) is made in Section 5 and Subsection 6.3, respectively.

## 4. LQG learning over an infinite horizon

To address the infinite-horizon case, we remove the final-stage cost $e'_N Q_N e_N$ (or equivalently, we assume $x_N = 0$ with probability 1, hence also $Q_N = 0$ with probability 1), and let $N \to +\infty$ (the precise formulation is provided later in this section).

**Assumption 6 (Identical distributions of the input examples)** *The random variables* $\{x_k\}$ *are identically distributed and have the same positive-definite covariance matrix, i.e.,*

$$\overline{Q}_k := \mathbb{E}_{x_k}\{x_k x'_k\} = \overline{Q}$$

*for every* $k = 0, 1, \ldots$. *Moreover, the common probability distribution has bounded support.*

Due to Assumption 6, the analysis has some similarities with the one of the optimal solution to the LQG problem performed, e.g., in [7, Section 5.2 and Appendix E.4]. We denote by $\overline{Q}^{1/2}$ a symmetric and positive-definite square root of $\overline{Q}$. As one can check directly from the definitions of reachability and observability[11] [3, Chapter 5], we observe that the

11. Given a discrete-time and time-invariant linear dynamical system of the form

$$\begin{cases} z_{t+1} = A z_t + B v_t , \\ \xi_t = C z_t + D v_t , \end{cases} \tag{38}$$

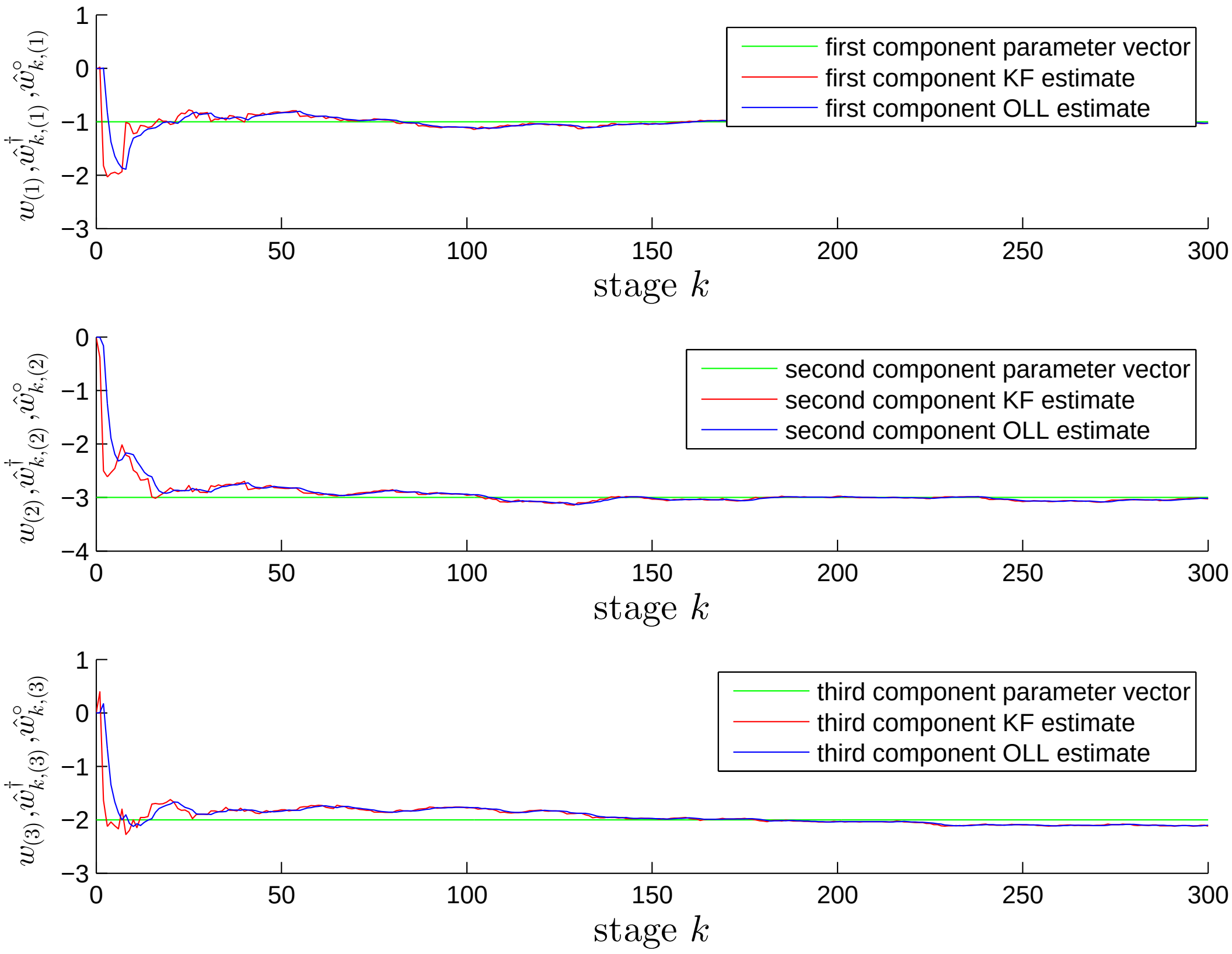

Figure 1: A comparison between the components of the OLL estimate $\hat{w}_k^\circ$ and of the KF estimate $\hat{w}_k^\dagger$. A three-dimensional case has been considered, with the realization $w = (-1, -3, -2)'$, and $N + 1 = 301$ online examples $(x_k, y_k)$ have been used to train the learning machine. The input examples have been generated with components mutually independent and uniformly distributed in $[-1, 1]$, whereas the covariance matrix $\Sigma_w$ of $w$ has been chosen to be diagonal with diagonal entries equal to 4, and the variance $\sigma_\varepsilon^2$ of the measurement noise is equal to 1, likewise the regularization parameter $\gamma$.

pair

$$(A, B) := (I, I)$$

is reachable[12], whereas the pair

$$(A, C) := (I, \overline{Q}^{1/2})$$

is observable. Hence, one can apply [7, Section 4.1, Proposition 4.1], from which it follows that the ARE (17) admits a stationary solution $\overline{K}$, associated with the two stationary matrices

$$L := -(\overline{K} + \gamma I)^{-1}\overline{K} \tag{39}$$

and

$$F := \overline{K}(\overline{K} + \gamma I)^{-1}\overline{K}$$

(see (16)). Moreover, by reversing the time-indices in (17) (i.e., setting $t := N - k$ and $P_t := \overline{K}_{N-k}$), the solution $P_{t+1}$ of the ARE

$$P_{t+1} = P_t - P_t(P_t + \gamma I)^{-1}P_t + \overline{Q} \tag{40}$$

(which is equivalent to (17)) converges to $\overline{K}$ for any initialization of the positive-semidefinite matrix $P_0$, still by [7, Section 4.1, Proposition 4.1].

The (average) learning functional over infinite horizon is defined as follows.

**Definition 2 (Average Learning functional over infinite horizon)** *Let $\gamma > 0$, and*

$$\begin{aligned} J_\gamma^\infty\left(\left\{u_k(\tilde{I}_k)\right\}_{k=0}^{\infty}\right) &:= \liminf_{N\to+\infty}\left(\frac{1}{N}\underset{w,\{x_k\}_{k=0}^{N-1},\{\varepsilon_k\}_{k=0}^{N-2}}{\mathbb{E}}\left\{\sum_{k=0}^{N-1}\left[e_k' Q_k e_k + \gamma u_k' u_k\right]\right\}\right) \\ &= \liminf_{N\to+\infty}\left(\frac{1}{N}\underset{w,\{x_k\}_{k=0}^{N-1},\{\varepsilon_k\}_{k=0}^{N-2}}{\mathbb{E}}\left\{\sum_{k=0}^{N-1}\left[(e_k' x_k)^2 + \gamma u_k' u_k\right]\right\}\right). \end{aligned} \tag{41}$$

**Problem OLL$_\gamma^\infty$ (On-Line Learning over infinite horizon)**. *Given the examples $(x_k, y_k)$ generated at each time instant $k = 0, 1, \ldots,$ according to the model defined by Assumptions 1 and 2, and the learning machine defined by Assumption 3, find the infinite sequence $u_0^\circ(I_0), u_1^\circ(I_1), \ldots,$ of* optimal updating functions *with the structure defined by Assumption 4, that minimizes the average learning functional (41).*

Likewise for Problem OLL$_\gamma^N$, for every $k = 0, 1, \ldots,$ we shall call $\hat{w}_k$ *online estimate* (*OLL estimate*, for short). In the following, we consider directly the LQG case, with identical distributions for the input examples.

---

where $z_t \in \mathbb{R}^n$, $v_t \in \mathbb{R}^m$, $\xi_t \in \mathbb{R}^p$, the pair $(A, B)$ is reachable if and only if, starting from any initial state, any other state can be reached at some subsequent finite time $\overline{t}$, by choosing an appropriate sequence of $\overline{t}$ controls. Moreover, the pair $(A, C)$ is said to be observable if and only if, given any sequence of measures $\xi_0, \ldots, \xi_{t-1}$ and applied controls $v_0, \ldots, v_{t-1}$ for $\hat{t} \geq 1$ sufficiently long, it is possible to determine exactly the initial state $z_0 \in \mathbb{R}^n$ of the dynamical system (38).

12. Actually, the expression used in [7, Section 4.1, Definition 1.1] for this situation is "controllable pair", but the definition provided therein is actually the one of "reachable pair" reported in footnote 11.

**Proposition 4 (Optimal updating functions and ARE)** *Let Assumptions 1, 2, 3, 4, 5, and 6 be satisfied. Then, updating functions that solve Problem* $OLL_\gamma^\infty$ *are given, for* $k = 0, 1, \ldots,$ *by*

$$u_k^\circ(\tilde{I}_k) = L \mathop{\mathbb{E}}_{e_k^\circ} \left\{ e_k^\circ \middle| \tilde{I}_k \right\} , \tag{42}$$

*where* $L$ *is defined in (39)*[13].

**Remark 10** A nice feature of the AREs (17) (with $\mathop{\mathbb{E}}_{Q_k}\{Q_k\} = \overline{Q}$) and (40) is that their common stationary solution $\overline{K}$ can be easily expressed in terms of the eigenvalues/eigenvectors of the matrix $\overline{Q}$, which can be useful from a computational point of view. Indeed, let us express $\overline{Q}$ as

$$\overline{Q} = U \Lambda_{\overline{Q}} U' ,$$

where $U$ is a basis of orthogonal unit-norm eigenvectors of $\overline{Q}$ (hence, $U' = U^{-1}$), and $\Lambda$ is a diagonal matrix collecting the corresponding positive eigenvalues. We recall that $\overline{K}$ satisfies

$$\overline{K} = \overline{K} - \overline{K}(\overline{K} + \gamma I)^{-1} \overline{K} + \overline{Q} ,$$

i.e,

$$\overline{K}(\overline{K} + \gamma I)^{-1} \overline{K} = \overline{Q} . \tag{43}$$

Now, $\overline{K}(\overline{K} + \gamma I)^{-1} \overline{K}$ has the same eigenvectors as $\overline{K}$. Hence, also $\overline{K}$ and $\overline{Q}$ have the same eigenvectors, so $\overline{K}$ is expressed as

$$\overline{K} = U \Lambda_{\overline{K}} U' , \tag{44}$$

where $\Lambda_{\overline{K}}$ is a suitable diagonal matrix, with positive eigenvalues. Due to (43), the diagonal elements $(\Lambda_{\overline{K}})_{(i,i)}$ and $(\Lambda_{\overline{Q}})_{(i,i)}$ are related by

$$(\Lambda_{\overline{K}})^2_{(i,i)} ((\Lambda_{\overline{K}})_{(i,i)} + \gamma)^{-1} = (\Lambda_{\overline{Q}})_{(i,i)} .$$

Hence, by the positiveness of $(\Lambda_{\overline{K}})_{(i,i)}$, one gets

$$(\Lambda_{\overline{K}})_{(i,i)} = \frac{(\Lambda_{\overline{Q}})_{(i,i)} + \sqrt{(\Lambda_{\overline{Q}})^2_{(i,i)} + 4\gamma (\Lambda_{\overline{Q}})_{(i,i)}}}{2} .$$

Similarly, the stationary matrix $L$ can be expressed as

$$L = U \Lambda_L U' , \tag{45}$$

where the elements $(\Lambda_L)_{(i,i)}$ of the diagonal matrix $\Lambda_L$ are

$$(\Lambda_L)_{(i,i)} = - \frac{(\Lambda_{\overline{K}})_{(i,i)}}{(\Lambda_{\overline{K}})_{(i,i)} + \gamma} = - \frac{(\Lambda_{\overline{Q}})_{(i,i)} + \sqrt{(\Lambda_{\overline{Q}})^2_{(i,i)} + 4\gamma (\Lambda_{\overline{Q}})_{(i,i)}}}{(\Lambda_{\overline{Q}})_{(i,i)} + \sqrt{(\Lambda_{\overline{Q}})^2_{(i,i)} + 4\gamma (\Lambda_{\overline{Q}})_{(i,i)}} + 2\gamma} . \tag{46}$$

13. Here, we recall that $\mathop{\mathbb{E}}_{e_k^\circ} \left\{ e_k^\circ \middle| \tilde{I}_k \right\}$ is the KF estimate of the error vector $e_k^\circ$ at the time $k$, based on the information vector $\tilde{I}_k$.

A particularly simple case occurs when the matrix $\overline{Q}$ is diagonal, hence one can choose the matrices $U$ and $U'$ as the identity matrix $I$, so also $\overline{K}$ and $L$ are diagonal, too, by formulas (44) and (45). Moreover, if $\overline{Q}$ is proportional to the identity matrix $I$, also $\overline{K}$ and $L$ are proportional to $I$. Finally, similar remarks hold also in the finite-horizon case for the matrices $\overline{K}_k$ and $L_k$, in case the matrices $\mathop{\mathbb{E}}_{Q_k}\{Q_k\}$ commute (this happens, e.g., when all the matrices $\mathop{\mathbb{E}}_{Q_k}\{Q_k\}$ are equal to $\overline{Q}$).

## 5. Convergence properties of the On-Line Learning estimates in terms of mean-square errors

Let us denote by

$$\mathrm{MSE}_k^{\dagger} := \mathop{\mathbb{E}}_{w,\hat{w}_k^{\dagger}} \left\{ \left(w - \hat{w}_k^{\dagger}\right)' \left(w - \hat{w}_k^{\dagger}\right) \right\}$$

the mean-square error of the KF estimate at time $k$. Differently from the LQG case detailed in [7], the expectation of $\Sigma_k$ is needed here, as $\Sigma_k$ depends on the sequence of random matrices $C_0, \dots, C_k$. Under the assumptions of Proposition 3 (iv), this converges to the 0 matrix as $k$ tends to $+\infty$ (see formula (34)).

**Proposition 5 (Convergence of the MSE of the KF estimate)** *Let Assumptions 1, 2, 3, 4, 5, and 6 be satisfied. Then the following hold.*

*(i)*

$$\mathrm{MSE}_k^{\dagger} = \mathrm{Tr}\left\{ \mathop{\mathbb{E}}_{\Sigma_k}\{\Sigma_k\} \right\}.$$

*(ii) For every $k = 1, 2, \dots$,*

$$\mathrm{Tr}\left\{ \mathop{\mathbb{E}}_{\Sigma_k}\{\Sigma_k\} \right\} \leq \sqrt{\frac{(c_1 + \sigma_\varepsilon^2)\, d\, \mathrm{Tr}\{\mathop{\mathbb{E}}_{\Sigma_0}\{\Sigma_0\}\}}{k \lambda_{\min}(\overline{Q})}}\,, \tag{47}$$

*where $c_1$ is a positive constant such that $C_{k+1}\Sigma_{-1}C'_{k+1} \leq c_1$ with a-priori probability 1. Moreover, $\lim_{k\to+\infty} \mathrm{MSE}_k^{\dagger} = 0$.*

*(iii)*

$$\lim_{k\to+\infty} \mathop{\mathbb{E}}_{H_k}\{H_k\} = 0\,. \tag{48}$$

*(iv) Every element $H_{k,(h,l)}$ of $H_k$ converges to 0 also in probability, i.e., for every $\delta > 0$,*

$$\lim_{k\to+\infty} \Pr\{|H_{k,(h,l)}| > \delta\} = 0\,. \tag{49}$$

Note that the upper bound in Proposition 5 (ii) provides for the convergence to 0 of the mean-square error of the KF estimate of $w$ at the time $k$, a rate of order $O(\sqrt{1/k})$. As to Proposition 5 (iv), an intuitive explanation is the following: as the parameter $w$ to be learned does not change in time, after a sufficiently large number of "good" examples, the learning machine has practically learned $w$, and future examples are practically not needed (of course, this holds in the case - considered so far - in which the parameter $w$ does not change with time; see Section 8 for a relaxation of this assumption).

Now, let

$$\mathrm{MSE}_k^{\circ} := \mathop{\mathbb{E}}_{w,\hat{w}_k^{\circ}} \left\{ (w - \hat{w}_k^{\circ})' \, (w - \hat{w}_k^{\circ}) \right\}$$

denote the mean-square error of the OLL estimate at time $k$. The next proposition provides a recursion to compute and bound from above $\mathrm{MSE}_k^{\circ}$, and states its convergence to 0. We refer to Remark 18 in the Appendix for a possible way to derive estimates of the associated rate of convergence.

Let

$$e_k^{\dagger} := \hat{w}_k^{\dagger} - w \,,$$

and denote by

$$\begin{aligned} \Sigma_{e_k^{\dagger}} &:= \mathop{\mathbb{E}}_{e_k^{\dagger}} \left\{ \left( e_k^{\dagger} - \mathop{\mathbb{E}}_{e_k^{\dagger}}\{e_k^{\dagger}\} \right) \left( e_k^{\dagger} - \mathop{\mathbb{E}}_{e_k^{\dagger}}\{e_k^{\dagger}\} \right)' \right\} \\ &= \mathop{\mathbb{E}}_{e_k^{\dagger}} \left\{ \left( e_k^{\dagger} \right) \left( e_k^{\dagger} \right)' \right\} \end{aligned} \tag{50}$$

the (unconditional) covariance matrix of $e_k^{\dagger}$. Moreover, we denote by

$$\begin{aligned} \Sigma_{e_k^{\circ}} &:= \mathop{\mathbb{E}}_{e_k^{\circ}} \left\{ \left( e_k^{\circ} - \mathop{\mathbb{E}}_{e_k^{\circ}}\{e_k^{\circ}\} \right) \left( e_k^{\circ} \mathop{\mathbb{E}}_{e_k^{\circ}}\{e_k^{\circ}\} \right)' \right\} \\ &= \mathop{\mathbb{E}}_{e_k^{\circ}} \left\{ (e_k^{\circ}) \, (e_k^{\circ})' \right\} \end{aligned} \tag{51}$$

the (unconditional) covariance matrix of $e_k^{\circ}$.

**Proposition 6 (Convergence of the MSE of the OLL estimate)** *Let Assumptions 1, 2, 3, 4, 5, and 6 be satisfied. Then the following hold.*

*(i)*

$$\mathrm{MSE}_k^{\circ} = \mathrm{Tr} \left\{ \Sigma_{e_k^{\circ}} \right\} \,,$$

*(ii)*

$$\begin{aligned} \mathrm{MSE}_k^{\circ} \;\; &\leq \;\; \mathrm{Tr} \left\{ (I + L_{k-1}) \Sigma_{e_{k-1}^{\circ}} (I + L_{k-1})' \right\} + \mathrm{Tr} \left\{ L_{k-1} \Sigma_{e_{k-1}^{\dagger}} L_{k-1}' \right\} \\ &+ \;\; 2 \sqrt{ \mathrm{Tr} \left\{ (I + L_{k-1}) \Sigma_{e_{k-1}^{\circ}} (I + L_{k-1})' \right\} \mathrm{Tr} \left\{ L_{k-1} \Sigma_{e_{k-1}^{\dagger}} L_{k-1}' \right\} } \,. \end{aligned}$$

*(iii) Under the assumptions made in Subsection 4, one has*

$$\lim_{k \to +\infty \mathrm{MSE}_k^{\circ}} = \lim_{k \to +\infty} \mathop{\mathbb{E}}_{w,\hat{w}_k^{\circ}} \left\{ (w - \hat{w}_k^{\circ})' \, (w - \hat{w}_k^{\circ}) \right\} = 0 \,. \tag{52}$$

## 6. Comparisons with other machine-learning techniques

In this section, some connections and comparisons are presented between the solutions to our learning paradigm, and machine-learning techniques such as average regret minimization (Subsection 6.1), stochastic gradient descent (Subsection 6.2), and Kalman-based estimates (Subsections 6.3 and 6.4).

### 6.1 Connections with average regret minimization

Likewise for the finite-horizon case, the minimization of the average learning functional (41) is equivalent to the minimization of the alternative learning functional

$$\liminf_{N\to+\infty}\left(\frac{1}{N}\mathop{\mathbb{E}}_{w,\{x_k\}_{k=0}^{N-1},\{\varepsilon_k\}_{k=0}^{N-2}}\left\{\sum_{k=0}^{N-1}\left[(\hat{w}_k'x_k-y_k)^2+\gamma u_k'u_k\right]\right\}\right)\,, \tag{53}$$

since (53) is just equal to

$$J_\gamma^\infty\left(\left\{u_k(\tilde{I}_k)\right\}_{k=0}^\infty\right)+\sigma_\varepsilon^2 \tag{54}$$

(see formula (6)). We now consider the limit case $\gamma=0$, denoting by $J_0^\infty$ the corresponding average learning functional. Then, observing that the following equality holds (due to the independence of the disturbance noises $\varepsilon_k$):

$$\sigma_\varepsilon^2=\liminf_{N\to+\infty}\left(\frac{1}{N}\mathop{\mathbb{E}}_{w,\{x_k\}_{k=0}^{N-1},\{\varepsilon_k\}_{k=0}^{N-2}}\left\{\sum_{k=0}^{N-1}\left[(w'x_k-y_k)^2\right]\right\}\right)\,, \tag{55}$$

we obtain

$$\begin{aligned}J_0^\infty\left(\left\{u_k(\tilde{I}_k)\right\}_{k=0}^\infty\right) &= \liminf_{N\to+\infty}\left(\frac{1}{N}\mathop{\mathbb{E}}_{w,\{x_k\}_{k=0}^{N-1},\{\varepsilon_k\}_{k=0}^{N-2}}\left\{\sum_{k=0}^{N-1}\left[(\hat{w}_k'x_k-y_k)^2\right]\right\}\right)\\ &\quad-\liminf_{N\to+\infty}\left(\frac{1}{N}\mathop{\mathbb{E}}_{w,\{x_k\}_{k=0}^{N-1},\{\varepsilon_k\}_{k=0}^{N-2}}\left\{\sum_{k=0}^{N-1}\left[(w'x_k-y_k)^2\right]\right\}\right)\,,\end{aligned} \tag{56}$$

which can be interpreted as an *average regret functional* [52]. Hence, under the assumptions of Proposition 4 (with $\gamma=0$), the sequence of KF estimates minimizes the average regret functional (56). Moreover, its minimum is 0 because, by Proposition 5 (ii), one has

$$\lim_{k\to+\infty}\mathrm{MSE}_k^\dagger=\lim_{k\to+\infty}\mathop{\mathbb{E}}_{w,\hat{w}_k^\dagger}\left\{\left(w-\hat{w}_k^\dagger\right)'\left(w-\hat{w}_k^\dagger\right)\right\}=\lim_{k\to+\infty}\mathop{\mathbb{E}}_{e_k^\dagger}\left\{\left(e_k^\dagger\right)'\left(e_k^\dagger\right)\right\}=0\,, \tag{57}$$

then one gets also

$$\begin{aligned}\lim_{k\to+\infty}\mathop{\mathbb{E}}_{e_k^\dagger,Q_k}\left\{\left(e_k^\dagger\right)'Q_k\left(e_k^\dagger\right)\right\} &= \lim_{k\to+\infty}\mathop{\mathbb{E}}_{e_k^\dagger}\left\{\left(e_k^\dagger\right)'\bar{Q}\left(e_k^\dagger\right)\right\}\\ &\leq \lambda_{\max}(\bar{Q})\lim_{k\to+\infty}\mathop{\mathbb{E}}_{e_k^\dagger}\left\{\left(e_k^\dagger\right)'\left(e_k^\dagger\right)\right\}\\ &= 0\,.\end{aligned} \tag{58}$$

Nevertheless, the next proposition shows that, for any $\gamma>0$, also the sequence of OLL estimates minimizes the average regret functional (56).

**Proposition 7** *For any $\gamma > 0$, under the assumptions of Proposition 4, the sequence of OLL estimates minimizes the average regret functional (56).*

### 6.2 Connections with KF and stochastic gradient descent

At each time $k+1$, our OLL estimate $\hat{w}^{\circ}_{k+1}$ of the parameter vector $w$ associated with the data-generation model has the following recursive form (see, e.g., the statement of Proposition 2):

$$\hat{w}^{\circ}_{k+1} = \hat{w}^{\circ}_{k} + L_k(\hat{w}^{\circ}_{k} - \hat{w}^{\dagger}_{k}) \,, \tag{59}$$

where $L_k$ is a suitable square matrix and $\hat{w}^{\dagger}_{k} := \mathbb{E}_{w}\{w|I_k\}$ is the Kalman-filter estimate of $w$ at time $k$, based on the vector $I_k$ that collects all the information available to the learning machine up to time $k$. Hence, it follows from (59) that our estimates are obtained from the Kalman-filter estimates through an additional smoothing step. The form of equation (59) is similar to the one of other online estimates obtained through various machine learning techniques, such as stochastic gradient descent [46, Chapter 3]. However, there is a substantial difference: we derive (59) as the optimal solution of a suitable optimal control/estimation problem, showing various interesting consequences of that in the paper, made possible by our use of optimal control/estimation techniques. We believe that this approach could be fruitfully applied also to other machine learning techniques used in online learning. Moreover, we offer a principled way to construct the matrix $L_k$ in (59) as the solution of a suitable Riccati equation.

### 6.3 Outperformance with respect to KF in terms of smoothness

For simplicity, we consider the finite-horizon case. The extension to the infinite-horizon case can be performed by a limit process, likewise in Section 4. The next proposition shows that the OLL estimates are smoother than the KF estimates, in the sense that the value of

$$\mathbb{E}_{\{u_k\}_{k=0}^{N-1}} \left\{ \sum_{k=0}^{N-1} u_k' u_k \right\}$$

when $\gamma > 0$ and the updates are generated by (13) is smaller than or equal to the corresponding value obtained when the updates are generated by (37) with "$\simeq$" replaced by "$=$" (i.e., in the limit $\gamma \to 0$). The limit problem obtained when $\gamma$ tends to 0 is just the Kalman-estimation problem, whose optimal sequence of updates is

$$u_k^{\dagger} := -\mathbb{E}_{e_k^{\dagger}} \left\{ e_k^{\dagger} \middle| \tilde{I}_k \right\} \,,$$

where

$$e_k^{\dagger} := w_k^{\dagger} - w_k \,.$$

**Proposition 8** *Let Assumptions 1, 2, 3, 4, and 5 be satisfied. Then*

$$\mathbb{E}_{\{u_k^{\circ}\}_{k=0}^{N-1}} \left\{ \sum_{k=0}^{N-1} \left[ (u_k^{\circ})'(u_k^{\circ}) \right] \right\} \leq \mathbb{E}_{\{u_k^{\dagger}\}_{k=0}^{N-1}} \left\{ \sum_{k=0}^{N-1} \left[ (u_k^{\dagger})'(u_k^{\dagger}) \right] \right\} .$$

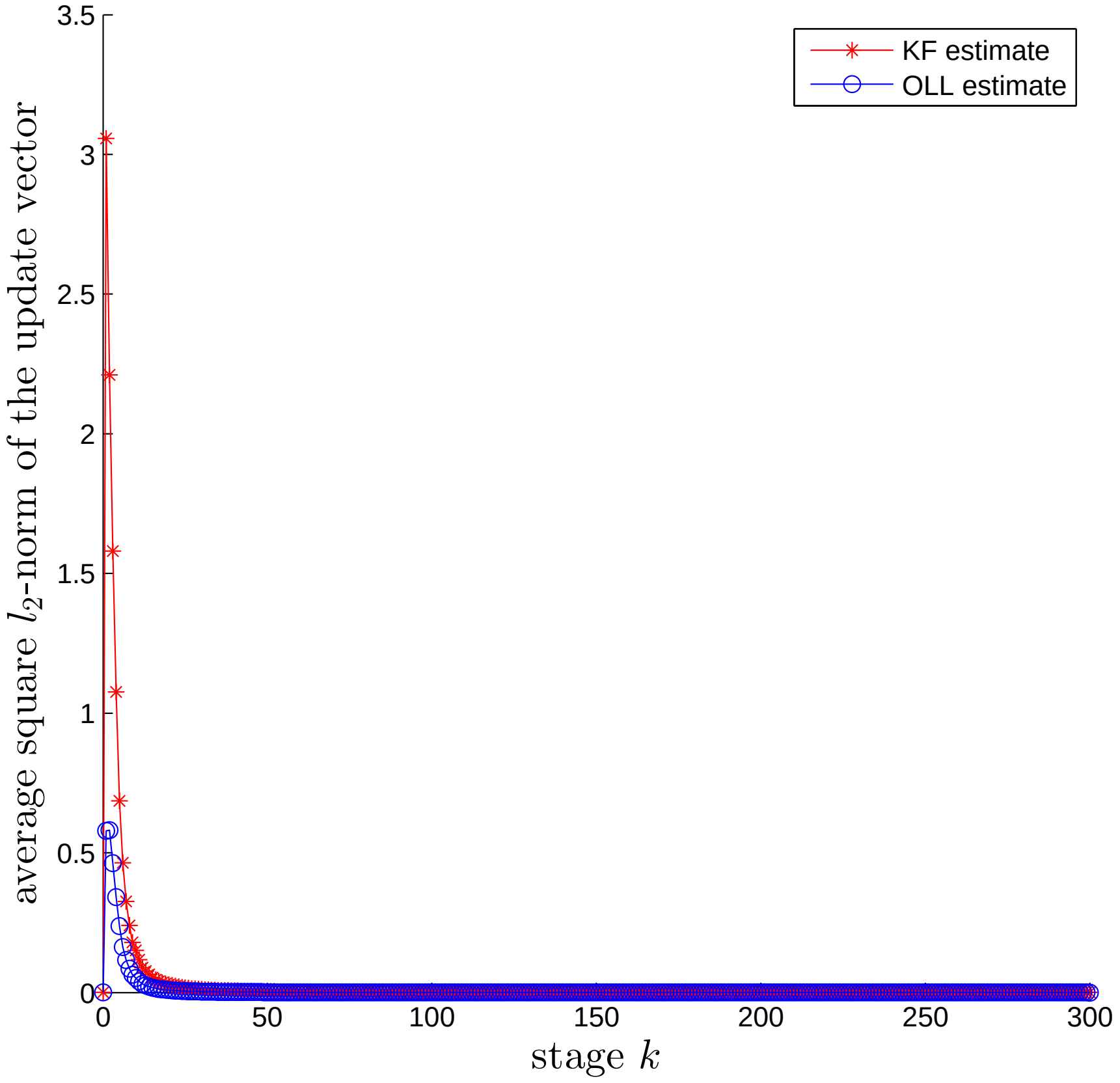

Figure 2: A comparison between the empirical averages of the square $l_2$-norms of the vectors of updates used to generate the OLL estimate $\hat{w}_k^\circ$ and the KF estimate $\hat{w}_k^\dagger$. The parameters are the same as in Figure 1, apart from $w$, which is generated according to a Gaussian distribution, with mean $(0, 0, 0)'$ and covariance matrix $\Sigma_w = 4I$. The empirical averages of the square $l_2$-norms have been computed by considering 10000 independent simulations.

The simulation results shown in Figure 2, which refers to a setup similar to the one of Figure 1, are in line with the result from Proposition 8. The figure suggests that, for every $k$, the stronger result

$$\mathop{\mathbb{E}}_{u_k^\circ} \left\{ \left[ (u_k^\circ)'(u_k^\circ) \right] \right\} \leq \mathop{\mathbb{E}}_{u_k^\dagger} \left\{ \left[ (u_k^\dagger)'(u_k^\dagger) \right] \right\} , \tag{60}$$

may also hold.

Finally, Figure 3 shows that both approaches are suitable also for parameter vectors of much larger dimension. Indeed, it refers to the case of $d = 100$, and $N+1 = 1000+1$ online examples. The figure reports the square $l_2$-norm of the error vector associated with the KF estimate and with the OLL estimate, respectively, at the generic stage $k$. The running time of such a simulation (whose code was written in MATLAB R2013, likewise for all the other simulations) was of about 28 seconds, on a notebook with a 1.40 GHz CPU and 4 GB of RAM. The figure also shows that, for this case of a time-invariant parameter vector, in general a smaller error is associated to the KF estimate with respect to the OLL estimate. However, the KF estimate is less smooth with respect to the time index $k$. In item *e.3)* of Section 8, it is shown that, for the case of a slowly time-varying parameter vector, the OLL estimate can achieve even a smaller error than the KF estimate, under a suitable periodic re-initialization of the matrices $\Sigma_k$ (see Figure 7 in Section 8).

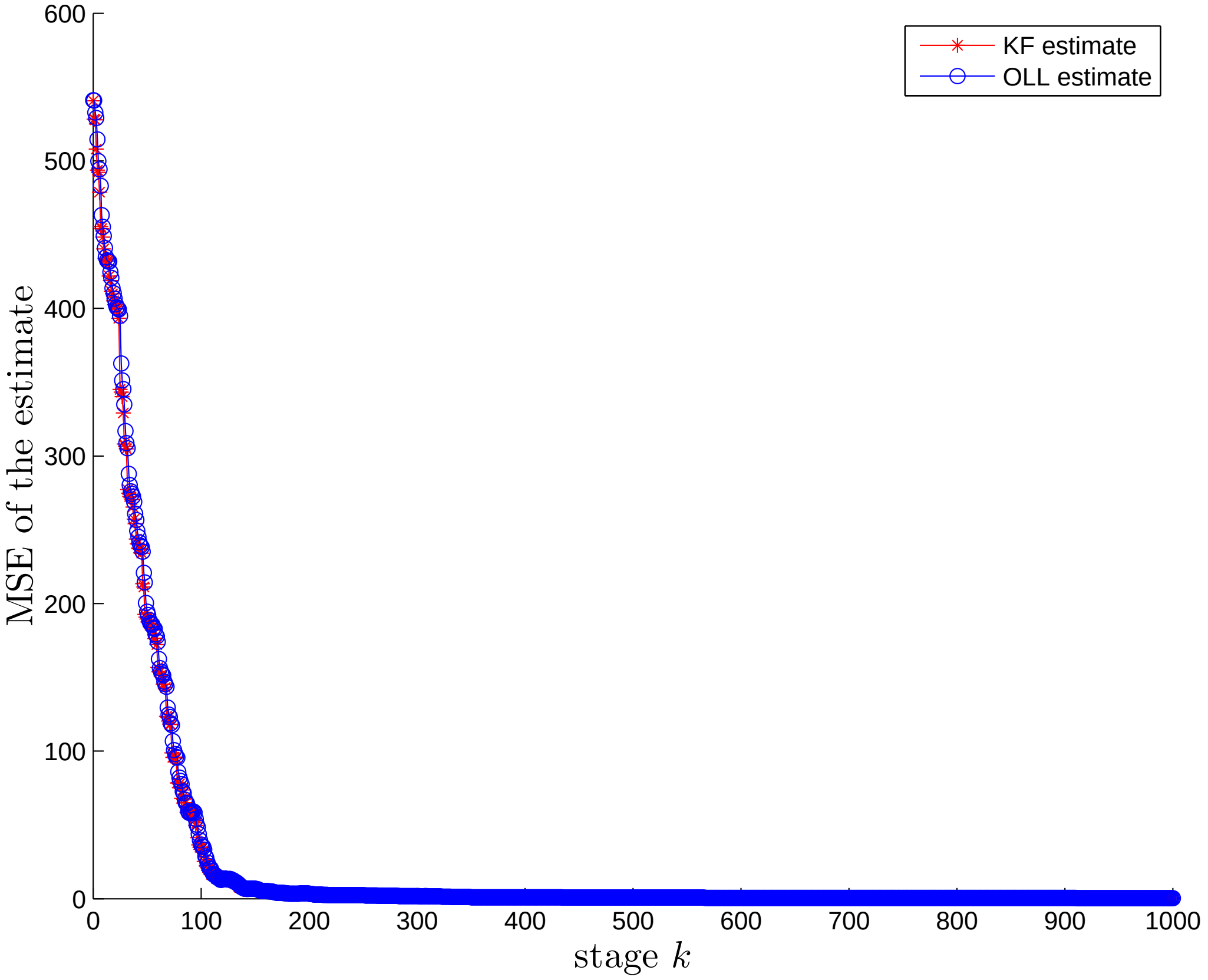

Figure 3: For a setup similar to the one of Figure 2, but with $d = 100$ and $N + 1 = 1000 + 1$ online examples: square $l_2$-norm of the error vector associated with the KF estimate and OLL estimate at the generic stage $k$.

### 6.4 Outperformance with respect to KF in terms of sensitivity to outliers

Here we further compare numerically the KF and OLL estimates, now in terms of their different sensitivity to outliers. To this end, we alter periodically the output data-perturbation model, choosing the disturbance $\varepsilon_k$ to be equal to a positive constant $z_1$ when $k$ is a multiple of some positive integer $Z$, otherwise equal to a negative constant $z_2$ when $k$ is not a multiple of $Z$. The two constants $z_1$ and $z_2$ are chosen in such a way that the empirical mean (over any time-window of duration $Z$) of the $\varepsilon_k$'s is 0 (i.e., the condition $z_1 + (Z-1)z_2 = 0$ is imposed), and the empirical variance (over the same time-window) is $\sigma_\varepsilon^2$ (i.e., the condition $\frac{z_1^2+(Z-1)z_2^2}{Z-1} = \sigma_\varepsilon^2$ is imposed). Hence, $z_1 = \frac{(Z-1)\sigma_\varepsilon}{\sqrt{Z}}$ and $z_2 = -\frac{\sigma_\varepsilon}{\sqrt{Z}}$ are obtained. Moreover, for a fair comparison with the KF estimate, this modified assumption on the output data-perturbation model is not included in the optimization problem producing the OLL estimates[14]. In other words, that knowledge is not provided to the learning machine.

The numerical results reported in Figure 4 show clearly the much smaller sensitivity to outliers of the OLL estimates with respect to the KF ones (details about the parameter choices are reported in the caption of the figure). This is ultimately due to the larger smoothness of the OLL estimates with respect to the time index.

An additional theoretical motivation for the smaller sensitivity to outliers of our OLL estimates is obtained by an inspection of formula (22) in Proposition 2. Limiting for simplicity of the analysis to the first OLL updates of the parameter vector, it follows by that formula and by $\hat{w}_{-1}^\circ = 0$ that $\hat{w}_0^\circ = 0$ and $\hat{w}_1^\circ = -L_0 \hat{w}_0^\dagger$, where $\hat{w}_0^\dagger$ is the first KF update, which is influenced only by the first example presented to the learning machine. Since $|\lambda|_{\max}(L_0) < 1$ (see formula (14)), it is evident that the OLL estimate is less influenced by the presence of a possible outlier. Moreover, such an influence decreases by increasing the regularization parameter $\gamma$, since, by formula (14), the larger $\gamma$, the smaller $|\lambda|_{\max}(L_0)$. Figure 5 confirms this advantage of the OLL estimates, showing that the $l_2$-norms of the differences between consecutive OLL estimates are typically much smaller than the $l_2$-norms of the differences between consecutive KF estimates. An additional significant advantage of the OLL estimates in the presence of time-varying parameter vectors is detailed in item *e)* of Section 8.

## 7. Nonlinear models of data-generation and application of kernel methods

An interesting extension of the model investigated in the sections above is obtained by mapping the input data $x_k$ preliminarily to another Euclidean space, then applying the model in the new input space. More precisely, one introduces a (possibly nonlinear) mapping $\phi : \mathbb{R}^d \to E$, where $E$ is an Euclidean space of dimension $d_E$, possibly larger than $d$ (or even infinite). Then, the measurement equation (1) becomes

$$y_k = w' \phi(x_k) + \varepsilon_k \, , \tag{61}$$

where the parameter vector $w$ belongs now to $E$. In this case, one can still apply all the techniques described in the paper taking $E$ as the new input space. Of course, in doing this, the dimensions of some matrices would in general increase: for instance, in case of a finite $d_E$, the matrices $\overline{K}$, $L$, and $\Sigma_k$ would become $d_E \times d_E$ matrices, whereas the Kalman gain matrix $H_k$ would become a $1 \times d_E$ matrix. In case of an infinite-dimensional Hilbert space $E$, they would be replaced by suitable infinite-dimensional linear operators.

Interestingly, as we show below, when doing such an extension, one can apply the so-called "kernel trick" of kernel machines [15]. More precisely, we show some circumstances under which, for every (possibly unseen) input $x$, one can express both $(\hat{w}_k^\dagger)' \phi(x)$ and

14. Nevertheless, it is worth observing that the sequence of measures generated in this way is still an admissible sequence of measures for the original Gaussian disturbance model.

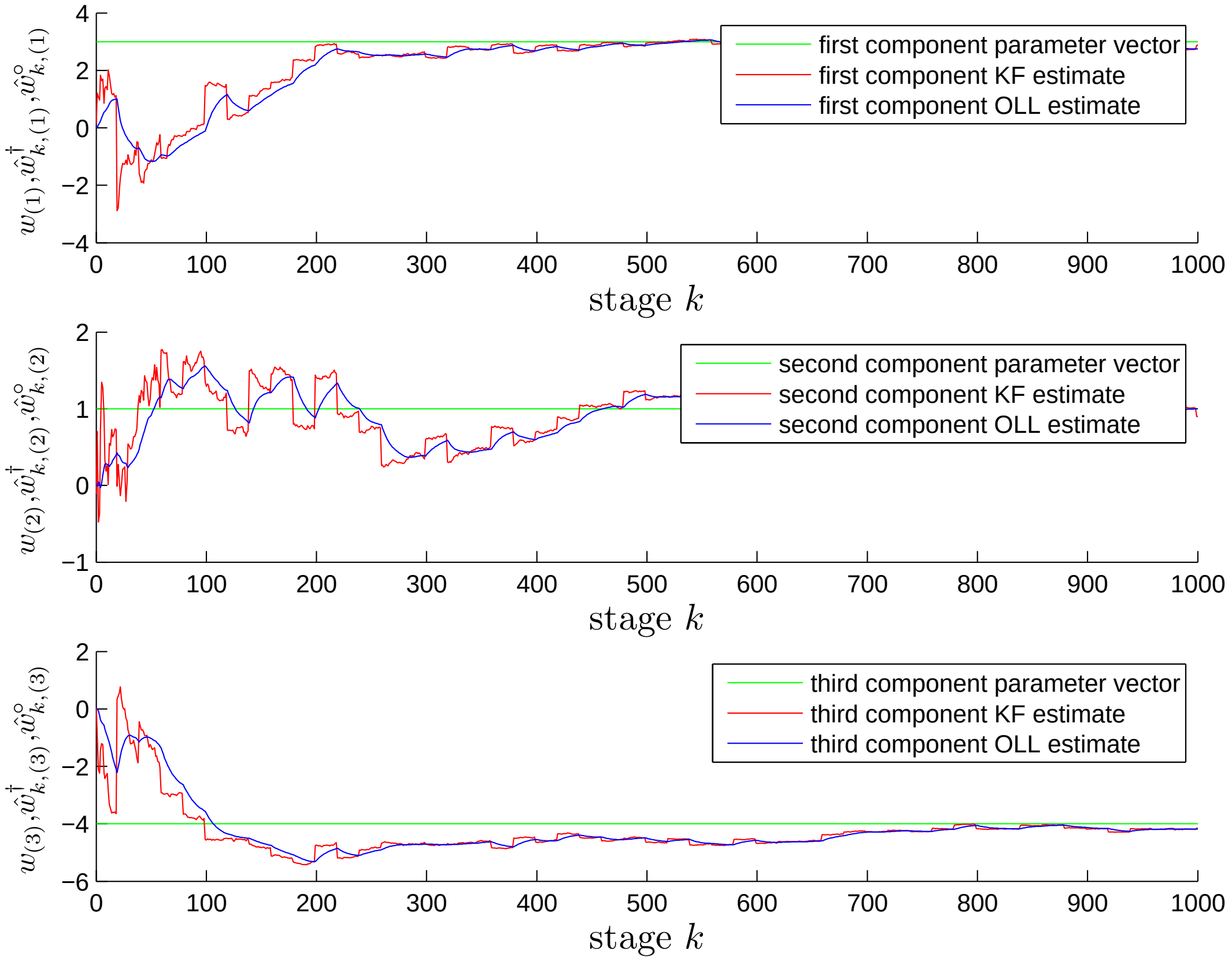

Figure 4: For a setup similar to the one of Figure 1, but choosing $\sigma_\varepsilon = 10$, $\gamma = 50$, $N + 1 = 1001$ examples, the diagonal entries of $\Sigma_w$ equal to 10, and the disturbance $\varepsilon_k$ equal to $z_1 = \frac{(Z-1)\sigma_\varepsilon}{\sqrt{Z}}$ when $k$ is a multiple of $Z = 20$, otherwise equal to $z_2 = -\frac{\sigma_\varepsilon}{\sqrt{Z}}$: comparison between the components of the OLL estimate $\hat{w}^{\circ}_k$ and of the KF estimate $\hat{w}^{\dagger}_k$.

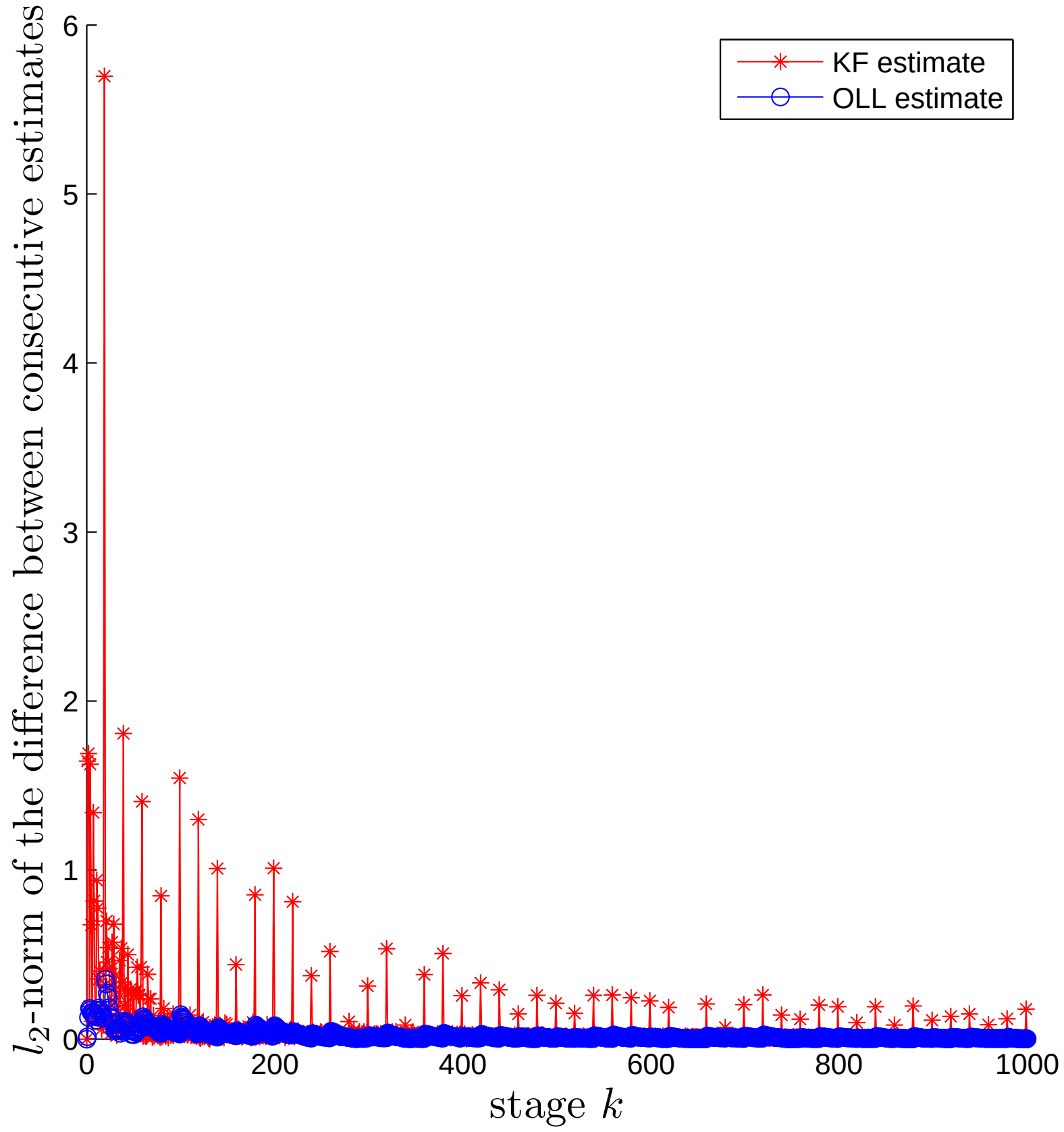

Figure 5: For the example in Figure 4: comparison between the $l_2$-norms of the differences between consecutive OLL estimates and the $l_2$-norms of the differences between consecutive KF estimates.

$(\hat{w}^\circ)'\phi(x)$ in terms of inner products of the form $\phi(x_j)'\phi(x)$, where $x_j$ is an input example already seen by the learning machine. Hence, if one is able to express $\phi(x_j)'\phi(x)$ in a simple way (e.g., through a symmetric kernel function $\mathcal{K} : \mathbb{R}^d \times \mathbb{R}^d \to \mathbb{R}$ such that $\phi(x_j)'\phi(x) = \mathcal{K}(x_j, x)$), one can compute $(\hat{w}_k^\dagger)'\phi(x)$ and $(\hat{w}^\circ)'\phi(x)$ even without knowing explicitly the expression of the mapping $\phi$.

**Remark 11** As an example of a mapping $\phi$ and its associated kernel $\mathcal{K}$, we consider the case $d = 2$ and the feature mapping $\phi : \mathbb{R}^2 \to E = \mathbb{R}^6$, defined as

$$\phi(x) := \left(1\,,\sqrt{2}\,x_{(1)}\,,\sqrt{2}\,x_{(2)}\,,\sqrt{2}\,x_{(1)}x_{(2)}\,,x_{(1)}^2\,,x_{(2)}^2\right)'\,,$$

where $x_{(1)}$ and $x_{(2)}$ are the two components of the vector $x$. Then, given any two input vectors $x, z \in \mathbb{R}^2$, the inner product $\phi(x)'\phi(z)$ is expressed as

$$\begin{aligned}\phi(x)'\phi(z) &= 1 + 2x_{(1)}z_{(1)} + 2x_{(2)}z_{(2)} + 2x_{(1)}x_{(2)}z_{(1)}z_{(2)} + x_{(1)}^2 z_{(1)}^2 + x_{(2)}^2 z_{(2)}^2 \\ &= (1 + x'z)^2 \\ &:= \mathcal{K}(x, z)\,,\end{aligned}$$

which is the so-called homogeneous polynomial kernel [15, Section 3.2] of order 2, whose evaluation involves only computations to be performed in the original input space $\mathbb{R}^2$.

We first consider how to compute $(\hat{w}_k^\dagger)'\phi(x)$ using kernels, then how to compute $(\hat{w}_k^\circ)'\phi(x)$, too. We make the following assumption.

**Assumption 7 (Covariance matrix of the measurement noise)** *Let*

$$\Sigma_w = \nu I_{d_E}\,, \tag{62}$$

*where $\nu > 0$ and $I_{d_E}$ denotes the (matrix associated with the) identity operator on $E$.*

The results presented in the next proposition for the kernel version of the KF estimate are essentially the same as the ones obtained in [33, Theorems 2 and 3], which shows also how to express the linear combinations inside such equations, through an application of the matrix inversion lemma (see, e.g., [41, Section 2.6]). However, their extension to the kernel version of the OLL estimate, provided in the next Proposition 10, is novel. In order to improve their readability, in the next formulas (63), (64), (65), (69), (70), (71), (72) only the functional form of the right-hand side is provided.

**Proposition 9** *Let Assumptions 1, 2, 3, 4, 5, and 7 be satisfied for the kernel version of the KF estimate (ie., with every $x_k$ replaced by $\phi(x_k)$). Then, for every $k = 0, 1, \ldots,$*

$$H_k = \Sigma_k (C_k^{(\phi)})'(\sigma_\varepsilon^2)^{-1} = \text{linear combination of } \phi(x_0), \ldots, \phi(x_k)\,, \tag{63}$$

$$\hat{w}_k^\dagger = \text{linear combination of } \phi(x_0), \ldots, \phi(x_k)\,, \tag{64}$$

*and*

$$(\hat{w}_k^\dagger)'\phi(x) = \text{linear combination of } \mathcal{K}(x_0, x), \ldots, \mathcal{K}(x_k, x)\,. \tag{65}$$

In case of a finite-dimensional space $E$, the convergence analysis is exactly the same as the one in Proposition 5, and a similar (even though more technical) analysis is expected to hold for the infinite-dimensional case. Finally, in case the (matrix associated with the) covariance operator

$$\overline{Q}^{(\phi)} := \mathbb{E}_{\phi(x)}\{\phi(x)\phi(x)'\} \tag{66}$$

is only positive-semidefinite but not positive-definite, one could still follow Remark 9 to prove the convergence of the estimate on the subspace on which the input data lie with probability 1. Such a subspace could be estimated, e.g., by an application of Kernel Principal Component Analysis (KPCA) [42]. Moreover, one could even redefine the problem taking that subspace as the new input space, making the operator $\overline{Q}^{(\phi)}$ be positive-definite when restricted on it.

After dealing with the kernel-version of the KF estimate of $w$, we now investigate the kernel-version of its OLL estimate. We make the following assumption.

**Assumption 8 (Covariance operator)** *Let one of the following hold.*

*(i) The covariance operator $\overline{Q}^{(\phi)}$ has the form*

$$\overline{Q}^{(\phi)} = qI_{d_E}\,, \tag{67}$$

*for some $q > 0$*

*(ii)*

$$\overline{Q}^{(\phi)} = \overline{Q}^{(\phi)}_{\text{emp}} := \frac{1}{l_U}\sum_{j=1}^{l_U} \phi(\tilde{x}_j)\phi(\tilde{x}_j)'\,, \tag{68}$$

*where $l_U$ is a given positive integer, and $\{\tilde{x}_j, j = 1, \ldots, l_U\}$ are some unsupervised examples (assumed here for simplicity to be available to the learning machine starting from the time $k = 0$).*

**Remark 12** Assumption 8 (i) refes to a particularly simple model for $\overline{Q}^{(\phi)}$, which is relaxed in Assumption 8 (ii), which refers to the case in which $\overline{Q}^{(\phi)}$ is modeled by an empirical estimate $\overline{Q}^{(\phi)}_{\text{emp}}$ obtained using the unsupervised examples.

**Proposition 10** *(i) Let Assumptions 1, 2, 3, 4, 5, 7, and 8 (i) be satisfied for the kernel version of the KF estimate (ie., with every $x_k$ replaced by $\phi(x_k)$). Then, for every $k = 0, 1, \ldots,$*

$$\hat{w}_k^\circ \quad = \quad \text{linear combination of } \phi(x_0), \ldots, \phi(x_{k-1}) \tag{69}$$

*and*

$$(\hat{w}_k^\circ)'\phi(x) \quad = \quad \text{linear combination of } \mathcal{K}(x_0, x), \ldots, \mathcal{K}(x_{k-1}, x)\,. \tag{70}$$

*(ii) If, instead, Assumption 8 (ii) is used, then, for every $k = 0, 1, \ldots,$*

$$\hat{w}_k^\circ \quad = \quad \text{linear combination of } \phi(\tilde{x}_1), \ldots, \phi(\tilde{x}_{l_U}) \tag{71}$$

*and*

$$(\hat{w}_k^\circ)'\phi(x) \quad = \quad \text{linear combination of } \mathcal{K}(\tilde{x}_1, x), \ldots, \mathcal{K}(\tilde{x}_{l_U}, x)\,. \tag{72}$$

**Remark 13** A significant advantage of the representation (70) over the ones (64) and (69) is that the vector $\hat{w}_k^\circ$ has dimension at most $l_U$.

**Remark 14** More generally, $\tilde{x}_1, \ldots, \tilde{x}_{l_U}$ in Assumption 8 (ii) could be previously seen input data, preferably not used by the learning machine in combination with labels, to reduce/avoid overtraining. So, their number could grow up as the learning machine acquires examples. Of course, after adding new empirical data in the estimate (68), one could also update accordingly the matrix $L_k$ (or, in the infinite-horizon case, the stationary matrix $L$), likewise in item *d.1)* of the next Section 8.

We conclude this section by mentioning that, in the nonlinear case, differently from techniques such as the extended KF [30], the kernel version of the OLL estimate has the advantage of solving an optimal control (or optimal estimation) problem. Other approaches to online learning with kernels are described, e.g., in the review paper [19].

## 8. Extensions

In this section, we illustrate some other extensions of the proposed OLL estimation scheme investigated in the paper.

*a)* **Nonzero mean of** $x_k$: in this case, no significant change in the analysis is required. The only difference is that $\mathop{\mathbb{E}}_{Q_k}\{Q_k\}$ and $\overline{Q}$ are now correlation matrices, instead than covariance matrices.

*b)* **Nonzero mean of** $w$: Propositions 5 and 6 still hold true if $\mathop{\mathbb{E}}_{w}\{w\} \neq 0$, and the KF estimate and the OLL estimate are initialized, respectively, by

$$\hat{w}_{-1}^\dagger = \mathop{\mathbb{E}}_{w}\{w\} \,,$$

and

$$\hat{w}_0^\circ = \hat{w}_0 = \mathop{\mathbb{E}}_{w}\{w\} \tag{73}$$

(notice that two different initialization indices have been used for the two estimates, where the subscript "$-1$" has been used to denote the "a-priori" KF estimate, i.e., the one obtained before the presentation of the first example, whereas the subscript "0" has been used for the initialization of the OLL estimate[15]). Indeed, in such a case one obtains similar expressions as in the Appendix for the matrices $\Sigma_k$ in (136), for the matrices $\Sigma_{e_k^\dagger}$, $\Sigma_{e_k^\circ}$, $\Sigma_{e_k^\circ, e_k^\dagger}$, $\Sigma_{e_k^\dagger, e_k^\circ}$ in (50), (51), (153) and (154), respectively, and the same equation (161), which is used therein to obtain the convergence result (52) through an analysis of the convergence of $\mathrm{Tr}\{\Sigma_{e_k^\circ}\}$ when $k$ tends to $+\infty$.

15. Recall that $\hat{w}_0^\circ$ refers to the OLL estimate obtained before seeing the first example, $\hat{w}_1^\circ$ refers to the OLL estimate obtained after seeing the first example but before seeing the second one, and so on. Instead, $\hat{w}_{-1}^\dagger$ refers to the KF estimate obtained before seeing the first example, whereas $\hat{w}_0^\dagger$ refers to the KF estimate obtained after seeing the first example but before seeing the second one, and so on. Hence, according to the current notation, there is a shift in the indices of the two estimates, the available information beeing the same. Of course, a more uniform notation could have been used, instead, at the expense of shifting and renaming the index for the KF estimate, but using a less common notation for it.

**Remark 15** The case $\mathop{\mathbb{E}}_{w}\{w\} \neq 0$ is important in practice, and - among other ones - it models the situation in which, after some number $\overline{k}$ of measures, the time index $k$ is shifted to the left (i.e., $k$ is replaced by $k - \overline{k}$, or equivalently, one reformulates Problem OLL using $\overline{k}$ instead of 0 as the initial index in the summation of its objective (5)), and the knowledge derived by the previous estimates (i.e., the one up to the time $\overline{k} - 1$) is used to generate the term $\mathop{\mathbb{E}}_{w}\{w\}$ (this is actually an "a-posteriori" knowledge, since it summarizes the knowledge deriving from the previous estimates, but becomes the new "a-priori" knowledge for the new problem with modified starting index in the summation).

Similar convergence results are obtained even if one replaces (73) with

$$\hat{w}_0^\circ = \hat{w}_0 \neq \mathop{\mathbb{E}}_{w}\{w\} \,. \tag{74}$$

Indeed, in such a case, the convergence analysis of $\mathrm{Tr}\{\Sigma_{e_k^\circ}\}$ made in the Appendix is still valid. The only difference is that now one has $\mathop{\mathbb{E}}_{e_k^\circ}\{e_k^\circ\} \neq 0$, but $\mathop{\mathbb{E}}_{e_k^\circ}\{e_k^\circ\}$ tends also to 0 exponentially fast as $k$ tends to $+\infty$, due to equation (152) in the Appendix with $L_j$ replaced by $L$, since the matrix $I + L$ has spectral radius smaller than 1.

*c)* **Introduction of a bias in the model**: instead of the measurement equation (1), one could consider the one

$$y_k = w' x_k + \varepsilon_k + b \,, \tag{75}$$

where $b$ is an additional parameter to be learned using the sequence of examples available online. This case can be reduced to (1) by replacing the input vector $x_k$ by $(x_k', 1)'$, and the parameter vector $w$ to be learned by $(w', 1)'$. As the last component of the new input vector $(x_k', 1)'$ has nonzero mean, one is also reduced to the case *a)* above. Moreover, the assumption of positive-definiteness of the correlation matrix of $(x_k', 1)'$ is satisfied automatically if it is satisfied by the covariance matrix of $x_k$.

*d)* **More complex models for the measurement errors**: the measurement errors $\varepsilon_k$ could be have nonzero means, nonidentical distributions, and/or be not mutually independent. The first two cases can be dealt with in a straightforward way: indeed, in the first case one has only to subtract the expectation of $\varepsilon_k$ from the measure $y_k$ before presenting it as an input to the KF[16], while in the second case one has to insert an additional index $k$ to $\sigma_\varepsilon^2$, using terms of the form $\sigma_{\varepsilon_k}^2$ in the Kalman-filter recursion scheme (25) and in the Kalman gain matrix (24). Finally, in the correlated case one could model the measurement noise as the output of an auxiliary uncontrolled linear dynamical system, which receives mutually independent noises as inputs. In this case, when the horizon tends to $+\infty$, the convergence of the solution of the ARE to a stationary solution could be more difficult to prove (or such a convergence could even not hold at all), since the reachability condition

16. I.e., by re-defining $\tilde{y}_{k+1}^\circ$ in equation (107) as

$$\tilde{y}_{k+1}^\circ = C_{k+1} \hat{w}_{k+1}^\circ - (y_{k+1} - \mathop{\mathbb{E}}_{\epsilon_{k+1}}\{\varepsilon_{k+1}\}) \,.$$

needed for the application of [7, Section 4.1, Proposition 4.1] would be violated.

*e)* **Time-varying models for some parameters**: when solving the "shifted version" of Problem OLL that uses $\overline{k}$ as the initial index (see Remark 15), one could exploit, for some of its parameters, time-varying models (which could be also estimated online), including the cases of:

*e.1)* slowly time-varying covariance matrices $\mathop{\mathbb{E}}_{Q_k}\{Q_k\}$ of the input examples $x_k$;

*e.2)* slowly time-varying variances $\sigma^2_{\varepsilon_k}$ of the measurement noises $\varepsilon_k$;

*e.3)* a slowly time-varying parameter vector in the data-generation model (1).

About the issue *e.1)* above, we notice that the covariance matrices $\mathop{\mathbb{E}}_{Q_k}\{Q_k\}$ may be estimated online from the already-observed input data[17] $x_k$ ($k = 0, \ldots, \overline{k} - 1$), when such covariance matrices do not depend on the time index (as assumed in the basic infinite-horizon version of the problem), but also if one makes the assumption, instead, that they are slowly time-varying (in such a case, one could give more weight in such estimates to the last input data using, e.g., a forgetting factor). Now, let us suppose that one is solving the infinite-horizon version of Problem OLL, under the assumptions of Section 4, replacing the "common"[18] covariance matrix $\overline{Q}$ of the future input data by its initial estimate, here denoted as $\overline{Q}(0)$. However, let us also suppose that the estimate at time $\overline{k}$ of $\overline{Q}$ derived from the previous input data (denoted in the following as $\overline{Q}(\overline{k})$) is significantly different from $\overline{Q}(0)$ (due, e.g., to slow changes with respect to time of the covariance matrices $\mathop{\mathbb{E}}_{Q_k}\{Q_k\}$, or simply due to a possibly bad initial estimate $\overline{Q}(0)$). Then, using the updated estimate $\overline{Q}(\overline{k})$ as the model for the "common" covariance matrix of the future input data, one could update also the stationary matrix $L$ of the proposed optimal controller, replacing $\overline{Q}$ by $\overline{Q}(\overline{k})$ in the average Riccati equation (17), and looking for its (new) "stationary" solution $\overline{K}$, hence deriving the (new) "stationary" matrix $L$ from (39). Such new solutions will be "stationary" as long as the future estimates of the "common" covariance matrix of the future input data will not change significantly with respect to $\overline{Q}(\overline{k})$.

Concerning the issue *e.2)* above, one may detect changes in $\sigma^2_{\varepsilon_k}$ from the already-observed input/output data $(x_k, y_k)$ ($k = 0, \ldots, \overline{k} - 1$), if one assumes that $\sigma^2_{\varepsilon_k}$ is also slowly time-varying. In particular, such changes would be easily detected if, inside that set, the learning machine has at its disposal some pairs $(x_k, y_k)$ with similar values for $x_k$, or even several groups of such pairs, each of which is characterized by similar values of $x_k$. In this way, indeed, one could generate, inside the $i$-th such group

$$G^{(i)} := \left\{ \left( x^{(i)}_{k_1}, y^{(i)}_{k_1} \right), \ldots, \left( x^{(i)}_{k_{|G^{(i)}|}}, y^{(i)}_{k_{|G^{(i)}|}} \right) \right\},$$

17. In practice, in order to reduce/avoid overtraining, one could use a subset of the input data available online up to time $\overline{k}$ to estimate the covariance matrices $\mathop{\mathbb{E}}_{Q_k}\{Q_k\}$, and include only the data associated with the other subset in the definition of the learning functional of Problem OLL. A similar remark holds for the online estimate of the variance of the measurement noise, which is described in item *e.2)*.

18. Even though, in this paragraph, such a covariance matrix is modeled as time-varying, one could assume that, on the basis of the knowledge available at time $\overline{k}$, the future input data will have a similar model as the present (time-varying) one.

the auxiliary $|G^{(i)}| \times |G^{(i)}|$ data matrix $Y^{\mathrm{aux},(i)}$, whose element $Y^{\mathrm{aux},(i)}_{(h,l)}$ is defined as

$$Y^{\mathrm{aux},(i)}_{(h,l)} := y^{(i)}_{k_h} - y^{(i)}_{k_l}\,,$$

and has a very small dependence on $w$, since

$$y^{(i)}_{k_h} - y^{(i)}_{k_l} = w'\left(x^{(i)}_{k_h} - x^{(i)}_{k_l}\right) + \varepsilon^{(i)}_{k_h} - \varepsilon^{(i)}_{k_l}\,,$$

and

$$x^{(i)}_{k_h} - x^{(i)}_{k_l} \simeq 0\,.$$

Then, one could use such matrices $G^{(i)}$ to estimate the variance of the measurement noise, giving more importance/weight to the last among such measures. The obtained estimate would be then used as the "common" variance $\sigma^2_\varepsilon$ in the "shifted version" of Problem OLL that uses $\overline{k}$ as the initial index, presented in Section 4.

Finally, about the issue *e.3)* above, we show in the following that a minor modification of the proposed OLL model can learn even a time-varying parameter vector $w$. We first focus on the case in which, at some time $\hat{k} \leq \overline{k}$, the parameter vector $w$ changes, then remains fixed. In this case, the mean-square error of the estimate of the new parameter vector still converges to 0 when $k$ tends to $+\infty$ (see formula (52), together with item *b)* in this section). However, the trace of $\Sigma_{\hat{k}}$ may be extremely small (see (47)), making also the trace of $\Sigma_k$ be extremely small, for every $k \geq \hat{k}$. Since the traces of such matrices are used to bound from above the trace of the covariance matrix of the error at time $k$ of the OLL estimate (see again item *b)* in this section), the convergence to 0 of the error with respect to the new parameter vector may be extremely small, for both the KF estimate and the OLL estimate. This issue could be solved, in both cases, by a re-initialization at time $\overline{k}$ of the covariance matrix $\Sigma_{\overline{k}}$ to the "a-priori" covariance matrix $\Sigma_w$. More generally, a periodic re-initialization could be used, to track a periodically (or continuosly) changing parameter vector $w$. In this case, the OLL estimate would have the advantage, with respect to the KF estimate, to change more slowly in time, making it more suitable to a slowly time-varying parameter vector. The next figures provide more insights about this last issue.

Figure 6 and 7 refer to a parameter that changes periodically, the change in its components being small and random. In Figure 6, there is no re-initialization to $\Sigma_w$ of any of the matrices $\Sigma_k$, and the convergence to the new parameter vector of both the KF estimate and the OLL estimate is slow. Figure 7, instead, refers to the case in which - the change in the parameter vector being the same - there is a periodic re-initialization of the matrices $\Sigma_k$ to $\Sigma_w$ (this second period has been chosen as different from the first one, just to avoid giving the learning machine the advantage of knowing when the parameter vector changes, in order to model the more realistic situation in which this knowledge is not available to the machine). In this case, both estimates are able to track the time-varying parameter vector $w$ in a better way, but the OLL estimate is smoother than the KF estimate, due to the presence of the regularization parameter $\gamma > 0$. So, in this context, the OLL estimate is preferable to the KF estimate if one knows that the parameter vector changes slowly with time. Figures 8 and 9 show the reason for which this happens: when there is no re-initialization to $\Sigma_w$ of the matrices $\Sigma_k$, the Frobenius norm[19] of the Kalman gain matrix $H_k$ is expected to be small for $k$ large (see formulas (24) and (47)), which is confirmed by Figure 8. Hence, even though the norm of the error $y_k - C_k \hat{w}_k^\dagger$ tends to increase when the parameter vector changes, the KF estimate of $w$ at time $k$ is not affected so much by this change (see formula (23)), hence also the OLL estimate does not change so much (see formula (22)). Instead, a re-initialization to $\Sigma_w$ tends to make the Frobenius norm of the Kalman gain matrix $H_k$ bigger (which is confirmed by Figure 9), and this amplifies the effect of the larger norm of $y_k - C_k \hat{w}_k^\dagger$ (due to the change in the parameter vector $w$) on the KF estimate of $w$ at time $k$. For the OLL estimate, the change in the estimate is expected to be smaller, due to the smoothing effect of formula (22). So, likewise the KF estimate, also the OLL estimate is able to track the change in the parameter vector, but with a smoother behavior.

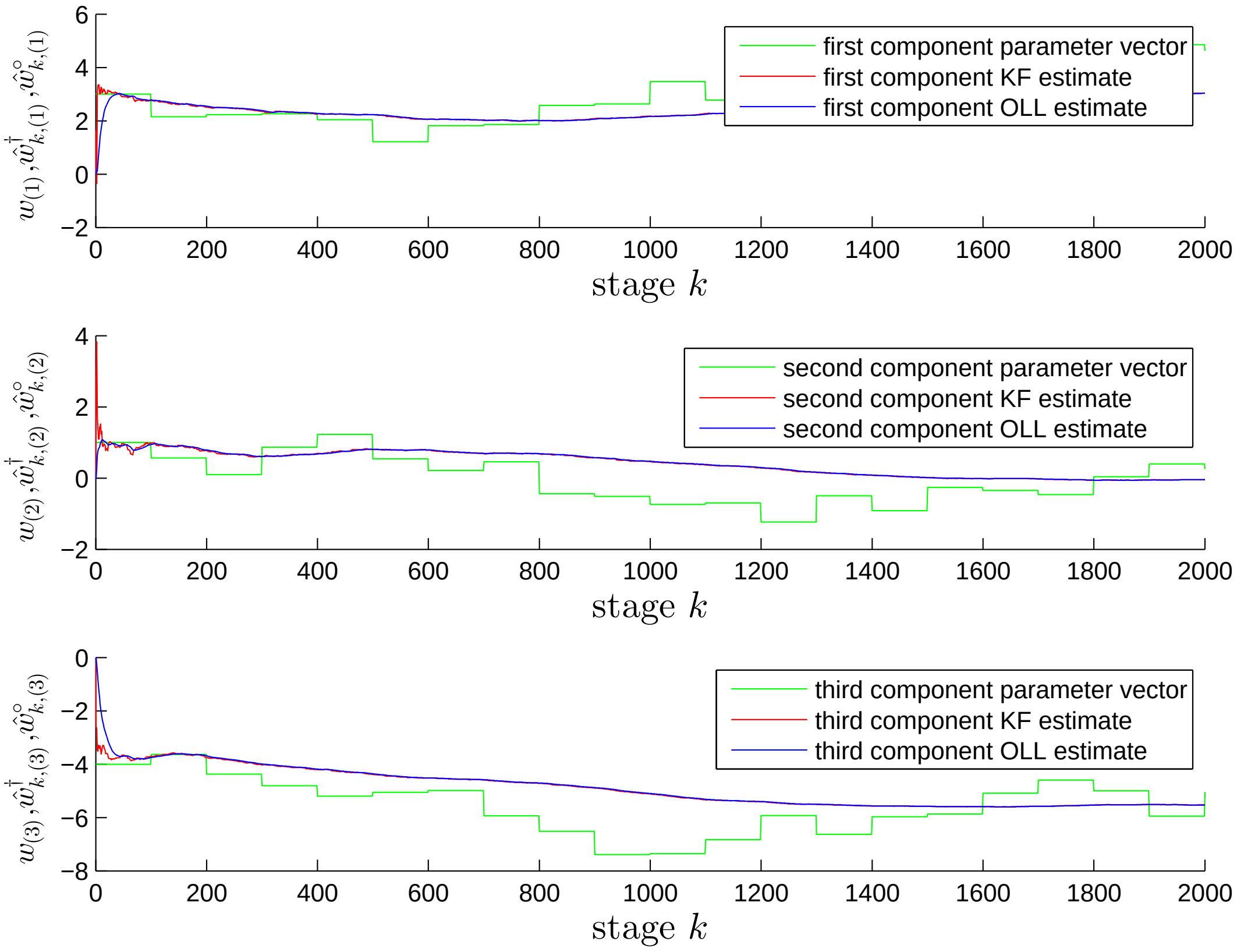

Figure 6: For the case of a time-varying parameter vector: a comparison between the components of the optimal estimate $\hat{w}_k^\circ$ at the time $k$ of the parameter vector $w$, obtained by solving Problem OLL modeling online learning, and the corresponding components of the estimate $\hat{w}_k^\dagger$ at the time $k$, obtained by applying the Kalman filter. A setting similar to the one of Figure 2 has been considered, but with $N+1 = 2000+1$ online examples, and the parameter vector $w$ randomly changed of a small amount every 100 online examples. The covariance matrix $\Sigma_w$ of the initial $w$ has been chosen to be diagonal with diagonal entries equal to 64, the variance $\sigma_\varepsilon^2$ of the measurement noise has been chosen to be equal to 1, and the regularization parameter $\gamma$ has been set to 30. No periodic re-initialization to $\Sigma_w$ of any of the matrices $\Sigma_k$ was performed.

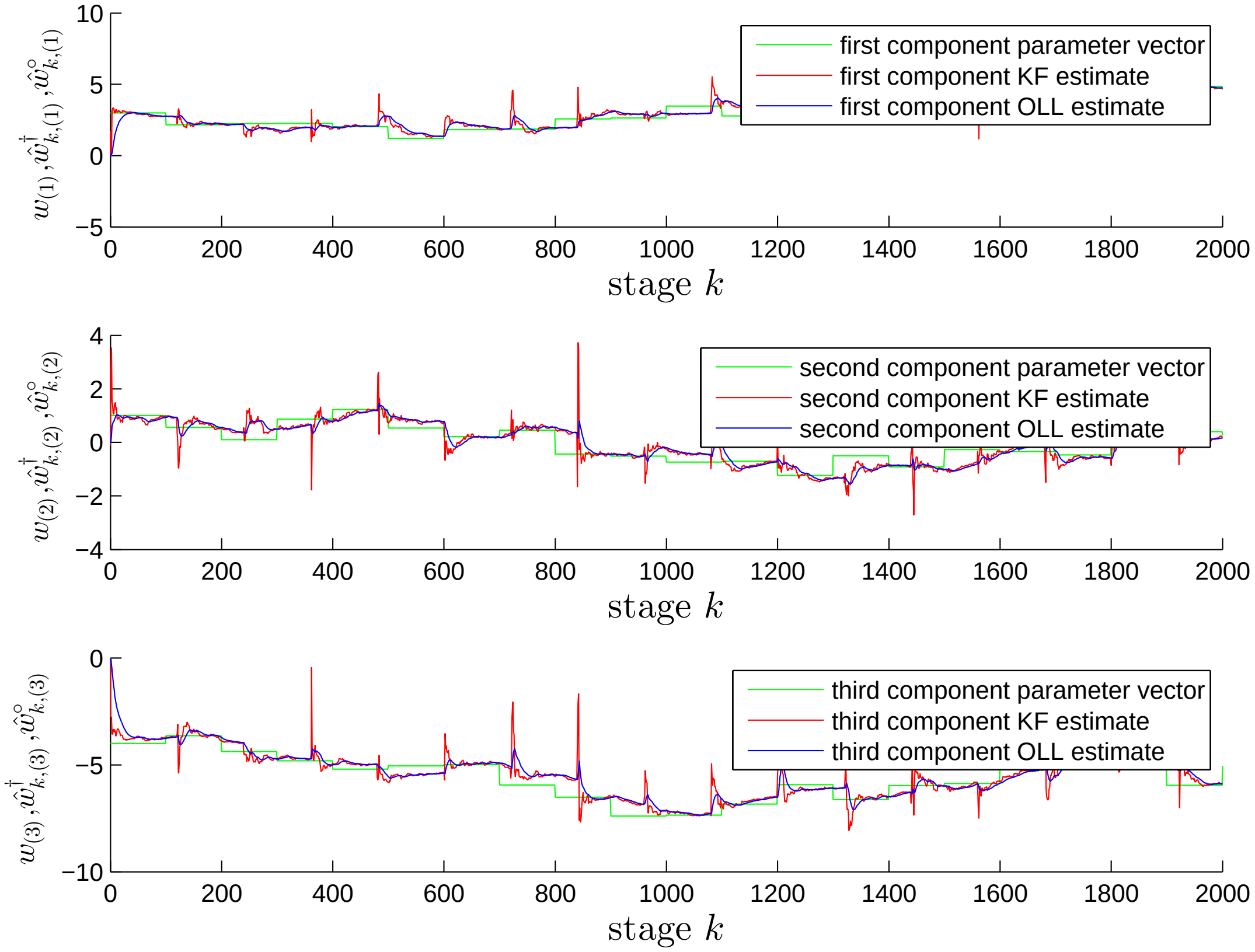

Figure 7: For the case of a time-varying parameter vector: a comparison similar to the one of Figure 6, with the same online examples and the same changes in $w$, but with a re-initialization to $\Sigma_w$ of the matrices $\Sigma_k$ performed every 120 online examples.

Figures 10 and 11 demonstrate also that a periodic re-initialization to $\Sigma_w$ of the matrices $\Sigma_k$ does not negatively affect so much the tracking of a *time-invariant* parameter vector $w$ in the case of the OLL estimate, whereas the KF estimate is more negatively affected. Again, this is due to the smoothing effect of formula (22).

Finally, we mention that various different approaches to deal with learning time-varying parameters online were presented, e.g., also in [13, 33, 35], in some cases in the context of state estimation of dynamical systems in the presence of outliers [1, 18]. As a possible extension, one could combine those approaches with the regularization of the updates included in our model, whose beneficial effects have been just demonstrated also in this time-varying case.

*f)* **Insertion of a discount factor in the problem**: one may be interested to give different weights to future expected errors, giving more importance to the present. This can be modeled by inserting a discount factor $\rho \in (0, 1)$, and modifying the learning functional (9) as follows:

$$\begin{aligned} J_{\gamma,\rho}^N\left(\left\{u_k(\tilde{I}_k)\right\}_{k=0}^{N-1}\right) &:= \mathop{\mathbb{E}}_{e_0,\{x_k\}_{k=0}^N,\{\tilde{\varepsilon}_k\}_{k=0}^{N-1}} \left\{\sum_{k=0}^{N-1} \rho^k \left[(e_k' x_k)^2 + \gamma u_k' u_k\right] + \rho^N (e_N' x_N)^2\right\} \\ &= \mathop{\mathbb{E}}_{e_0,\{x_k\}_{k=0}^N,\{\tilde{\varepsilon}_k\}_{k=0}^{N-1}} \left\{\sum_{k=0}^{N-1} \rho^k \left[e_k' Q_k e_k + \gamma u_k' u_k\right] + \rho^N e_N' Q_N e_N\right\}. \end{aligned} \tag{77}$$

Then, the resulting problem is just a variation, with random matrices, of the discounted LQ/LQG problem (see [17, Section 6.3] for the version with deterministic matrices).

In practice, the modification (77) changes only slightly the Bellman equations for the cost-to-go functions, with the introduction of the discount factor $\rho$. Particularly, the updates

19. We recall that, for a $n \times m$ real matrix $S$, its Frobenius norm is defined as

$$\|S\|_{\text{Frobenius}} := \sqrt{\sum_{i=1}^n \sum_{j=1}^m S_{(i,j)}^2}\,. \tag{76}$$

In this problem, as being $H_k$ a column vector, its Frobenius norm coincides with its $l_2$-norm.

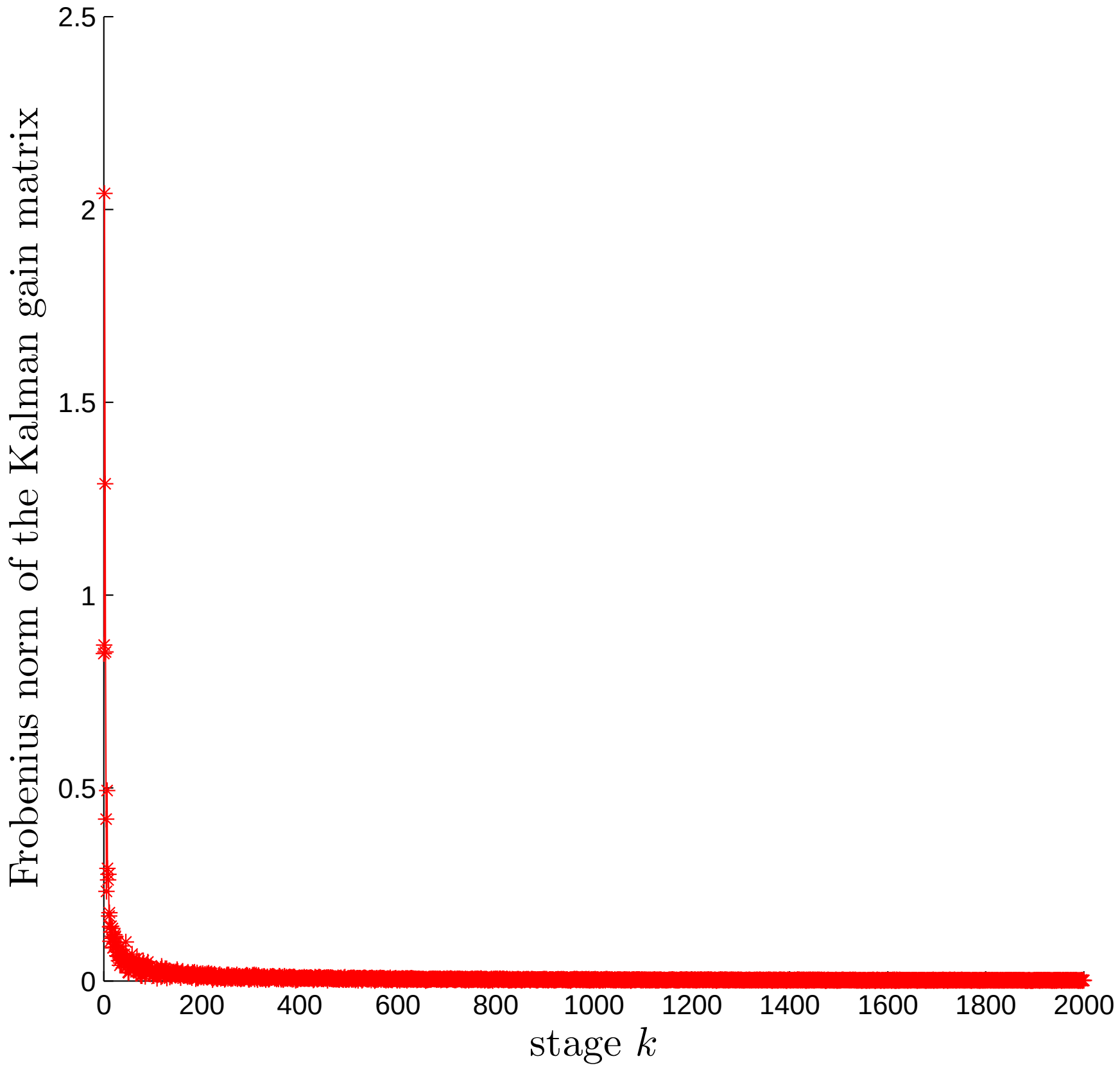

Figure 8: Frobenius norm of the Kalman gain matrix $H_k$ for the experiment shown in Figure 6 ( highly "dense" regions correspond to oscillations in the Frobenius norm with respect to the time index $k$).

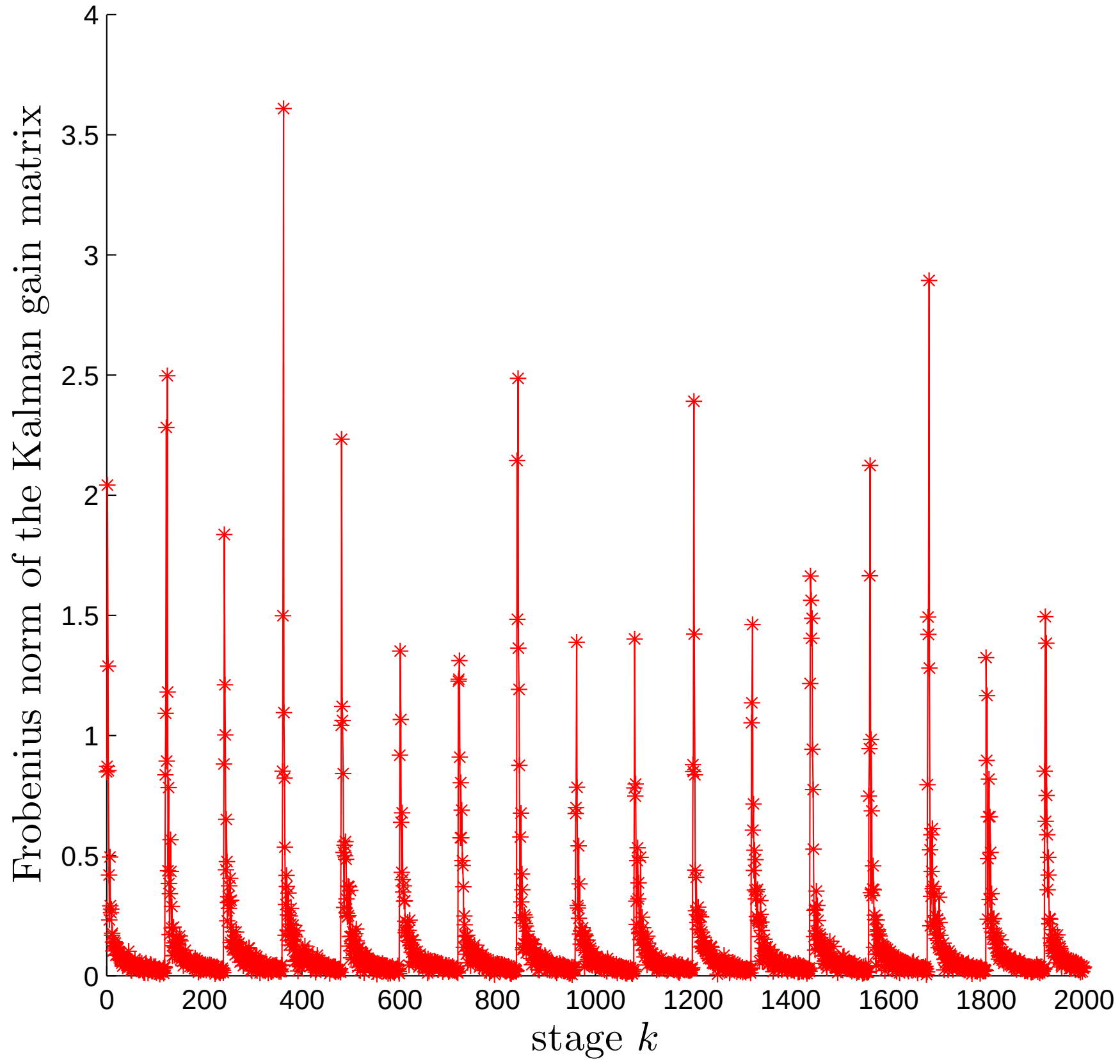

Figure 9: Frobenius norm of the Kalman gain matrix $H_k$ for the experiment shown in Figure 7.

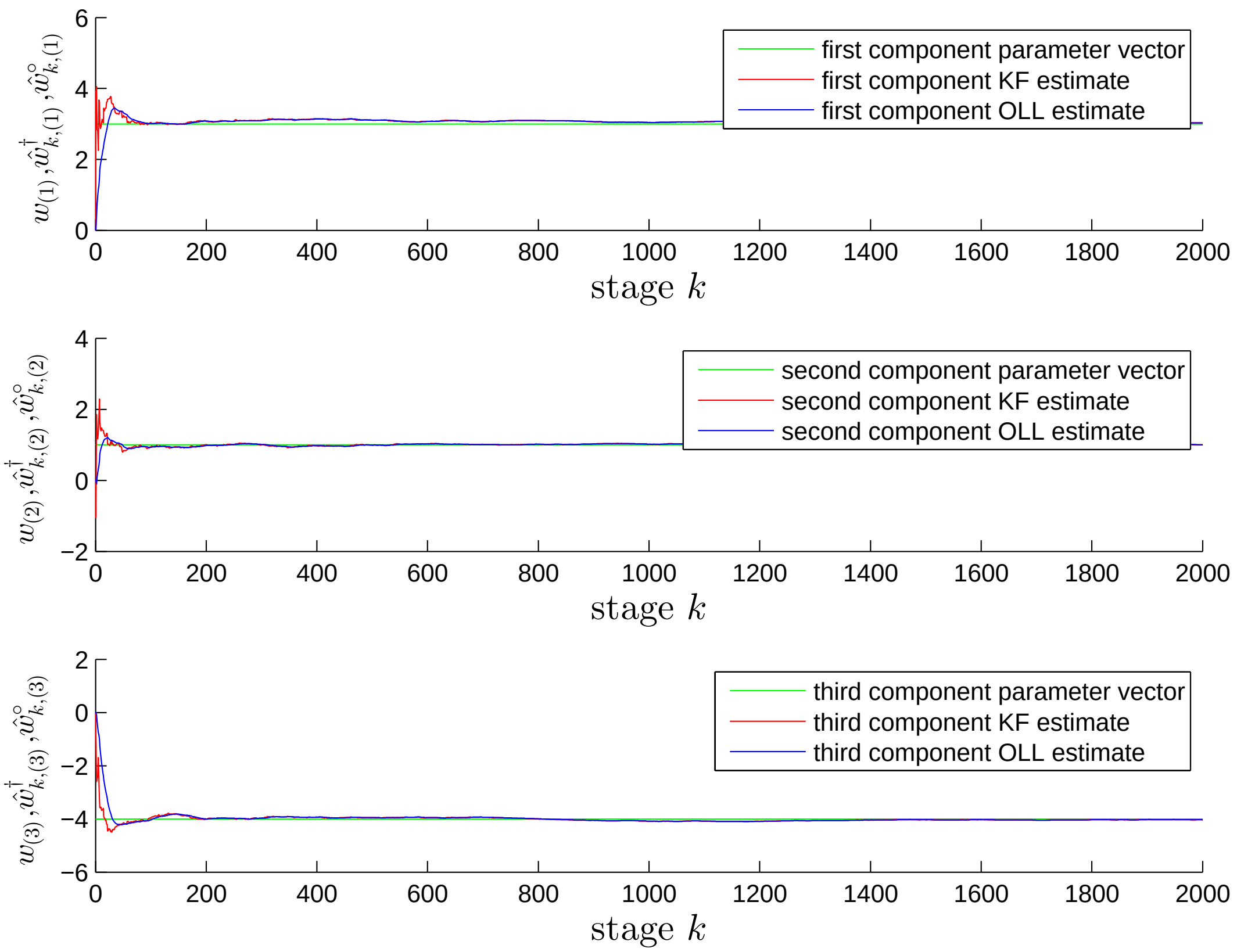

Figure 10: An experimental setting similar to the one of Figure 6, but with no change in the parameter vector.

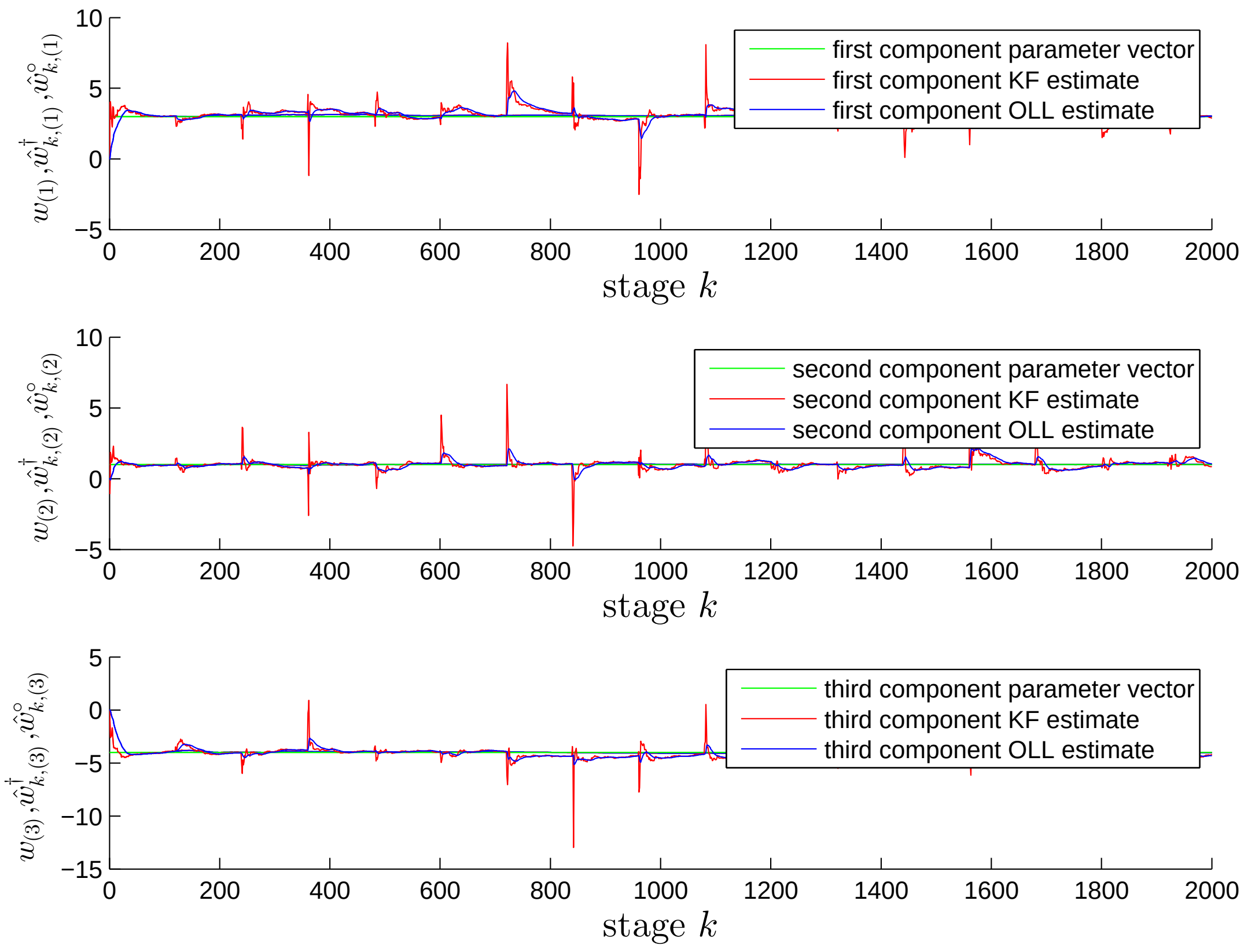

Figure 11: A comparison similar to the one of Figure 10, with the same online examples, and with a re-initialization to $\Sigma_w$ of the matrices $\Sigma_k$ performed every 120 online examples.

(14), (15) and (16), and the ARE (17) become, respectively,

$$L_k := -(\rho\overline{K}_{k+1} + \gamma I)^{-1}\rho\overline{K}_{k+1}\,, \tag{78}$$

$$K_k := \rho\overline{K}_{k+1} - \rho\overline{K}_{k+1}(\rho\overline{K}_{k+1} + \gamma I)^{-1}\rho\overline{K}_{k+1} + Q_k\,, \tag{79}$$

$$F_k := \rho\overline{K}_{k+1}(\rho\overline{K}_{k+1} + \gamma I)^{-1}\rho\overline{K}_{k+1}\,, \tag{80}$$

and

$$\overline{K}_k := \mathop{\mathbb{E}}_{K_k}\{K_k\} = \rho\overline{K}_{k+1} - \rho\overline{K}_{k+1}(\rho\overline{K}_{k+1} + \gamma I)^{-1}\rho\overline{K}_{k+1} + \mathop{\mathbb{E}}_{Q_k}\{Q_k\}\,, \tag{81}$$

whereas the SRE (25) is not changed at all. For the infinite-horizon case, the learning functional is the limit of (77) for $N \to +\infty$ (so, no "average" learning functional similar to (41) is needed). Moreover, the stationary matrix

$$L := -(\rho\overline{K} + \gamma I)^{-1}\rho\overline{K} \tag{82}$$

associated with the stationary solution $\overline{K}$ of the stationary version

$$\overline{K} = \rho\overline{K} - \rho\overline{K}(\rho\overline{K} + \gamma I)^{-1}\rho\overline{K} + \overline{Q}\,, \tag{83}$$

of the ARE (81), is symmetric and such that $(I + L)$ has all its eigenvalues inside the unit circle, as it follows from the positive-definiteness of $\overline{K}$, which is derived directly from (83), assuming that $\overline{Q}$ is positive-definite, likewise in Section 4. Concluding, these analogies with the undiscounted case would allow one to extend the properties stated in Proposition 6 to the discounted case.

As regards the modification (77), one can also observe that, likewise for the basic version of Problem OLL, the past updates of the estimate of the parameter vector $w$ are not modified when a new example arrives, and that, in any case, such updates have no influence on the cost-to-go functions at later stages. Moreover, such a discount factor should be not confused with a forgetting factor[20] (which gives less importance to the past errors, to determine the current estimate).

*g)* **Introduction of additional regularizations on the estimates of** $w$: given a sequence of additional regularization parameters $\gamma_{\hat{w}_k} > 0$, one could insert terms of the form $\gamma_{\hat{w}_k}\hat{w}_k'\hat{w}_k$ in the per-stage cost of the learning functional (5), hence replacing it by

$$\begin{aligned}
& J^N_{\gamma,\gamma_{\hat{w}}}\left(\{u_k(I_k)\}_{k=0}^{N-1}\right) \\
:= & \mathop{\mathbb{E}}_{e_0,\{x_k\}_{k=0}^{N},\{\varepsilon_k\}_{k=0}^{N-1}}\left\{\sum_{k=0}^{N-1}\left[\left((\hat{w}_k - w_k)'x_k\right)^2 + \gamma_{\hat{w}_k}\hat{w}_k'\hat{w}_k + \gamma u_k'u_k\right] + \left((\hat{w}_N - w_N)'x_N\right)^2\right\} \\
= & \mathop{\mathbb{E}}_{e_0,\{x_k\}_{k=0}^{N},\{\varepsilon_k\}_{k=0}^{N-1}}\left\{\sum_{k=0}^{N-1}[(\hat{w}_k - w_k)'Q_k(\hat{w}_k - w_k) + \gamma_{\hat{w}_k}\hat{w}_k'\hat{w}_k + \gamma u_k'u_k]\right. \\
& \left. +(\hat{w}_N - w_N)'Q_N(\hat{w}_N - w_N)\right\}.
\end{aligned} \tag{84}$$

20. Also the case $\rho > 1$, not studied here because the resulting problem could be not well-defined in the infinite-horizon case, does not correspond to a forgetting factor, for the same reasons as above.

In this case, one would still obtain an ARE, but referred to the original dynamical system (3) instead than the reduced one (7), due to the presence of the terms $\gamma_{\hat{w}_k}\hat{w}_k'\hat{w}_k$ in the objective functional (84). So, the matrix $L_k$ would have in this case the size $2d \times 2d$. Instead, the SRE would be exactly the same as (25). A suitable choice for the sequence of the regularization parameters $\gamma_{\hat{w}_k}$ could be a sequence that decreases monotonically to 0. In this case, in the first stages - when the machine has seen a small number of examples - it would be preferable to make the a-priori knowledge on $w$ dominate the one associated with the small number of examples seen so far, whereas in the successive stages - having seen a much larger number of examples - it would be preferable to make the a-posteriori knowledge coming from the examples be predominant. This would be also in accordance with common choices of the regularization parameter (as a function of the number of available examples, plus other parameters) for batch learning with regularization [16]. For the infinite-horizon version of the problem, a constant $\gamma_{w_k}$ (denoted by $\gamma_w$) could be used.

Another extension has to do with the insertion in the per-stage cost of terms related to the previous updates. For instance, limiting for simplicity to the last previous update and adding another regularization parameter $\gamma_{-1} > 0$, one could extend the definition of the learning functional (84) to

$$
\begin{aligned}
& J^N_{\gamma,\gamma_{\hat{w}},\gamma_{-1}}\left(\{u_k(I_k)\}_{k=0}^{N-1}\right) \\
:= & \mathop{\mathbb{E}}_{e_0,\{x_k\}_{k=0}^{N},\{\varepsilon_k\}_{k=0}^{N-1}} \Bigg\{ \sum_{k=0}^{N-1} \Bigg[ \left((\hat{w}_k - w_k)'x_k\right)^2 + \gamma_{\hat{w}_k}\hat{w}_k'\hat{w}_k \\
& \qquad + \gamma_{-1}(u_k - (\hat{w}_k - \hat{w}_{k-1}))'(u_k - (\hat{w}_k - \hat{w}_{k-1})) + \gamma u_k' u_k \Bigg] \\
& \qquad + \left((\hat{w}_N - w_N)'x_N\right)^2 \Bigg\} \\
= & \mathop{\mathbb{E}}_{e_0,\{x_k\}_{k=0}^{N},\{\varepsilon_k\}_{k=0}^{N-1}} \Bigg\{ \sum_{k=0}^{N-1} \Bigg[ (\hat{w}_k - w_k)'Q_k(\hat{w}_k - w_k) + \gamma_{\hat{w}_k}\hat{w}_k'\hat{w}_k \\
& \qquad + \gamma_{-1}(u_k - (\hat{w}_k - \hat{w}_{k-1}))'(u_k - (\hat{w}_k - \hat{w}_{k-1})) + \gamma u_k' u_k \Bigg] \\
& \qquad + (\hat{w}_N - w_N)'Q_N(\hat{w}_N - w_N) \Bigg\}, \qquad (85)
\end{aligned}
$$

with $\hat{w}_{-1} := 0$. Here, each component $(u_k - (\hat{w}_k - \hat{w}_{k-1}))_{(j)}$ of the vector $u_k - (\hat{w}_k - \hat{w}_{k-1})$ can be intepreted as a central-difference approximation (with discretization step $\Delta_k = 1$) of the second derivative of a continuous-time function (or more precisely, of the second derivative of a stochastic process [37, Appendix 10A]) $u_{(j)}(t)$ such that $u_{(j)}(k) = u_{k,(j)}$. Indeed, one has

$$
\frac{d}{dt}u_{(j)}(t)\Big|_{t=k} \simeq \frac{\frac{(\hat{w}_{k+1}-\hat{w}_k)_{(j)}}{\Delta_k} - \frac{(\hat{w}_k - \hat{w}_{k-1})_{(j)}}{\Delta_k}}{\Delta_k} = (u_k - (\hat{w}_k - \hat{w}_{k-1}))_{(j)}\,.
$$

In order to optimize the learning functional (85), one could at first extend the definition of the state vector of the dynamical system (3), including also the previous estimate $\hat{w}_{k-1}$ in the extended state at the time $k$. Then, tools similar to the ones used in this work could be still used for the optimization.

One can see that, by choosing a suitable value for $\gamma_{-1}$, one could give different importance to the previous estimates of $w$, in order to generate the optimal update $u_k^\circ$ at each time $k$. Indeed, when $\gamma_{-1} \simeq 0$, one would expect the estimate $\hat{w}_{k-1}$ not to be practically taken into account to generate $u_k^\circ$, whereas, for $\gamma_{-1}$ extremely large, one would expect the learning machine to penalize, for each $j$-th component, a change in the "slope" of $u_{k,(j)}$ more than its absolute value $|(u_{k,(j)}|$. In other words, to generate $u_{k,(j)}^\circ$, the difference $(\hat{w}_k - \hat{w}_{k-1})_{(j)}$ would be taken into account more than $|(u_{k,(j)}|$. For intermediate values of $\gamma_{-1}$, both $(\hat{w}_k - \hat{w}_{k-1})_{(j)}$ and $|(u_{k,(j)}|$ would be expected to be taken into account significantly to generate $u_{k,(j)}^\circ$. A more rigorous analysis of these cases could be done by solving the ARE for the specific problem, and is outside the scope of this work.

Finally, a similar technique could be used to include terms approximating derivatives of larger order of $u_{(j)}(t)$, extending the definition of the state vector of the dynamical system (3) in such a way to include all the previous estimates of $w$ that are used to approximate such derivatives.

*h)* **Extension of the problem formulation to the continuous-time case**: as already mentioned in the item *g)*, some terms in the learning functional (84) can be interpreted as approximations of terms arising in a continuous-time formulation of the problem. Also the basic version of Problem OLL can be considered as a discrete-time version of a continuous-time problem whose learning functional is

$$\begin{aligned}
& J_{\gamma,c}^{T}\left(u(\cdot,\cdot)\right) \\
:= & \mathop{\mathbb{E}}_{e_0,\{x(\cdot)\},\{\varepsilon(\cdot)\}_{k=0}^{N-1}} \left\{ \int_0^T \left[ \left( (\hat{w}(t) - w(t))' x(t) \right)^2 + \gamma (u(t,I_t))'(u(t,I_t)) \right] dt + \left( (\hat{w}(T) - w(T))' x(T) \right)^2 \right\} \\
= & \mathop{\mathbb{E}}_{e_0,\{x(\cdot)\},\{\varepsilon(\cdot)\}_{k=0}^{N-1}} \left\{ \int_0^T \left[ (\hat{w}(t) - w(t))' Q(t) (\hat{w}(t) - w(t)) + \gamma (u(t,I_t))'(u(t,I_t)) \right] dt \right. \\
& \qquad \left. + (\hat{w}(T) - w(T)' Q(T) (\hat{w}(T) - w(T)) \right\},
\end{aligned} \tag{86}$$

with $T = N/\Delta_k$ and[21] $\Delta_k = 1$, where $u, w, \hat{w}, x, Q$ are suitable stochastic processes, and the system dynamics are described by the stochastic differential equation

$$\begin{cases} dw = 0\,, \\ d\hat{w} = u\,dt\,. \end{cases} \tag{87}$$

Moreover, the measurement process would be modeled by the stochastic process

$$y = Cw + \varepsilon\,, \tag{88}$$

21. A discretization of (86) better than (5) would be obtained using a smaller value of $\Delta_k$, but it would have exactly the same form as (5).

where $\varepsilon$ is another stochastic process (white noise). Finally, the causality constraint that the current update depends only on the "history" of the measurement and decision processes up to the current time would be enforced imposing that the stochastic process $u$ is non-anticipative. As it would be rather technical, a rigorous analysis of this case is outside the scope of the work, but it could be done using classical tools, such as the ones used to solve the LQG optimal control problem in continuous time [34, Chapter 14].

*i)* **Introduction of constraints on the updates**: as a possible extension, one could insert constraints on the variable $u_k$, such as

$$\|u_k\|_2 \leq B_k$$

for some $B_k > 0$ (where $\|\cdot\|_2$ denotes the $l_2$-norm), or constraints of the form

$$|u_{k,(j)}| \leq B_j$$

for some $B_j > 0$, $j = 1, \ldots, d$. In particular, the latter are linear constraints, as they can be written as the union of the constraints

$$u_{k,(j)} \leq B_j$$

and

$$-u_{k,(j)} \leq B_j$$

(we refer to [25–27] for other examples of constraints in machine learning problems, although presented in a batch framework therein). From a theoretical point of view, one could still search for the optimal solution of the resulting constrained optimization problem by solving Bellman equations, provided that one is able to determine the conditional probability distribution of $e_k$ given $\tilde{I}_k$. Under Assumption 5, such a conditional probability distribution is still Gaussian, so one needs to know only the conditional mean $\mathbb{E}_{e_k}\{e_k | \tilde{I}_k\}$ and the conditional covariance matrix $\Sigma_k = \mathbb{E}_{e_k}\left\{(e_k - \mathbb{E}_{e_k}\{e_k|\tilde{I}_k\})(e_k - \mathbb{E}_{e_k}\{e_k|\tilde{I}_k\})'\right\}$ (which are provided by the Kalman-filter recursion scheme) to determine such a conditional probability distribution completely. However, in this case, solving Bellman equations would be not reduced to solving suitable AREs. In practice, being able to solve the problem optimally may be the exception, particularly, in the case of a large finite horizon (or of an infinite horizon), since the complexity of the structure of the optimal-cost-to-go functions - e.g., the number of "pieces" in case of their possible piecewise-quadraticity - may grow up when performing the backward phase of dynamic programming. So, in practice one may be forced to give up searching for an optimal solution, and look for good suboptimal solutions instead.

*j)* **Modification of the per-stage cost**: besides the changes already discussed in the item *g)*, one could also modify the per-stage cost by inserting additive non-quadratic but still convex terms, e.g., a LASSO (Least Absolute Shrinkage and Selection Operator) term of the form $\gamma_L \|u_k\|_1$, where $\gamma_L > 0$ and $\|\cdot\|_1$ denotes the $l_1$-norm. The goal of the LASSO is to enforce sparsity of the update $u_k$ at optimality, which is possible due to geometrical properties of the $l_1$-norm [49]. Here, we observe that, likewise the case discussed in the item *i)* above, at least from a theoretical point of view it is still possible to solve the problem by an application of dynamic programming, if one is able to express the conditional probability distribution of $e_k$ given $\tilde{I}_k$. This is the case, because, under Assumption 5, such a conditional probability distribution is still Gaussian, and can be computed efficiently through the Kalman-filter recursion scheme.

*k)* **Introduction of constraints on the parameter vector $w$ to be learned**: the parameter vector $w$ could be subject to other constraints (e.g., the non-negativity constraint $w \geq 0$) modeling an additional prior knowledge on it. Again, dynamic programming can be still applied in such a variation of Problem OLL, if one is able to express the conditional probability distribution of $e_k$ given both $\tilde{I}_k$ and the constraint on $w$. In practice, this is the case, since, under Assumption 5, such a distribution is obtained from the Gaussian conditional probability distribution of $e_k$ given $\tilde{I}_k$ only, by imposing the constraint $w \geq 0$, then doing a successive renormalization of that Gaussian conditional probability distribution.

*l)* **Application of moving-horizon techniques**: in practice, in certain situations, it could be very difficult to solve exactly the modified problems discussed in the items *i)*, *j)*, and *k)*. In such cases, one could resort to variations of such problems, obtained following a moving-horizon approach [11], which uses a sliding optimization window of constant width, and the current estimate as the "initial" estimate at the left extreme of the current optimization window. Such an approach would assign no importance to far-in-the-future expected errors, outside the current optimization window. Moreover, as it is typical of moving horizon approaches, it may allow one to find possibly good and stable suboptimal solutions to the original infinite-horizon problems; such solutions could be computed by solving - possibly in real-time - simpler (and still convex) optimization problems, especially if the width of the optimization window is small.

*m)* **Reducing the complexity of the estimate, and downdating**: a similar variation as the one discussed in item *l)*, which also uses a sliding optimization window - but with a re-initialization at 0 (or at a fixed vector) of the estimate of the parameter vector $w$ at the left extreme of each window - has to do with the case in which one wants to forget completely the "old" examples, making the current estimate depend only on the examples contained in the current window. For a small width of the window, this would have the advantage of limiting the "complexity" of the current estimate of the parameter vector $w$[22]. Moreover, when shifting the sliding window one unit to the right (hence, inserting a new example, and removing the oldest one), a recursive approach could be used to generate the optimal solution of the resulting optimization problem, starting from the one of the previous problem. This approach - called "downdating" in the literature (as opposed to "updating")

22. As shown in Section 7, for every $k = 1, 2, \ldots$, the vector $\hat{w}_k^\circ$ belongs to a finite-dimensional subspace of dimension at most $k$. So, considering only the examples contained in the current sliding window would make each estimate belong to a finite-dimensional subspace of dimension at most equal to the width of the window. However, it has also to be taken into account that, in the basic version of Problem OLL (the one without the nonlinear mapping $\phi$), the maximal dimension of the subspace is equal to the dimension $d$ of the input space, whereas in its kernel version it is equal to $d_E$. So, this sliding-window approach is expected to reduce significantly the "complexity" of the estimate only when the size of the window is small compared to $d$ in the basic version of the problem, and to $d_E$ in its kernel version.

- works, e.g., for recursive least squares (see, e.g., [50, Section 5.4.1]) and may be extended to the present problem.

*n)* **Active online learning**: one could give to the learning machine the capability of influencing (at least partially) the choice of the sequence of input data. For instance, the learning machine could try to generate examples similar to the ones already seen, possibly making it easier to estimate the statistics of the measurement noises $\varepsilon_k$, or it could focus - at least in an initial learning stage - on certain components of the parameter vector $w$ (e.g., components that seem to be easier to be learned), generating examples with 0 values for all the other components. Focusing on such components in the initial stages could improve the convergence of the estimate to the true parameter vector $w$, making the machine learn the "more difficult" components of $w$ in a second phase. However, in doing this, one should give to the learning machine not an excessive freedom to generate its input examples (e.g., giving it enough freedom only at certain time instants) - in order to explore the state space enough, and to focus not only on tasks that appear easy to be learned. Finally, some of the input examples could be chosen deterministically and even presented periodically to the machine, with small changes in the analysis.

*o)* **Extension of the problem formulation through techniques from robust estimation and control**: as a last possible extension, we discuss briefly an adversarial framework (e.g., both the input examples and the output disturbance noise - the latter not modeled anymore as Gaussian random vectors - could be chosen in an adversarial way), in which one still wants to have the ability to learn the parameter vector even in a (suitably defined) worst-case setting. To study this possible extension, techniques from robust estimation/control could be used, particularly, the ones from $\mathcal{H}_\infty$-filtering/control [4], which - incidentally - is also based on suitable Riccati equations. Still, a direct application of such methods to online machine learning would be not trivial, e.g., in case not only the disturbance noise, but also also the input examples were chosen in an adversarial way. Compared with the LQ/LQG online learning framework investigated in the paper, a possible extension to online robust estimation/control could have the advantage of further decrease the sensitivity to outliers, since the worst case would be considered explicitly to generate the estimates of the parameter vector. However, differently from the present setting, closed-form optimal solutions may not be available for such an extension.

## 9. Discussion

We have proposed and investigated an optimal-control approach to online learning from supervised examples, modeled as the online-estimation of an unknown parameter relating the input examples $x_k$ with their outputs $y_k$. We have shown the connections of the proposed problem with the classical LQ and LQG optimal control problems, of which the former is a non-trivial variation, as it involves random matrices. We have also compared the optimal solution with the KF estimate, showing cases in which the latter has advantages on it (e.g., more smoothness and less sensitivity to outliers). We have also described, and in some cases, developed in details, some extensions of the basic model, including

a) the infinite-horizon case, with convergence results (in particular, convergence to 0 of the mean-square estimation error of the OLL estimate, when the time index goes to infinity);

b) nonlinear models, exploiting kernel methods;

c) nonzero-mean random variables;

d) more complex models for the measurement errors;

e) online estimates of some covariance matrices;

f) a slowly time-varying parameter vector to be learned from the sequence of supervised examples;

g) discounted problems;

h) higher-order regularizations of the estimates of $w$;

i) continuous time;

j) active online learning.

## Appendix: proofs

The proof of Proposition 1 is based on the following lemma.

**Lemma 1** *For $k = N-2, \ldots, 0$ one has*

$$\begin{aligned}&\mathop{\mathbb{E}}_{e_{k+1},\tilde{I}_{k+1}} \left\{ \left( e_{k+1} - \mathop{\mathbb{E}}_{e_{k+1}} \left\{ e_{k+1} \middle| \tilde{I}_{k+1} \right\} \right)' F_{k+1} \left( e_{k+1} - \mathop{\mathbb{E}}_{e_{k+1}} \left\{ e_{k+1} \middle| \tilde{I}_{k+1} \right\} \right) \middle| \tilde{I}_k, u_k \right\} \\ = \ & \mathrm{Tr} \left\{ F_{k+1} \mathop{\mathbb{E}}_{C_{k+1}} \left\{ \Sigma_k - \Sigma_k C'_{k+1} (C_{k+1} \Sigma_k C'_{k+1} + \sigma_\varepsilon^2)^{-1} C_{k+1} \Sigma_k \middle| C_0, \ldots, C_k \right\} \right\} .\end{aligned}$$

**Proof.** By the law of iterated expectations [10, Section 3.2], we get

$$\begin{aligned}&\mathop{\mathbb{E}}_{e_{k+1},\tilde{I}_{k+1}} \left\{ \left( e_{k+1} - \mathop{\mathbb{E}}_{e_{k+1}} \left\{ e_{k+1} \middle| \tilde{I}_{k+1} \right\} \right)' F_{k+1} \left( e_{k+1} - \mathop{\mathbb{E}}_{e_{k+1}} \left\{ e_{k+1} \middle| \tilde{I}_{k+1} \right\} \right) \middle| \tilde{I}_k, u_k \right\} \\ = \ & \mathop{\mathbb{E}}_{\tilde{I}_{k+1}} \left\{ \mathop{\mathbb{E}}_{e_{k+1}} \left\{ \left( e_{k+1} - \mathop{\mathbb{E}}_{e_{k+1}} \left\{ e_{k+1} \middle| \tilde{I}_{k+1} \right\} \right)' F_{k+1} \left( e_{k+1} - \mathop{\mathbb{E}}_{e_{k+1}} \left\{ e_{k+1} \middle| \tilde{I}_{k+1} \right\} \right) \middle| \tilde{I}_{k+1} \right\} \middle| \tilde{I}_k, u_k \right\} .\end{aligned} \tag{89}$$

Now, by properties of the trace and the linearity of the trace and of the expectation operator, one has

$$
\begin{aligned}
& \mathop{\mathbb{E}}_{e_{k+1}} \left\{ \left( e_{k+1} - \mathop{\mathbb{E}}_{e_{k+1}} \left\{ e_{k+1} \middle| \tilde{I}_{k+1} \right\} \right)' F_{k+1} \left( e_{k+1} - \mathop{\mathbb{E}}_{e_{k+1}} \left\{ e_{k+1} \middle| \tilde{I}_{k+1} \right\} \right) \middle| \tilde{I}_{k+1} \right\} \\
= & \mathop{\mathbb{E}}_{e_{k+1}} \left\{ \mathrm{Tr} \left\{ \left( e_{k+1} - \mathop{\mathbb{E}}_{e_{k+1}} \left\{ e_{k+1} \middle| \tilde{I}_{k+1} \right\} \right)' F_{k+1} \left( e_{k+1} - \mathop{\mathbb{E}}_{e_{k+1}} \left\{ e_{k+1} \middle| \tilde{I}_{k+1} \right\} \right) \right\} \middle| \tilde{I}_{k+1} \right\} \\
= & \mathop{\mathbb{E}}_{e_{k+1}} \left\{ \mathrm{Tr} \left\{ F_{k+1} \left( e_{k+1} - \mathop{\mathbb{E}}_{e_{k+1}} \left\{ e_{k+1} \middle| \tilde{I}_{k+1} \right\} \right) \left( e_{k+1} - \mathop{\mathbb{E}}_{e_{k+1}} \left\{ e_{k+1} \middle| \tilde{I}_{k+1} \right\} \right)' \right\} \middle| \tilde{I}_{k+1} \right\} \\
= & \mathrm{Tr} \left\{ F_{k+1} \mathop{\mathbb{E}}_{e_{k+1}} \left\{ \left( e_{k+1} - \mathop{\mathbb{E}}_{e_{k+1}} \left\{ e_{k+1} \middle| \tilde{I}_{k+1} \right\} \right) \left( e_{k+1} - \mathop{\mathbb{E}}_{e_{k+1}} \left\{ e_{k+1} \middle| \tilde{I}_{k+1} \right\} \right)' \middle| \tilde{I}_{k+1} \right\} \right\} \\
= & \mathrm{Tr} \left\{ F_{k+1} \Sigma_{k+1} \right\} . \qquad (90)
\end{aligned}
$$

Due to equations (15), (16), (17) and their initializations (18) and (19), the last expression in (90) is a function of the form $\tilde{f}_{k+1}(\{C_j\}_{j=0}^{k+1})$. Finally, by combining (25), (89), and (90), one gets

$$
\begin{aligned}
& \mathop{\mathbb{E}}_{e_{k+1}} \left\{ \left( e_{k+1} - \mathop{\mathbb{E}}_{e_{k+1}} \left\{ e_{k+1} \middle| \tilde{I}_{k+1} \right\} \right)' F_{k+1} \left( e_{k+1} - \mathop{\mathbb{E}}_{e_{k+1}} \left\{ e_{k+1} \middle| \tilde{I}_{k+1} \right\} \right) | \tilde{I}_k, u_k \right\} \\
= & \mathop{\mathbb{E}}_{\Sigma_{k+1}} \left\{ \mathrm{Tr} \left\{ F_{k+1} \Sigma_{k+1} \right\} \middle| \tilde{I}_k, u_k \right\} \\
= & \mathrm{Tr} \left\{ F_{k+1} \mathop{\mathbb{E}}_{C_{k+1}} \left\{ \Sigma_k - \Sigma_k C'_{k+1} (C_{k+1} \Sigma_k C'_{k+1} + \sigma_\varepsilon^2)^{-1} C_{k+1} \Sigma_k \middle| C_0, \ldots, C_k \right\} \right\} . \qquad (91)
\end{aligned}
$$

■

According to Lemma 1, the terms

$$
\mathop{\mathbb{E}}_{e_{k+1}, \tilde{I}_{k+1}} \left\{ \left( e_{k+1} - \mathop{\mathbb{E}}_{e_{k+1}} \left\{ e_{k+1} \middle| \tilde{I}_{k+1} \right\} \right)' F_{k+1} \left( e_{k+1} - \mathop{\mathbb{E}}_{e_{k+1}} \left\{ e_{k+1} \middle| \tilde{I}_{k+1} \right\} \right) \middle| \tilde{I}_k, u_k \right\} , \qquad (92)
$$

are functions of the form $f_k(\{C_j\}_{j=0}^{k})$ of the random matrices $C_j$, for $j = 0, \ldots, k$

**Proof of Proposition 1**

Likewise in the classical derivations of the cost-to-go functions for the LQ problem shown in [7, Section 4.1], first we solve the Bellman equation (12) for $k = N-1$ and $k = N-2$, then we infer the form of its solution for $k = N-3, \ldots, 0$. For $k = N-1$, one obtains

$$
\begin{aligned}
& \tilde{J}^\circ_{N-1}(\tilde{I}_{N-1}) \\
= & \inf_{u_{N-1} \in \mathbb{R}^d} \mathop{\mathbb{E}}_{e_{N-1}, \tilde{I}_N} \left\{ e'_{N-1} Q_{N-1} e_{N-1} + \gamma u'_{N-1} u_{N-1} + \tilde{J}^\circ_N(\tilde{I}_N) \middle| \tilde{I}_{N-1}, u_{N-1} \right\} \\
= & \inf_{u_{N-1} \in \mathbb{R}^d} \mathop{\mathbb{E}}_{e_{N-1}, Q_N} \left\{ e'_{N-1} Q_{N-1} e_{N-1} + \gamma u'_{N-1} u_{N-1} + (e_{N-1} + u_{N-1})' Q_N (e_{N-1} + u_{N-1}) \middle| \tilde{I}_{N-1}, u_{N-1} \right\} \\
= & \mathop{\mathbb{E}}_{e_{N-1}} \left\{ e'_{N-1} Q_{N-1} e_{N-1} \middle| \tilde{I}_{N-1} \right\} \\
& + \inf_{u_{N-1} \in \mathbb{R}^d} \mathop{\mathbb{E}}_{e_{N-1}, Q_N} \left\{ \gamma u'_{N-1} u_{N-1} + (e_{N-1} + u_{N-1})' Q_N (e_{N-1} + u_{N-1}) \middle| \tilde{I}_{N-1}, u_{N-1} \right\} . \qquad (93)
\end{aligned}
$$

For uniformity of notation with some of the next equations, from now on we set $K_N := Q_N$. Now, we observe that $K_N$ is conditionally mutually independent from $e_{N-1}$ and $u_{N-1}$ given $\tilde{I}_{N-1}$ and $u_{N-1}$ and is mutually independent from $\tilde{I}_{N-1}$ and $u_{N-1}$. Hence, by setting

$$\overline{K}_N := \mathop{\mathbb{E}}_{K_N}\{K_N\}$$

(which is a symmetric and positive-semidefinite matrix), one gets

$$\begin{aligned} &\mathop{\mathbb{E}}_{e_{N-1},K_N}\{(e_{N-1}+u_{N-1})'K_N(e_{N-1}+u_{N-1})\big|\tilde{I}_{N-1},u_{N-1}\} \\ = &\mathop{\mathbb{E}}_{e_{N-1}}\{(e_{N-1}+u_{N-1})'\overline{K}_N(e_{N-1}+u_{N-1})\big|\tilde{I}_{N-1},u_{N-1}\}. \end{aligned} \tag{94}$$

Combining (93) and (94), one has

$$\begin{aligned} &\tilde{J}^{\circ}_{N-1}(\tilde{I}_{N-1}) \\ = &\mathop{\mathbb{E}}_{e_{N-1}}\{e'_{N-1}(Q_{N-1}+\overline{K}_N)e_{N-1}\big|\tilde{I}_{N-1}\} \\ &+\inf_{u_{N-1}\in\mathbb{R}^d}\left[u'_{N-1}(\overline{K}_N+\gamma I)u_{N-1}+2(\overline{K}_N\mathop{\mathbb{E}}_{e_{N-1}}\{e_{N-1}\big|\tilde{I}_{N-1}\})'u_{N-1}\right], \end{aligned} \tag{95}$$

where $I$ denotes the $d\times d$ identity matrix. Now, the matrix $\overline{K}_N+\gamma I$ is symmetric and positive-definite, hence, by the first-order optimal condition, the optimal updating function $u^{\circ}_{N-1}(\tilde{I}_{N-1})$ in equation (95) is given by

$$u^{\circ}_{N-1}(\tilde{I}_{N-1}) = L_{N-1}\mathop{\mathbb{E}}_{e_{N-1}}\left\{e_{N-1}\big|\tilde{I}_{N-1}\right\},$$

where

$$L_{N-1} := -(\overline{K}_N+\gamma I)^{-1}\overline{K}_N\,.$$

Moreover, by putting $u_{N-1}=u^{\circ}_{N-1}(\tilde{I}_{N-1})$ into (95), one obtains

$$\begin{aligned} &\tilde{J}^{\circ}_{N-1}(\tilde{I}_{N-1}) \\ = &\mathop{\mathbb{E}}_{e_{N-1}}\left\{e'_{N-1}K_{N-1}e_{N-1}\big|\tilde{I}_{N-1}\right\} \\ &+\mathop{\mathbb{E}}_{e_{N-1}}\left\{\left(e_{N-1}-\mathop{\mathbb{E}}_{e_{N-1}}\left\{e_{N-1}\big|\tilde{I}_{N-1}\right\}\right)'F_{N-1}\left(e_{N-1}-\mathop{\mathbb{E}}_{e_{N-1}}\{e_{N-1}\big|\tilde{I}_{N-1}\right)\Big|\tilde{I}_{N-1}\right\}, \end{aligned} \tag{96}$$

where

$$K_{N-1} := \overline{K}_N-\overline{K}_N(\overline{K}_N+\gamma I)^{-1}\overline{K}_N+Q_{N-1}\,, \tag{97}$$

and

$$F_{N-1} := \overline{K}_N(\overline{K}_N+\gamma I)^{-1}\overline{K}_N \tag{98}$$

are symmetric and positive-semidefinite matrices. Similarly, for the stage $k=N-2$, the Bellman equation (12) becomes

$$\begin{aligned}
&\tilde{J}^{\circ}_{N-2}(\tilde{I}_{N-2}) \\
= \; & \inf_{u_{N-2}\in\mathbb{R}^d} \mathop{\mathbb{E}}_{e_{N-2},\tilde{I}_{N-1}} \left\{ e'_{N-2}Q_{N-2}e_{N-2} + \gamma u'_{N-2}u_{N-2} + \tilde{J}^{\circ}_{N-1}(\tilde{I}_{N-1}) \middle| \tilde{I}_{N-2}, u_{N-2} \right\} \\
= \; & \mathop{\mathbb{E}}_{e_{N-2}} \left\{ e'_{N-2}Q_{N-2}e_{N-2} \middle| \tilde{I}_{N-2} \right\} \\
& + \inf_{u_{N-2}\in\mathbb{R}^d} \Bigg[ \mathop{\mathbb{E}}_{e_{N-1},K_{N-1}} \left\{ \gamma u'_{N-2}u_{N-2} + e'_{N-1}K_{N-1}e_{N-1} \middle| \tilde{I}_{N-2}, u_{N-2} \right\} \\
& \quad + \mathop{\mathbb{E}}_{e_{N-1},\tilde{I}_{N-1}} \left\{ (e_{N-1} - \mathop{\mathbb{E}}_{e_{N-1}}\{e_{N-1}|\tilde{I}_{N-1}\})' F_{N-1} (e_{N-1} - \mathop{\mathbb{E}}_{e_{N-1}}\{e_{N-1}|\tilde{I}_{N-1}\}) \middle| \tilde{I}_{N-2}, u_{N-2} \right\} \Bigg] .
\end{aligned} \tag{99}$$

Now, by [7, Section 5.2, Lemma 2.1], the term

$$\mathop{\mathbb{E}}_{e_{N-1},\tilde{I}_{N-1}} \left\{ \left( e_{N-1} - \mathop{\mathbb{E}}_{e_{N-1}} \left\{ e_{N-1} \middle| \tilde{I}_{N-1} \right\} \right)' F_{N-1} \left( e_{N-1} - \mathop{\mathbb{E}}_{e_{N-1}} \left\{ e_{N-1} \middle| \tilde{I}_{N-1} \right\} \right) \middle| \tilde{I}_{N-2}, u_{N-2} \right\} \tag{100}$$

does not depend on $u_{N-2}$, neither on the sequence of updates applied up to the time $N-2$. This is basically due to the linearity of the dynamical system and of the measurement equation[23]. By Lemma 1, the term (100) is a function $f_{N-2}(\{C_j\}_{j=0}^{N-2})$ of the random matrices $C_j$, for $j = 0, \ldots, N-2$, whose realizations can be derived directly from the information vector $\tilde{I}_{N-2}$. Hence, the term (100) does not influence the search for an optimal update at the time $N-2$. So, one obtains

$$\begin{aligned}
&\tilde{J}^{\circ}_{N-2}(\tilde{I}_{N-2}) \\
= \; & \inf_{u_{N-2}\in\mathbb{R}^d} \mathop{\mathbb{E}}_{e_{N-2},K_{N-1}} \{ \gamma u'_{N-2}u_{N-2} + (e_{N-2}+u_{N-2})'K_{N-1}(e_{N-2}+u_{N-2}) | \tilde{I}_{N-2}, u_{N-2} \} \\
& + f_{N-2}(\{C_j\}_{j=0}^{N-2}) \\
= \; & \inf_{u_{N-2}\in\mathbb{R}^d} \mathop{\mathbb{E}}_{e_{N-2},K_{N-1}} \left\{ \gamma u'_{N-2}u_{N-2} + (e_{N-2}+u_{N-2})'K_{N-1}(e_{N-2}+u_{N-2}) \middle| \tilde{I}_{N-2}, u_{N-2} \right\} \\
& \qquad + \text{a term that does not depend on } u_{N-2} .
\end{aligned} \tag{101}$$

Such an optimization problem has the same nature as the one in (93). Hence, by setting

$$\overline{K}_{N-1} := \mathop{\mathbb{E}}_{K_{N-1}} \{K_{N-1}\} = \overline{K}_N - \overline{K}_N(\overline{K}_N + \gamma I)^{-1}\overline{K}_N + \mathop{\mathbb{E}}_{Q_{N-1}} \{Q_{N-1}\} , \tag{102}$$

the optimal updating function at the time $k = N-2$ is

$$u^{\circ}_{N-2}(\tilde{I}_{N-2}) = L_{N-2} \mathop{\mathbb{E}}_{e_{N-2}} \left\{ e_{N-2} \middle| \tilde{I}_{N-2} \right\} ,$$

where

$$L_{N-2} := -(\overline{K}_{N-1} + \gamma I)^{-1}\overline{K}_{N-1} .$$

23. To prove [7, Section 5.2, Lemma 2.1], it is assumed therein that the matrices $C_k$ are deterministic. However, inspection of the proof shows that the result still holds when they are random matrices, generated according to the model described in this paper.

Moreover, by putting $u_{N-2} = u^{\circ}_{N-2}(\tilde{I}_{N-2})$ into (95), one obtains

$$\begin{aligned}
&\tilde{J}^{\circ}_{N-2}(\tilde{I}_{N-2}) \\
= \; & \mathop{\mathbb{E}}_{e_{N-2}} \left\{ e'_{N-2} K_{N-2} e_{N-2} \middle| \tilde{I}_{N-2} \right\} \\
& + \mathop{\mathbb{E}}_{e_{N-2}} \left\{ \left( e_{N-2} - \mathop{\mathbb{E}}_{e_{N-2}} \left\{ e_{N-2} \middle| \tilde{I}_{N-2} \right\} \right)' F_{N-2} \left( e_{N-2} - \mathop{\mathbb{E}}_{e_{N-2}} \{ e_{N-2} | \tilde{I}_{N-2} \right) \middle| \tilde{I}_{N-2} \right\} \\
& + f_{N-2}(\{C_j\}_{j=0}^{N-2}) \, , \qquad (103)
\end{aligned}$$

where

$$K_{N-2} := \overline{K}_{N-1} - \overline{K}_{N-1}(\overline{K}_{N-1} + \gamma I)^{-1} \overline{K}_{N-1} + Q_{N-2} \, , \qquad (104)$$

and

$$F_{N-2} := \overline{K}_{N-1}(\overline{K}_{N-1} + \gamma I)^{-1} \overline{K}_{N-1}$$

are symmetric and positive-semidefinite matrices. Finally, by (104) the matrix

$$\overline{K}_{N-2} := \mathop{\mathbb{E}}_{K_{N-2}} \{K_{N-2}\} = \overline{K}_{N-1} - \overline{K}_{N-1}(\overline{K}_{N-1} + \gamma I)^{-1} \overline{K}_{N-1} + \mathop{\mathbb{E}}_{Q_{N-2}} \{Q_{N-2}\} \, ,$$

is symmetric and positive-semidefinite.

**Remark 16** By Lemma 1, for $k = N-3, \dots, 0$, also the terms

$$\mathop{\mathbb{E}}_{e_{k+1}, \tilde{I}_{k+1}} \{ (e_{k+1} - \mathop{\mathbb{E}}_{e_{k+1}} \{ e_{k+1} | \tilde{I}_{k+1} \})' F_{k+1} (e_{k+1} - \mathop{\mathbb{E}}_{e_{k+1}} \{ e_{k+1} | \tilde{I}_{k+1} \}) | \tilde{I}_k, u_k \} \qquad (105)$$

do not depend on $u_k$, neither on the sequence of updates applied up to the time $k$. Moreover, Lemma 1 shows that they are functions $f_k(\{C_j\}_{j=0}^{k})$ of the random matrices $C_j$, for $j = 0, \dots, k$. Again, their realizations can be derived directly from the information vector $\tilde{I}_k$.

The same arguments used for the stages $k = N-1$ and $k = N-2$ can be applied to $k = N-3, \dots, 0$. Proceeding in such a way, one gets the following recursion for $k = N-1, \dots, 0$:

$$u^{\circ}_k(\tilde{I}_k) = L_k \mathop{\mathbb{E}}_{e_k} \left\{ e_k \middle| \tilde{I}_k \right\} \, ,$$

where we recall that

$$L_k := -(\overline{K}_{k+1} + \gamma I)^{-1} \overline{K}_{k+1} \, ,$$

$$\begin{aligned}
\tilde{J}^{\circ}_k(\tilde{I}_k) \; = \; & \mathop{\mathbb{E}}_{e_k} \left\{ e'_k K_k e_k \middle| \tilde{I}_k \right\} + \mathop{\mathbb{E}}_{e_k} \left\{ \left( e_k - \mathop{\mathbb{E}}_{e_k} \left\{ e_k \middle| \tilde{I}_k \right\} \right)' F_k \left( e_k - \mathop{\mathbb{E}}_{e_k} \{ e_k | \tilde{I}_k \right) | \tilde{I}_k \right\} \\
& + \mathop{\mathbb{E}}_{\{C_j\}_{j=k+1}^{N-2}} \left\{ \sum_{h=k}^{N-2} f_h(\{C_j\}_{j=0}^{h}) \middle| \tilde{I}_k \right\} \, ,
\end{aligned} \qquad (106)$$

and the matrices

$$K_k := \overline{K}_{k+1} - \overline{K}_{k+1}(\overline{K}_{k+1} + \gamma I)^{-1}\overline{K}_{k+1} + Q_k \,,$$
$$F_k := \overline{K}_{k+1}(\overline{K}_{k+1} + \gamma I)^{-1}\overline{K}_{k+1}$$

and

$$\overline{K}_k := \mathop{\mathbb{E}}_{K_k}\{K_k\} = \overline{K}_{k+1} - \overline{K}_{k+1}(\overline{K}_{k+1} + \gamma I)^{-1}\overline{K}_{k+1} + \mathop{\mathbb{E}}_{Q_k}\{Q_k\}$$

are symmetric and positive-semidefinite. ■

### Proof of Proposition 2

Let us first show that the recursion (25) holds. To this end, we exploit the classical Kalman-filter recursion scheme to the specific problem; see, e.g., [7, Appendix E.3]. This can be done since, at the time $k$, the realization of the random matrix $C_k$ becomes known to the learning machine, hence one can apply the Kalman recursion to compute $\mathop{\mathbb{E}}_{e_k}\left\{e_k \middle| \tilde{I}_k\right\}$. Indeed, such a recursion requires the knowledge of such a matrix at the time $k$, not before. Note that, differently from the analysis in [7, Appendix E.3], the (conditional) covariance matrix $\Sigma_k$ in (20) depends actually on $\tilde{I}_k$ through the realizations of the random matrices $C_j$, for $j = 0, \ldots, k$. Instead, in the deterministic case, there is no dependence of such a covariance matrix on the information vector, so it is just an unconditional covariance matrix.

Now, let $\hat{e}_k^{\circ,\dagger} := \mathop{\mathbb{E}}_{e_k^\circ}\left\{e_k^\circ \middle| \tilde{I}_k\right\}$ . By [7, Appendix E.3], the KF estimate of $e_k^\circ$ at the time $k$, based on the information vector $\tilde{I}_k$, is given by

$$\hat{e}_{k+1}^{\circ,\dagger} = \hat{e}_k^{\circ,\dagger} + L_k\hat{e}_k^{\circ,\dagger} + H_{k+1}\left(\tilde{y}_{k+1}^\circ - C_{k+1}(\hat{e}_k^{\circ,\dagger} + L_k\hat{e}_k^{\circ,\dagger})\right) \,, \tag{107}$$

which is initialized by $\hat{e}_{-1}^{\circ,\dagger} := 0$, where the Kalman gain matrix $H_{k+1}$ is defined in (24). Moreover, since $e_k^\circ = \hat{w}_k^\circ - w$, $e_{k+1}^\circ = \hat{w}_{k+1}^\circ - w$, $\tilde{y}_{k+1} = C_{k+1}\hat{w}_{k+1}^\circ - y_{k+1}$, and $\hat{w}_k^\circ$ is known at the time $k$, the KF estimate $\hat{w}_k^\dagger$ of $w$ at the time $k$, based on the information vector $I_k$, satisfies[24]

$$\hat{e}_k^{\circ,\dagger} = \hat{w}_k^\circ - \hat{w}_k^\dagger \,. \tag{108}$$

Similarly, since $\hat{w}_{k+1}^\circ$, $C_{k+1}$ and $y_{k+1}$ are known at the time $k+1$, the KF estimate $\hat{w}_{k+1}^\dagger$ of $w$ at the time $k+1$, based on the information vector $I_{k+1}$, satisfies

$$\hat{e}_{k+1}^{\circ,\dagger} = \hat{w}_{k+1}^\circ - \hat{w}_{k+1}^\dagger \,. \tag{109}$$

So, $\hat{w}_{k+1}^\dagger$ is derived from (107) by replacement of (108) and (109), obtaining

$$\begin{aligned}
& (\hat{w}_{k+1}^\circ - \hat{w}_{k+1}^\dagger) \\
= \;& (\hat{w}_k^\circ - \hat{w}_k^\dagger) + L_k(\hat{w}_k^\circ - \hat{w}_k^\dagger) \\
& + H_{k+1}\left((C_{k+1}\hat{w}_{k+1}^\circ - y_{k+1}) - C_{k+1}\left((\hat{w}_k^\circ - \hat{w}_k^\dagger) + L_k(\hat{w}_k^\circ - \hat{w}_k^\dagger)\right)\right) \,.
\end{aligned} \tag{110}$$

24. To make the notation uniform, let $\hat{w}_{-1}^\circ := 0$ in (108).

By equations (13) and (108) and the definition of the error vector $e_k := \hat{w}_k - w$, one obtains

$$\hat{w}^{\circ}_{k+1} = \hat{w}^{\circ}_k + L_k(\hat{w}^{\circ}_k - \hat{w}^{\dagger}_k) = \hat{w}^{\circ}_k + L_k(\hat{e}^{\circ}_k - \hat{e}^{\circ,\dagger}_k)\,. \tag{111}$$

This, combined with (110), provides

$$\hat{w}^{\dagger}_{k+1} = \hat{w}^{\dagger}_k + H_{k+1}(y_{k+1} - C_{k+1}\hat{w}^{\dagger}_k)\,,$$

which is initialized by $\hat{w}^{\dagger}_{-1} := 0$. Finally, $L_{-1}$ can be chosen arbitrarily, since it multiplies vectors with all-zero components (e.g., one can choose $L_{-1} := -\left(\bar{K}_0 + \gamma I\right)^{-1}\bar{K}_0$, as in the statement of the proposition). ■

**Remark 17** Note that the update (23) does not depend on the sequence of applied updates. This could have been obtained more directly by considering the evolution of the dynamical system

$$w_{k+1} = w_k \tag{112}$$

only, together with the initial condition $w_0 := w$, and the measurement equation (4).

**Proof of Proposition 3**

The first bound in (29) follows from the fact that $\Sigma_{k+1}$ is a (conditional) covariance matrix. Let $\Sigma_k^{1/2}$ be the symmetric and positive-semidefinite square root of the matrix $\Sigma_k$ and $M_{k+1} := C_{k+1}\Sigma_k^{1/2}$. The second bound follows by (25) and the fact that

$$\begin{aligned}\Sigma_k - \Sigma_{k+1} &= \Sigma_k^{1/2}\left(\Sigma_k^{1/2}C'_{k+1}(C_{k+1}\Sigma_k^{1/2}\Sigma_k^{1/2}C'_{k+1} + \sigma_\varepsilon^2)^{-1}C_{k+1}\Sigma_k^{1/2}\right)\Sigma_k^{1/2}\,,\\ &= \Sigma_k^{1/2}M'_{k+1}(M_{k+1}M'_{k+1} + \sigma_\varepsilon^2)^{-1}M_{k+1}\Sigma_k^{1/2}\,.\end{aligned} \tag{113}$$

Since $M'_{k+1}(M_{k+1}M'_{k+1} + \sigma_\varepsilon^2)^{-1}M_{k+1}$ is symmetric and positive-semidefinite, by (113) $\Sigma_k - \Sigma_{k+1}$ is symmetric and positive-semidefinite, too.

(ii) Defining $N_{k+1} := \Sigma_k - \Sigma_{k+1}$, one has obviously

$$\Sigma_k = \Sigma_{k+1} + N_{k+1}\,, \tag{114}$$

where by (29) all the matrices involved in (114) are symmetric and positive-semidefinite. So, one gets

$$\mathrm{Tr}\{\Sigma_k\} = \mathrm{Tr}\{\Sigma_{k+1}\} + \mathrm{Tr}\{N_{k+1}\} \geq \mathrm{Tr}\{\Sigma_{k+1}\} \geq 0\,,$$

which proves (30). Moreover, by Weyl's inequalities[25] of matrix-perturbation theory, if one orders the eigenvalues of each of the three matrices $\Sigma_k$, $\Sigma_{k+1}$, and $N_{k+1}$ in nondecreasing order taking into account their multiplicities, then for every $j = 1, \ldots, d$ one gets

$$\lambda_j(\Sigma_{k+1}) + \lambda_d(\Sigma_{N_{k+1}}) \leq \lambda_j(\Sigma_k)\,,$$

25. Let $S_1, S_2 \in \mathbb{R}^{d\times d}$ and symmetric, and let their eigenvalues be ordered nondecreasingly with their multiplicities as

$$\lambda_1(S_1) \leq \lambda_2(S_1) \leq \ldots \leq \lambda_j(S_1) \leq \ldots \leq \lambda_d(S_1)\,,$$

By the positive-semidefiniteness of $N_{k+1}$, for every $j = 1, \dots, d$ one has also

$$\lambda_j(\Sigma_{k+1}) \leq \lambda_j(\Sigma_k)\,,$$

which implies (31).

(iii) By [31, Theorem II.5.4], every bounded and monotonic sequence of self-adjoint operators on a Hilbert space converges strongly to a self-adjoint operator. Then, formula (32) is obtained as a finite-dimensional case of such a result.

(iv) This part of the proof is based on the investigation of the limit behavior of equation (113) for $k \to +\infty$, using also the expectation and trace operators. First, we exploit the assumption that the common probability distribution of the random vectors $x_k$ has bounded support. This, together with the definition $M_{k+1} := C_{k+1}\Sigma_k^{1/2}$ and the bound (29), proves the existence of a positive constant $c_1$ such that

$$M_{k+1}M'_{k+1} = C_{k+1}\Sigma_k C'_{k+1} \leq C_{k+1}\Sigma_{-1}C'_{k+1} \leq c_1 \tag{115}$$

with a-priori probability 1. Then, one obtains

$$M'_{k+1}(M_{k+1}M'_{k+1} + \sigma_\varepsilon^2)^{-1}M_{k+1} \succeq (c_1 + \sigma_\varepsilon^2)^{-1}M'_{k+1}M_{k+1}$$

with a-priori probability 1. Moreover,

$$\begin{aligned}\Sigma_k^{1/2}M'_{k+1}(M_{k+1}M'_{k+1} + \sigma_\varepsilon^2)^{-1}M_{k+1}\Sigma_k^{1/2} &\succeq (c_1 + \sigma_\varepsilon^2)^{-1}\Sigma_k^{1/2}M'_{k+1}M_{k+1}\Sigma_k^{1/2} \\ &= (c_1 + \sigma_\varepsilon^2)^{-1}\Sigma_k C'_{k+1}C_{k+1}\Sigma_k\,, \end{aligned} \tag{116}$$

$$\lambda_1(S_2) \leq \lambda_2(S_2) \leq \ldots \leq \lambda_j(S_2) \leq \ldots \leq \lambda_d(S_2)\,.$$

Then, in their simplest form, Weyl's inequalities (see, e.g., [6, Theorem 8.4.11]) state that, for every $j = 1, \dots, d$, one has

$$\lambda_1(S_1) + \lambda_j(S_2) \leq \lambda_j(S_1 + S_2) \leq \lambda_d(S_1) + \lambda_j(S_2)\,.$$

where all the steps in (116) hold with a-priori probability 1[26]. Hence, exploiting properties of the trace operator and the independence between $C'_{k+1}C_{k+1}$ and $\Sigma_k^2$, one obtains

$$\begin{aligned}
& \mathrm{Tr}\left\{\mathop{\mathbb{E}}_{\Sigma_k, M_{k+1}}\{\Sigma_k^{1/2} M'_{k+1}(M_{k+1}M'_{k+1}+\sigma_\varepsilon^2)^{-1}M_{k+1}\Sigma_k^{1/2}\}\right\} \\
\geq\ & (c_1+\sigma_\varepsilon^2)^{-1}\mathrm{Tr}\left\{\mathop{\mathbb{E}}_{\Sigma_k, C_{k+1}}\{\Sigma_k C'_{k+1}C_{k+1}\Sigma_k\}\right\} \\
=\ & (c_1+\sigma_\varepsilon^2)^{-1}\mathrm{Tr}\left\{\mathop{\mathbb{E}}_{\Sigma_k, C_{k+1}}\{\Sigma_k^2 C'_{k+1}C_{k+1}\}\right\} \\
=\ & (c_1+\sigma_\varepsilon^2)^{-1}\mathrm{Tr}\left\{\mathop{\mathbb{E}}_{\Sigma_k}\{\Sigma_k^2\}\mathop{\mathbb{E}}_{C_{k+1}}\{C'_{k+1}C_{k+1}\}\right\} \\
=\ & (c_1+\sigma_\varepsilon^2)^{-1}\mathrm{Tr}\left\{\mathop{\mathbb{E}}_{\Sigma_k}\{\Sigma_k^2\}\mathop{\mathbb{E}}_{Q_{k+1}}\{Q_{k+1}\}\right\} \\
=\ & (c_1+\sigma_\varepsilon^2)^{-1}\mathrm{Tr}\left\{\mathop{\mathbb{E}}_{\Sigma_k}\{\Sigma_k^2\}\overline{Q}\right\} \\
\geq\ & (c_1+\sigma_\varepsilon^2)^{-1}\mathrm{Tr}\left\{\mathop{\mathbb{E}}_{\Sigma_k}\{\Sigma_k^2\}\right\}\lambda_{\min}(\overline{Q})\,,
\end{aligned} \tag{117}$$

where $\lambda_{\min}(\overline{Q})$ denotes the minimum eigenvalue of $\overline{Q}$, which is positive by the assumed positive-definiteness of $\overline{Q}$, and the last inequality in (117) follows by [20, Theorem 1].

At this point, we recall that, for two symmetric and positive-semidefinite $d \times d$ matrices $S_1$ and $S_2$, one has

$$|\mathrm{Tr}\{S_1 S_2\}| \leq \sqrt{\mathrm{Tr}\{S_1^2\}\mathrm{Tr}\{S_2^2\}}\,, \tag{118}$$

which is the Cauchy-Schwarz inequality for the Hilbert-Schmidt norm $\sqrt{\mathrm{Tr}\{S^2\}}$ [14, Chapter IX]. Hence, when $S_1 = \Sigma_k$ and $S_2 = I$, one obtains

$$|\mathrm{Tr}\{\Sigma_k\}| = |\mathrm{Tr}\{\Sigma_k I\}| \leq \sqrt{\mathrm{Tr}\{\Sigma_k^2\}\mathrm{Tr}\{I^2\}} = \sqrt{d\,\mathrm{Tr}\{\Sigma_k^2\}}\,. \tag{119}$$

---

26. To obtain the generalized inequality in (116), we have exploited the fact that, if $S_1$ and $S_2$ are symmetric and positive-semidefinite $d \times d$ matrices such that

$$S_1 \succeq S_2\,,$$

then, for every $d \times d$ matrix $\Sigma$, one has also

$$\Sigma' S_1 \Sigma \succeq \Sigma' S_2 \Sigma\,.$$

This is proved observing that, for every $y \in \mathbb{R}^d$, one has

$$y'\Sigma' S_1 \Sigma y = (\Sigma y)' S_1 (\Sigma y) \geq (\Sigma y)' S_2 (\Sigma y) = y' \Sigma' S_2 \Sigma y\,.$$

Then,

$$
\begin{aligned}
& \left(\mathrm{Tr}\{\overline{\Sigma}\right)^2\} \\
= \;& \left(\lim_{k\to+\infty} \mathop{\mathbb{E}}_{\Sigma_k}\{\mathrm{Tr}\{\Sigma_k\}\}\right)^2 \\
= \;& \lim_{k\to+\infty} \left(\mathop{\mathbb{E}}_{\Sigma_k}\{\mathrm{Tr}\{\Sigma_k\}\}\right)^2 \\
\leq \;& \liminf_{k\to+\infty} \mathop{\mathbb{E}}_{\Sigma_k}\{(\mathrm{Tr}\{\Sigma_k\})^2\} \\
\leq \;& \liminf_{k\to+\infty} \mathop{\mathbb{E}}_{\Sigma_k}\{d\,\mathrm{Tr}\{\Sigma_k^2\}\}\,, \qquad (120)
\end{aligned}
$$

where the last two inequalities derive, respectively, from the convexity of the function $\sqrt{(\cdot)}$ and Jensen's inequality [40, Theorem 3.3], and from (119).

Now, taking traces and expectations, making $k$ tend to $+\infty$, and exploiting equations (113), (117), and (120), we get

$$
\begin{aligned}
0 \;=\;& \mathrm{Tr}\{\overline{\Sigma}\} - \mathrm{Tr}\{\overline{\Sigma}\} \\
=\;& \lim_{k\to+\infty} \mathrm{Tr}\left\{\mathop{\mathbb{E}}_{\Sigma_k}\{\Sigma_k\}\right\} - \lim_{k\to+\infty} \mathrm{Tr}\left\{\mathop{\mathbb{E}}_{\Sigma_{k+1}}\{\Sigma_{k+1}\}\right\} \\
=\;& \lim_{k\to+\infty} \mathrm{Tr}\left\{\mathop{\mathbb{E}}_{\Sigma_k, C_{k+1}}\left\{\Sigma_k^{1/2}\left(\Sigma_k^{1/2}C'_{k+1}(C_{k+1}\Sigma_k^{1/2}\Sigma_k^{1/2}C'_{k+1}+\sigma_\varepsilon^2)^{-1}C_{k+1}\Sigma_k^{1/2}\right)\Sigma_k^{1/2}\right\}\right\} \\
\geq\;& (c_1+\sigma_\varepsilon^2)^{-1}\liminf_{k\to+\infty}\mathrm{Tr}\{\mathop{\mathbb{E}}_{\Sigma_k}\{\Sigma_k^2\}\}\lambda_{\min}(\overline{Q}) \\
\geq\;& (c_1+\sigma_\varepsilon^2)^{-1}\left(\mathrm{Tr}\left\{\overline{\Sigma}\right\}\right)^2 d^{-1}\lambda_{\min}(\overline{Q})\,. \qquad (121)
\end{aligned}
$$

Hence, since $(c_1+\sigma_\varepsilon^2)^{-1}$, $d^{-1}$, and $\lambda_{\min}(\overline{Q})$ are different from 0, (121) implies

$$
\mathrm{Tr}\{\overline{\Sigma}\} = 0\,, \qquad (122)
$$

and also

$$
\overline{\Sigma} = 0\,,
$$

since $\overline{\Sigma}$ is symmetric and positive-semidefinite. This concludes the proof of (34).

(v) Let us denote by $F_{k+1}^{1/2}$ and $N_{k+1}^{1/2}$, respectively, a symmetric and positive-semidefinite square root of the symmetric and positive-semidefinite matrix $F_{k+1}$, and a symmetric and positive-semidefinite square root of the symmetric and positive-semidefinite matrix $N_{k+1} := \Sigma_k - \Sigma_{k+1}$. Hence, one gets

$$
\begin{aligned}
\mathrm{Tr}\{F_{k+1}(\Sigma_k-\Sigma_{k+1})\} \;=\;& \mathrm{Tr}\{(F_{k+1}^{1/2})^2(N_{k+1}^{1/2})^2\} \\
=\;& \mathrm{Tr}\{N_{k+1}^{1/2}(F_{k+1}^{1/2})^2 N_{k+1}^{1/2}\} \\
=\;& \mathrm{Tr}\{(N_{k+1}^{1/2}F_{k+1}^{1/2})(N_{k+1}^{1/2}F_{k+1}^{1/2})'\} \\
\geq\;& 0 \qquad (123)
\end{aligned}
$$

for all the realizations of the random matrices involved. Similarly, for all the realizations of the random matrices involved, and for all $k = -1, 0, 1, 2, \ldots$, one obtains

$$\mathrm{Tr}\{F_{k+1}\Sigma_{k+1}\} \leq \mathrm{Tr}\{F_{k+1}\Sigma_k\} \leq \ldots \leq \mathrm{Tr}\{F_{k+1}\Sigma_{-1}\} \tag{124}$$

which is (36). ■

**Proof of Proposition 4**

In this proof, we use the superscript "$(N)$" to denote expressions obtained for the finite-horizon case with horizon $N$ and $Q_N = 0$ with probability 1, and assuming that (33) holds for the other values of $k$. Due to (106), for any sequence of feasible updates and any finite horizon $N$, one has

$$\begin{aligned}
& \mathop{\mathbb{E}}_{e_0,\{x_k\}_{k=0}^{N-1},\{\tilde{\varepsilon}_k\}_{k=0}^{N-1}} \left\{ \sum_{k=0}^{N-1} \left[ e_k' Q_k e_k + \gamma u_k' u_k \right] \right\} \\
\geq\ & \mathop{\mathbb{E}}_{\tilde{I}_0} \{ \tilde{J}_0^{\circ,(N)}(\tilde{I}_0) \} \\
=\ & \mathop{\mathbb{E}}_{\tilde{I}_0} \left\{ \mathop{\mathbb{E}}_{e_0,K_0^{(N)}} \left\{ e_0' K_0^{(N)} e_0 \middle| \tilde{I}_0 \right\} + \mathop{\mathbb{E}}_{e_0} \left\{ \left( e_0 - \mathop{\mathbb{E}}_{e_0} \left\{ e_0 \middle| \tilde{I}_0 \right\} \right)' F_0^{(N)} \left( e_0 - \mathop{\mathbb{E}}_{e_0}\{ e_0 \middle| \tilde{I}_0 \right) \middle| \tilde{I}_0 \right\} \right\} \\
& + \mathop{\mathbb{E}}_{\{C_j\}_{j=0}^{N-2}} \left\{ \sum_{h=0}^{N-2} f_h^{(N)}(\{C_j\}_{j=0}^h) \right\}.
\end{aligned} \tag{125}$$

Using steps similar to the ones made to obtain (90) and (91), and observing from (91) that the functions $f_h^{(N)}$ can be written as

$$f_h^{(N)}(\{C_j\}_{j=0}^h) = \mathop{\mathbb{E}}_{\Sigma_{h+1}} \left\{ \mathrm{Tr}\left\{ F_{h+1}^{(N)} \Sigma_{h+1} \right\} \middle| C_0, \ldots, C_h \right\},$$

one gets [27]

$$\begin{aligned}
& \mathop{\mathbb{E}}_{\tilde{I}_0} \left\{ \mathop{\mathbb{E}}_{e_0,K_0^{(N)}} \left\{ e_0' K_0^{(N)} e_0 \middle| \tilde{I}_0 \right\} + \mathop{\mathbb{E}}_{e_0} \left\{ \left( e_0 - \mathop{\mathbb{E}}_{e_0} \left\{ e_0 \middle| \tilde{I}_0 \right\} \right)' F_0^{(N)} \left( e_0 - \mathop{\mathbb{E}}_{e_0}\{ e_0 \middle| \tilde{I}_0 \right) \middle| \tilde{I}_0 \right\} \right\} \\
& + \mathop{\mathbb{E}}_{\{C_j\}_{j=0}^{N-2}} \left\{ \sum_{h=0}^{N-2} f_h^{(N)}(\{C_j\}_{j=0}^h) \right\} \\
=\ & \mathrm{Tr}\{\overline{K}_0^{(N)} \Sigma_w\} + \mathop{\mathbb{E}}_{\Sigma_0} \left\{ \mathrm{Tr}\{F_0^{(N)} \Sigma_0\} \right\} + \sum_{h=0}^{N-2} \mathop{\mathbb{E}}_{\Sigma_{h+1}} \left\{ \mathrm{Tr}\left\{ F_{h+1}^{(N)} \Sigma_{h+1} \right\} \right\} \\
=\ & \mathrm{Tr}\{\overline{K}_0^{(N)} \Sigma_w\} + \sum_{h=0}^{N-1} \mathop{\mathbb{E}}_{\Sigma_h} \left\{ \mathrm{Tr}\left\{ F_h^{(N)} \Sigma_h \right\} \right\} \\
=\ & \mathrm{Tr}\{\overline{K}_0^{(N)} \Sigma_w\} + \sum_{h=0}^{N-1} \mathrm{Tr}\left\{ F_h^{(N)} \mathop{\mathbb{E}}_{\Sigma_h} \{\Sigma_h\} \right\}.
\end{aligned} \tag{126}$$

27. Here, one applies also the law of iterated expectations, together with the fact that the matrices $\overline{K}_0^{(N)}$ and $F_h^{(N)}$ are deterministic. Moreover, $e_0 e_0'$ and $Q_0$ are independent (which justifies the appearance of $\overline{K}_0^{(N)}$ and $\Sigma_w$ in formula (126)).

Hence, from (125), one gets, for any feasible sequence of updates,

$$\mathop{\mathbb{E}}_{e_0,\{x_k\}_{k=0}^{N-1},\{\tilde{\varepsilon}_k\}_{k=0}^{N-1}} \left\{ \sum_{k=0}^{N-1} \left[ e_k' Q_k e_k + \gamma u_k' u_k \right] \right\} \geq \mathrm{Tr}\{\overline{K}_0^{(N)} \Sigma_w\} + \sum_{h=0}^{N-1} \mathrm{Tr} \left\{ F_h^{(N)} \mathop{\mathbb{E}}_{\Sigma_h} \{\Sigma_h\} \right\}, \tag{127}$$

then

$$\begin{aligned} & \liminf_{N \to +\infty} \left( \frac{1}{N} \mathop{\mathbb{E}}_{e_0,\{x_k\}_{k=0}^{N-1},\{\tilde{\varepsilon}_k\}_{k=0}^{N-1}} \left\{ \sum_{k=0}^{N-1} \left[ e_k' Q_k e_k + \gamma u_k' u_k \right] \right\} \right) \\ \geq & \liminf_{N \to +\infty} \left( \frac{1}{N} \left( \mathrm{Tr}\{\overline{K}_0^{(N)} \Sigma_w\} + \sum_{h=0}^{N-1} \mathrm{Tr} \left\{ F_h^{(N)} \mathop{\mathbb{E}}_{\Sigma_h} \{\Sigma_h\} \right\} \right) \right). \end{aligned} \tag{128}$$

Now, we show that the second "liminf" in (128) is actually a "lim", and that

$$\lim_{N \to +\infty} \left( \frac{1}{N} \left( \mathrm{Tr}\{\overline{K}_0^{(N)} \Sigma_w\} + \sum_{h=0}^{N-1} \mathrm{Tr} \left\{ F_h^{(N)} \mathop{\mathbb{E}}_{\Sigma_h} \{\Sigma_h\} \right\} \right) \right) = 0\,. \tag{129}$$

First, as $N$ tends to $+\infty$, one gets

$$\lim_{N \to +\infty} \mathrm{Tr}\{\overline{K}_0^{(N)} \Sigma_w\} = \mathrm{Tr}\{\overline{K} \Sigma_w\} \tag{130}$$

by the convergence of $\overline{K}_0^{(N)}$ to the stationary solution $\overline{K}$ of the ARE. Moreover, due to the definition of $F_h^{(N)}$, one gets

$$F_{h+1}^{(N+1)} = F_h^{(N)} = \ldots = F_0^{(N-h)}\,, \tag{131}$$

for every $N$ and $h = 0, \ldots, N-1$. From (131) and the convergence of $F_0^{(N-h)}$ to $\overline{F}$ as $N$ tends to $+\infty$ (for fixed $h$), it follows that, for every $\delta > 0$ there exists $t(\delta) \in \mathbb{N}_0$ such that, for every $N > t(\delta)$, one gets

$$-\delta I \preceq F_h^{(N)} - \overline{F} \preceq \delta I$$

for every $h = 0, \ldots, N - t(\delta) - 1$, whereas

$$\begin{aligned} F_{N-t(\delta)}^{(N)} &= F_0^{(t(\delta))}\,, \\ F_{N-t(\delta)+1}^{(N)} &= F_0^{(t(\delta)-1)}\,, \\ &\ldots\,, \\ F_{N-1}^{(N)} &= F_0^1 \end{aligned}$$

are a finite number $t(\delta)$ of fixed (i.e., independent from $N$) symmetric and positive-semidefinite matrices. Hence, recalling (21), (30), (36), and (123), and the linearity of the trace operator,

one obtains

$$
\begin{aligned}
0 &\leq \sum_{h=0}^{N-1} \mathrm{Tr}\left\{F_h^{(N)} \mathop{\mathbb{E}}_{\Sigma_h}\{\Sigma_h\}\right\} \\
&\leq \sum_{h=0}^{N-t(\delta)-1} \mathrm{Tr}\left\{(\overline{F}+\delta I)\mathop{\mathbb{E}}_{\Sigma_h}\{\Sigma_h\}\right\} + \sum_{h=N-t(\delta)}^{N-1} \mathrm{Tr}\left\{F_h^{(N)} \mathop{\mathbb{E}}_{\Sigma_h}\{\Sigma_h\}\right\} \\
&\leq (N-t(\delta))\delta \mathrm{Tr}\left\{\mathop{\mathbb{E}}_{\Sigma_{-1}}\{\Sigma_{-1}\}\right\} + \sum_{h=0}^{N-t(\delta)-1} \mathrm{Tr}\left\{\overline{F}\mathop{\mathbb{E}}_{\Sigma_h}\{\Sigma_h\}\right\} + \sum_{t=1}^{t(\delta)} \mathrm{Tr}\left\{F_0^{(t)} \mathop{\mathbb{E}}_{\Sigma_{-1}}\{\Sigma_{-1}\}\right\} \\
&= (N-t(\delta))\delta \mathrm{Tr}\{\Sigma_w\} + \sum_{h=0}^{N-t(\delta)-1} \mathrm{Tr}\left\{\overline{F}\mathop{\mathbb{E}}_{\Sigma_h}\{\Sigma_h\}\right\} + \sum_{t=1}^{t(\delta)} \mathrm{Tr}\left\{F_0^{(t)}\Sigma_w\right\}
\end{aligned}
$$

This, combined with the finiteness of $\mathrm{Tr}\left\{\overline{F}\mathop{\mathbb{E}}_{\Sigma_0}\{\Sigma_0\}\right\}$ and

$$
\mathrm{Tr}\left\{\overline{F}\mathop{\mathbb{E}}_{\Sigma_{h+1}}\{\Sigma_{h+1}\}\right\} \leq \mathrm{Tr}\left\{\overline{F}\mathop{\mathbb{E}}_{\Sigma_h}\{\Sigma_h\}\right\}
$$

for every $h = 0, 1, \ldots$, shows that

$$
\begin{aligned}
0 &\leq \liminf_{N\to+\infty}\left(\frac{1}{N}\sum_{h=0}^{N-1} \mathrm{Tr}\left\{F_h^{(N)} \mathop{\mathbb{E}}_{\Sigma_h}\{\Sigma_h\}\right\}\right) \\
&\leq \liminf_{N\to+\infty}\left(\mathrm{Tr}\left\{\overline{F}\mathop{\mathbb{E}}_{\Sigma_N}\{\Sigma_N\}\right\}\right) + \delta \mathrm{Tr}\{\Sigma_w\}
\end{aligned}
$$

(see also the derivation of (123) for a similar proof). Moreover, since this holds for every $\delta > 0$, one obtains

$$
\begin{aligned}
0 &\leq \liminf_{N\to+\infty}\left(\frac{1}{N}\sum_{h=0}^{N-1} \mathrm{Tr}\left\{F_h^{(N)} \mathop{\mathbb{E}}_{\Sigma_h}\{\Sigma_h\}\right\}\right) \\
&\leq \liminf_{N\to+\infty}\left(\mathrm{Tr}\left\{\overline{F}\mathop{\mathbb{E}}_{\Sigma_N}\{\Sigma_N\}\right\}\right) \\
&= \mathrm{Tr}\left\{\overline{F}\liminf_{N\to+\infty}\mathop{\mathbb{E}}_{\Sigma_N}\{\Sigma_N\}\right\}, \\
&= \mathrm{Tr}\left\{\overline{F}\lim_{N\to+\infty}\mathop{\mathbb{E}}_{\Sigma_N}\{\Sigma_N\}\right\}, \\
&= 0, \qquad (132)
\end{aligned}
$$

where the last three steps are due to the linearity of the trace operator, to (32), and to (34). This, combined with (130), proves (129). This means that the infimum of the average learning functional (41) to be optimized is 0.

Now, we show that the sequence of updating functions (42) minimizes the average learning functional (41) with respect to all feasible sequences of updates. To see this, we use the superscript "$(N,\overline{K})$" to denote expressions obtained for the finite-horizon case with horizon $N$ and with $\overline{K}_N = \overline{K}$ (i.e., assuming that $\overline{K}_N$ is equal to its stationary value $\overline{K}$), and assuming that (33) holds for the other values of $k$. Then, due to (106), when the updates are generated by the sequence of updating functions (42), one has, for any finite horizon $N$,

$$
\begin{aligned}
& \mathop{\mathbb{E}}_{e_0,\{x_k\}_{k=0}^{N-1},\{\tilde{\varepsilon}_k\}_{k=0}^{N-1}} \left\{ \sum_{k=0}^{N-1} [e_k' Q_k e_k + \gamma u_k' u_k] \right\} \\
\leq \; & \mathop{\mathbb{E}}_{e_0,\{x_k\}_{k=0}^{N},\{\tilde{\varepsilon}_k\}_{k=0}^{N-1}} \left\{ \sum_{k=0}^{N-1} [e_k' Q_k e_k + \gamma u_k' u_k] + e_N' \overline{K} e_N \right\} \\
\leq \; & \mathop{\mathbb{E}}_{\tilde{I}_0} \{ \tilde{J}_0^{\circ,(N,\overline{K})}(\tilde{I}_0) \} \\
= \; & \mathop{\mathbb{E}}_{\tilde{I}_0} \left\{ \mathop{\mathbb{E}}_{e_0} \left\{ e_0' \overline{K} e_0 \middle| \tilde{I}_0 \right\} + \mathop{\mathbb{E}}_{e_0} \left\{ \left( e_0 - \mathop{\mathbb{E}}_{e_0} \left\{ e_0 \middle| \tilde{I}_0 \right\} \right)' \overline{F} \left( e_0 - \mathop{\mathbb{E}}_{e_0} \left\{ e_0 \middle| \tilde{I}_0 \right\} \right) \middle| \tilde{I}_0 \right\} \right\} \\
& + \mathop{\mathbb{E}}_{\{C_j\}_{j=0}^{N-2}} \left\{ \sum_{h=0}^{N-2} f_h^{(N,\overline{K})}(\{C_j\}_{j=0}^{h}) \right\} . \qquad (133)
\end{aligned}
$$

Moreover, likewise in the derivations of (126)-(132), one gets

$$
\begin{aligned}
& \mathop{\mathbb{E}}_{\tilde{I}_0} \left\{ \mathop{\mathbb{E}}_{e_0} \left\{ e_0' \overline{K} e_0 \middle| \tilde{I}_0 \right\} + \mathop{\mathbb{E}}_{e_0} \left\{ \left( e_0 - \mathop{\mathbb{E}}_{e_0} \left\{ e_0 \middle| \tilde{I}_0 \right\} \right)' \overline{F} \left( e_0 - \mathop{\mathbb{E}}_{e_0} \left\{ e_0 \middle| \tilde{I}_0 \right\} \right) \middle| \tilde{I}_0 \right\} \right\} \\
& + \mathop{\mathbb{E}}_{\{C_j\}_{j=0}^{N-2}} \left\{ \sum_{h=0}^{N-2} f_h^{(N,\overline{K})}(\{C_j\}_{j=0}^{h}) \right\} . \\
= \; & \mathrm{Tr}\{\overline{K}\Sigma_w\} + \mathop{\mathbb{E}}_{\Sigma_0} \left\{ \mathrm{Tr}\{\overline{F}\Sigma_0\} \right\} + \sum_{h=0}^{N-2} \mathop{\mathbb{E}}_{\Sigma_{h+1}} \left\{ \mathrm{Tr}\left\{ \overline{F}\Sigma_{h+1} \right\} \right\} \\
= \; & \mathrm{Tr}\{\overline{K}\Sigma_w\} + \sum_{h=0}^{N-1} \mathop{\mathbb{E}}_{\Sigma_h} \left\{ \mathrm{Tr}\left\{ \overline{F}\Sigma_h \right\} \right\} \\
= \; & \mathrm{Tr}\{\overline{K}\Sigma_w\} + \sum_{h=0}^{N-1} \mathrm{Tr} \left\{ \overline{F} \mathop{\mathbb{E}}_{\Sigma_h} \{\Sigma_h\} \right\} , \qquad (134)
\end{aligned}
$$

and

$$
\begin{aligned}
0 &\leq \liminf_{N\to+\infty}\left(\frac{1}{N}\left(\mathrm{Tr}\{\overline{K}\Sigma_w\}+\sum_{h=0}^{N-1}\mathrm{Tr}\left\{\overline{F}\,\mathop{\mathbb{E}}_{\Sigma_h}\{\Sigma_h\}\right\}\right)\right)\\
&= \lim_{N\to+\infty}\left(\frac{1}{N}\left(\mathrm{Tr}\{\overline{K}\Sigma_w\}+\sum_{h=0}^{N-1}\mathrm{Tr}\left\{\overline{F}\,\mathop{\mathbb{E}}_{\Sigma_h}\{\Sigma_h\}\right\}\right)\right)\\
&= \lim_{N\to+\infty}\left(\frac{1}{N}\sum_{h=0}^{N-1}\mathrm{Tr}\left\{\overline{F}\,\mathop{\mathbb{E}}_{\Sigma_h}\{\Sigma_h\}\right\}\right)\\
&\leq \lim_{N\to+\infty}\mathrm{Tr}\left\{\overline{F}\,\mathop{\mathbb{E}}_{\Sigma_N}\{\Sigma_N\}\right\}\\
&= 0\,. \qquad (135)
\end{aligned}
$$

This, combined with (129), (133), and (134), proves the optimality of the sequence of updating functions (42) with respect to the average learning functional (41) and all sequences of feasible updating functions. ■

**Proof of Proposition 5**

(i) One has

$$
\begin{aligned}
\mathrm{MSE}_k^{\dagger} &= \mathrm{Tr}\left\{\mathop{\mathbb{E}}_{w,\hat{w}_k^{\dagger}}\left\{\left(w-\hat{w}_k^{\dagger}\right)'\left(w-\hat{w}_k^{\dagger}\right)\right\}\right\}\\
&= \mathrm{Tr}\left\{\mathop{\mathbb{E}}_{w,\hat{w}_k^{\dagger}}\left\{\left(w-\hat{w}_k^{\dagger}\right)\left(w-\hat{w}_k^{\dagger}\right)'\right\}\right\}\\
&= \mathrm{Tr}\left\{\mathop{\mathbb{E}}_{w_k,\tilde{I}_k}\left\{\left(w_k-\mathop{\mathbb{E}}_{w_k}\left\{w_k\middle|I_k\right\}\right)\left(w_k-\mathop{\mathbb{E}}_{w_k}\left\{w_k\middle|I_k\right\}\right)'\right\}\right\}\\
&= \mathrm{Tr}\left\{\mathop{\mathbb{E}}_{I_k}\left\{\mathop{\mathbb{E}}_{w_k}\left\{\left(w_k-\mathop{\mathbb{E}}_{w_k}\left\{w_k\middle|I_k\right\}\right)\left(w_k-\mathop{\mathbb{E}}_{w_k}\left\{w_k\middle|I_k\right\}\right)'\middle|\tilde{I}_k\right\}\right\}\right\}\\
&= \mathrm{Tr}\left\{\mathop{\mathbb{E}}_{\tilde{I}_k}\left\{\mathop{\mathbb{E}}_{e_k}\left\{\left(-\left(e_k-\mathop{\mathbb{E}}_{e_k}\left\{e_k\middle|\tilde{I}_k\right\}\right)\right)\left(-\left(e_k-\mathop{\mathbb{E}}_{e_k}\left\{e_k\middle|\tilde{I}_k\right\}\right)\right)'\middle|\tilde{I}_k\right\}\right\}\right\}\\
&= \mathrm{Tr}\left\{\mathop{\mathbb{E}}_{\Sigma_k}\{\Sigma_k\}\right\}\,. \qquad (136)
\end{aligned}
$$

(ii) To obtain the rate of convergence to 0 of the mean-square error of the KF estimate of $w$ at the time $k$, we observe that, using similar steps as in the derivation of (117) and (121), one gets

$$
\mathrm{Tr}\{\mathop{\mathbb{E}}_{\Sigma_k}\{\Sigma_k\}\}-\mathrm{Tr}\{\mathop{\mathbb{E}}_{\Sigma_{k+1}}\{\Sigma_{k+1}\}\}\geq(c_1+\sigma_\varepsilon^2)^{-1}\mathrm{Tr}\{\mathop{\mathbb{E}}_{\Sigma_k}\{\Sigma_k^2\}\}\lambda_{\min}(\overline{Q})\,. \qquad (137)
$$

By iterating (137), we obtain

$$\mathrm{Tr}\{\mathop{\mathbb{E}}_{\Sigma_0}\{\Sigma_0\}\} - \mathrm{Tr}\{\mathop{\mathbb{E}}_{\Sigma_{k+1}}\{\Sigma_{k+1}\}\} \geq (c_1 + \sigma_\varepsilon^2)^{-1} \left( \sum_{j=0}^{k} \mathrm{Tr}\{\mathop{\mathbb{E}}_{\Sigma_j}\{\Sigma_j^2\}\} \right) \lambda_{\min}(\overline{Q})\,. \tag{138}$$

This, combined with equation (31), provides

$$\begin{aligned} \mathrm{Tr}\{\mathop{\mathbb{E}}_{\Sigma_0}\{\Sigma_0\}\} - \mathrm{Tr}\{\mathop{\mathbb{E}}_{\Sigma_{k+1}}\{\Sigma_{k+1}\}\} &\geq (c_1 + \sigma_\varepsilon^2)^{-1}(k+1)\,\mathrm{Tr}\{\mathop{\mathbb{E}}_{\Sigma_k}\{\Sigma_k^2\}\}\lambda_{\min}(\overline{Q}) \\ &\geq (c_1 + \sigma_\varepsilon^2)^{-1}(k+1)\,\mathrm{Tr}\{\mathop{\mathbb{E}}_{\Sigma_{k+1}}\{\Sigma_{k+1}^2\}\}\lambda_{\min}(\overline{Q})\,. \end{aligned} \tag{139}$$

Now, we apply the property (118) with $S_1 = \Sigma_{k+1}$ and $S_2 = I$, obtaining

$$|\mathrm{Tr}\{\Sigma_{k+1}\}| \leq \sqrt{\mathrm{Tr}\{\Sigma_{k+1}^2\}d}\,,$$

hence

$$(\mathrm{Tr}\{\Sigma_{k+1}\})^2 \leq \mathrm{Tr}\{\Sigma_{k+1}^2\}d$$

and

$$\mathop{\mathbb{E}}_{\Sigma_{k+1}} \left\{ (\mathrm{Tr}\{\Sigma_{k+1}\})^2 \right\} \leq \mathop{\mathbb{E}}_{\Sigma_{k+1}} \left\{ \mathrm{Tr}\{\Sigma_{k+1}^2\} \right\} d = \mathrm{Tr} \left\{ \mathop{\mathbb{E}}_{\Sigma_{k+1}} \{\Sigma_{k+1}^2\} \right\} d\,. \tag{140}$$

Moreover, by the convexity of the square function $(\cdot)^2$ and Jensen's inequality, one gets

$$\mathop{\mathbb{E}}_{\Sigma_{k+1}} \left\{ (\mathrm{Tr}\{\Sigma_{k+1}\})^2 \right\} \geq \left( \mathop{\mathbb{E}}_{\Sigma_{k+1}} \{\mathrm{Tr}\{\Sigma_{k+1}\}\} \right)^2 = \left( \mathrm{Tr} \left\{ \mathop{\mathbb{E}}_{\Sigma_{k+1}} \{\Sigma_{k+1}\} \right\} \right)^2 . \tag{141}$$

Then, combining equations (139), (140), and (141), one obtains

$$\mathrm{Tr}\{\mathop{\mathbb{E}}_{\Sigma_0}\{\Sigma_0\}\} - \mathrm{Tr}\{\mathop{\mathbb{E}}_{\Sigma_{k+1}}\{\Sigma_{k+1}\}\} \geq (c_1 + \sigma_\varepsilon^2)^{-1}(k+1) \left( \mathrm{Tr} \left\{ \mathop{\mathbb{E}}_{\Sigma_{k+1}} \{\Sigma_{k+1}\} \right\} \right)^2 d^{-1}\lambda_{\min}(\overline{Q})\,. \tag{142}$$

Now, for a given tolerance $\eta > 0$, we use (142) to find an upper bound, as a function of $\eta$, on the maximal value $(k+1)(\eta)$ for $k+1$ for which $\mathrm{Tr}\{\mathop{\mathbb{E}}_{\Sigma_{k+1}}\{\Sigma_{k+1}\}\} \geq \eta$. Since, of course, one has

$$\mathrm{Tr}\{\mathop{\mathbb{E}}_{\Sigma_0}\{\Sigma_0\}\} \geq \mathrm{Tr}\{\mathop{\mathbb{E}}_{\Sigma_0}\{\Sigma_0\}\} - \mathrm{Tr}\{\mathop{\mathbb{E}}_{\Sigma_{k+1}}\{\Sigma_{k+1}\}\}\,, \tag{143}$$

by combining the inequalities (142) and (143), and the definition of $(k+1)(\eta)$, one obtains

$$\begin{aligned} \mathrm{Tr}\{\mathop{\mathbb{E}}_{\Sigma_0}\{\Sigma_0\}\} &\geq (c_1 + \sigma_\varepsilon^2)^{-1}(k+1)(\eta) \left( \mathrm{Tr} \left\{ \mathop{\mathbb{E}}_{\Sigma_{(k+1)(\eta)}} \{\Sigma_{(k+1)(\eta)}\} \right\} \right)^2 d^{-1}\lambda_{\min}(\overline{Q}) \\ &\geq (c_1 + \sigma_\varepsilon^2)^{-1}(k+1)(\eta)\,\eta^2 d^{-1}\lambda_{\min}(\overline{Q}) \\ &\geq 0\,. \end{aligned} \tag{144}$$

Now, equation (144) can hold only if

$$(k+1)(\eta) \leq \frac{(c_1 + \sigma_\varepsilon^2)\, d\, \mathrm{Tr}\{\mathop{\mathbb{E}}_{\Sigma_0}\{\Sigma_0\}\}}{\eta^2 \lambda_{\min}(\overline{Q})}\,. \tag{145}$$

Hence, renaming the index $k+1$ still by $k$ ($k = 1, 2, \ldots$), one obtains

$$k(\eta) \leq \frac{(c_1 + \sigma_\varepsilon^2)\, d\, \mathrm{Tr}\{\underset{\Sigma_0}{\mathbb{E}}\{\Sigma_0\}\}}{\eta^2 \lambda_{\min}(\overline{Q})} \,. \tag{146}$$

Similarly, for $k = 1, 2, \ldots$, denoting by $\eta(k)$ the maximal value of $\eta > 0$ for which $\mathrm{Tr}\{\underset{\Sigma_k}{\mathbb{E}}\{\Sigma_k\}\} \geq \eta$, one obtains

$$\eta(k) \leq \sqrt{\frac{(c_1 + \sigma_\varepsilon^2)\, d\, \mathrm{Tr}\{\underset{\Sigma_0}{\mathbb{E}}\{\Sigma_0\}\}}{k \lambda_{\min}(\overline{Q})}} \,, \tag{147}$$

which provides the desired rate of convergence (47). Finally, the last part of (ii) follows from such a rate of convergence, and from the definition of $\mathrm{MSE}_k^\dagger$.

(iii) By (24), (32), (34), and the assumed (uniform on $k$) almost-sure boundedness of $x_k = C_k'$, we get

$$\begin{aligned}
\lim_{k \to +\infty} \underset{H_k}{\mathbb{E}}\{H_k\} &= \lim_{(k+1) \to +\infty} \underset{H_{k+1}}{\mathbb{E}}\{H_{k+1}\} \\
&= \lim_{k \to +\infty} \underset{\Sigma_{k+1}, C_{k+1}}{\mathbb{E}}\{\Sigma_{k+1} C_{k+1}' (\sigma_\varepsilon^2)^{-1}\} \\
&= \lim_{k \to +\infty} \underset{\Sigma_{k+1}, C_{k+1}}{\mathbb{E}}\{\Sigma_{k+1} C_{k+1}'\} (\sigma_\varepsilon^2)^{-1} \\
&= 0 \cdot (\sigma_\varepsilon^2)^{-1} \\
&= 0 \,.
\end{aligned}$$

(iv) An application of (47), of the assumed (uniform on $k$) almost-sure boundedness of $x_k$, of [20, Theorem 1], and of the matrix version of Markov's inequality[28], shows that, for every $\delta > 0$, one has

$$\lim_{k \to +\infty} \Pr\{|H_{k,(h,l)}| > \delta\} = 0 \,.$$

■

28. The matrix version of Markov's inequality [39, Theorem A.1] states that, for any symmetric and positive-semidefinite random matrix $X$, and any fixed symmetric and positive-definite matrix $M$ of the same dimension, one has

$$\Pr\{X \not\preceq M\} \leq \mathrm{Tr}\{\underset{X}{\mathbb{E}}\{X\} M^{-1}\} \,. \tag{148}$$

Moreover, by [20, Theorem 2], one has also

$$\mathrm{Tr}\{\underset{X}{\mathbb{E}}\{X\} M^{-1}\} \leq \mathrm{Tr}\{\underset{X}{\mathbb{E}}\{X\}\} \frac{1}{\lambda_{\max}(M)} \,, \tag{149}$$

where $\lambda_{\max}(M)$ is the largest eigenvalue of $M$. Then, one applies (148) and (149) with $X$ replaced by $\Sigma_k$, and $M$ by a sequence of positive-definite matrices $M_j$ of the same dimension, with $\lambda_{\max}(M_j)$ decreasing and tending to 0, when $j$ tends to $+\infty$.

**Proof of Proposition 6**

(i) First, we observe that

$$\mathop{\mathbb{E}}_{e_k^\dagger}\{e_k^\dagger\} = \mathop{\mathbb{E}}_{w_k^\dagger, w}\{\hat{w}_k^\dagger - w\} = \mathop{\mathbb{E}}_{w, I_k}\left\{\mathop{\mathbb{E}}_{w_k}\left\{w_k \middle| I_k\right\} - w\right\} = \mathop{\mathbb{E}}_{w_k}\{w_k\} - \mathop{\mathbb{E}}_{w}\{w\} = \mathop{\mathbb{E}}_{w}\{w\} - \mathop{\mathbb{E}}_{w}\{w\} = 0\,. \tag{150}$$

Proceeding as in the proof of (136), we get

$$\Sigma_{e_k^\dagger} = \mathop{\mathbb{E}}_{\Sigma_k}\left\{\Sigma_k\right\}\,, \tag{151}$$

Similarly, recalling that

$$e_k^\circ := \hat{w}_k^\circ - w\,,$$

and using (22) and (150), we get

$$\begin{aligned}
\mathop{\mathbb{E}}_{e_{k+1}^\circ}\left\{e_{k+1}^\circ\right\} &= \mathop{\mathbb{E}}_{\hat{w}_{k+1}^\circ, w}\{\hat{w}_{k+1}^\circ - w\} \\
&= \mathop{\mathbb{E}}_{\hat{w}_k^\circ}\{\hat{w}_k^\circ\} + L_k \mathop{\mathbb{E}}_{\hat{w}_k^\circ, \hat{w}_k^\dagger}\{\hat{w}_k^\circ - \hat{w}_k^\dagger\} - \mathop{\mathbb{E}}_{w}\{w\} \\
&= \mathop{\mathbb{E}}_{\hat{w}_k^\circ}\{\hat{w}_k^\circ\} + L_k \mathop{\mathbb{E}}_{\hat{w}_k^\circ, w}\{\hat{w}_k^\circ - w\} - \mathop{\mathbb{E}}_{w}\{w\} \\
&= (I + L_k) \mathop{\mathbb{E}}_{\hat{w}_k^\circ, w}\{\hat{w}_k^\circ - w\} \\
&= (I + L_k) \mathop{\mathbb{E}}_{\hat{e}_k^\circ}\{\hat{e}_k^\circ\} \\
&= \Pi_{j=0}^k (I + L_j) \mathop{\mathbb{E}}_{\hat{w}_0^\circ, w}\{\hat{w}_0^\circ - w\} \\
&= \Pi_{j=0}^k (I + L_j) \mathop{\mathbb{E}}_{\hat{w}_0^\circ, w}\{\hat{w}_0 - w\} \\
&= 0\,.
\end{aligned} \tag{152}$$

Let

$$\begin{aligned}
\Sigma_{e_k^\circ, e_k^\dagger} &:= \mathop{\mathbb{E}}_{e_k^\circ, e_k^\dagger}\left\{\left(e_k^\circ - \mathop{\mathbb{E}}_{e_k^\circ}\{e_k^\circ\}\right)\left(e_k^\dagger - \mathop{\mathbb{E}}_{e_k^\dagger}\{e_k^\dagger\}\right)'\right\} \\
&= \mathop{\mathbb{E}}_{e_k^\circ, e_k^\dagger}\left\{\left(e_k^\circ\right)\left(e_k^\dagger\right)'\right\}
\end{aligned} \tag{153}$$

and denote by

$$\begin{aligned}
\Sigma_{e_k^\dagger, e_k^\circ} &:= \mathop{\mathbb{E}}_{e_k^\dagger, e_k^\circ}\left\{\left(e_k^\dagger - \mathop{\mathbb{E}}_{e_k^\dagger}\{e_k^\dagger\}\right)\left(e_k^\circ - \mathop{\mathbb{E}}_{e_k^\circ}\{e_k^\circ\}\right)'\right\} \\
&= \mathop{\mathbb{E}}_{e_k^\dagger, e_k^\circ}\left\{\left(e_k^\dagger\right)\left(e_k^\circ\right)'\right\}
\end{aligned} \tag{154}$$

and

$$\Sigma_{e_k^\circ, e_k^\dagger} \quad = \quad \Sigma'_{e_k^\dagger, e_k^\circ} \tag{155}$$

the two (unconditional) cross-covariance matrices of $e_k^\circ$ and $e_k^\dagger$. Finally, we consider the vector

$$E_k := ((e_k^\circ)', (e_k^\dagger)')' ,$$

and we denote by

$$\Sigma_{E_k} \quad = \quad \begin{pmatrix} \Sigma_{e_k^\circ} & \Sigma_{e_k^\circ, e_k^\dagger} \\ \Sigma_{e_k^\dagger, e_k^\circ} & \Sigma_{e_k^\dagger} \end{pmatrix}$$

its (unconditional) covariance matrix. As they are needed in the following analysis, we also provide the following upper bounds on the traces of the two matrices $(I + L_k)\Sigma_{e_k^\circ, e_k^\dagger} L'_k$ and $L_k \Sigma_{e_k^\dagger, e_k^\circ}(I + L_k)'$:

$$\begin{aligned}
& \mathrm{Tr}\{(I + L_k)\Sigma_{e_k^\circ, e_k^\dagger} L'_k\} \\
= \; & \mathrm{Tr}\left\{ \underset{e_k^\circ, e_k^\dagger}{\mathbb{E}} \left\{ ((I + L_k)e_k^\circ) \left( L_k e_k^\dagger \right)' \right\} \right\} \\
= \; & \mathrm{Tr}\left\{ \underset{e_k^\dagger, e_k^\circ}{\mathbb{E}} \left\{ \left( L_k e_k^\dagger \right)' ((I + L_k)e_k^\circ) \right\} \right\} \\
= \; & \underset{e_k^\dagger, e_k^\circ}{\mathbb{E}} \left\{ \left( L_k e_k^\dagger \right)' ((I + L_k)e_k^\circ) \right\} \\
\leq \; & \sqrt{ \underset{e_k^\dagger}{\mathbb{E}} \left\{ \left( L_k e_k^\dagger \right)' \left( L_k e_k^\dagger \right) \right\} \underset{e_k^\circ}{\mathbb{E}} \left\{ ((I + L_k)e_k^\circ)' ((I + L_k)e_k^\circ) \right\} } \\
= \; & \sqrt{ \mathrm{Tr}\left\{ \underset{e_k^\dagger}{\mathbb{E}} \left\{ \left( L_k e_k^\dagger \right)' \left( L_k e_k^\dagger \right) \right\} \right\} \mathrm{Tr}\left\{ \underset{e_k^\circ}{\mathbb{E}} \left\{ ((I + L_k)e_k^\circ)' ((I + L_k)e_k^\circ) \right\} \right\} } \\
= \; & \sqrt{ \mathrm{Tr}\left\{ \underset{e_k^\dagger}{\mathbb{E}} \left\{ \left( L_k e_k^\dagger \right) \left( L_k e_k^\dagger \right)' \right\} \right\} \mathrm{Tr}\left\{ \underset{e_k^\circ}{\mathbb{E}} \left\{ ((I + L_k)e_k^\circ) ((I + L_k)e_k^\circ)' \right\} \right\} } \\
= \; & \sqrt{ \mathrm{Tr}\{L_k \Sigma_{e_k^\dagger} L'_k\} \mathrm{Tr}\{(I + L_k)\Sigma_{e_k^\circ}(I + L_k)'\} } \, ,
\end{aligned} \tag{156}$$

and similarly,

$$\begin{aligned}
\mathrm{Tr}\{L_k \Sigma_{e_k^\dagger, e_k^\circ}(I + L_k)'\} \quad & = \quad \mathrm{Tr}\{ \left( (I + L_k)\Sigma_{e_k^\circ, e_k^\dagger} L'_k \right)' \} \\
& \leq \quad \sqrt{ \mathrm{Tr}\{L_k \Sigma_{e_k^\dagger} L'_k\} \mathrm{Tr}\{(I + L_k)\Sigma_{e_k^\circ}(I + L_k)'\} } \, .
\end{aligned} \tag{157}$$

The mean-square error of the OLL estimate $\hat{w}_k^\circ$ of $w$ at the time $k$ is given by

$$\begin{aligned}
\text{MSE}_k^\circ &= \mathop{\mathbb{E}}_{w,\hat{w}_k^\circ} \left\{ (w-\hat{w}_k^\circ)'\,(w-\hat{w}_k^\circ) \right\} \\
&= \mathop{\mathbb{E}}_{e_k^\circ} \left\{ (-e_k^\circ)'\,(-e_k^\circ) \right\} = \mathop{\mathbb{E}}_{e_k^\circ} \left\{ (e_k^\circ)'\,(e_k^\circ) \right\} = \text{Tr}\left\{ \mathop{\mathbb{E}}_{e_k^\circ} \left\{ (e_k^\circ)'\,(e_k^\circ) \right\} \right\} = \text{Tr}\left\{ \mathop{\mathbb{E}}_{e_k^\circ} \left\{ (e_k^\circ)\,(e_k^\circ)' \right\} \right\} \\
&= \text{Tr}\left\{ \Sigma_{e_k^\circ} \right\}, \qquad (158)
\end{aligned}$$

which proves (i).

(ii) By (22), one has the following recursion:

$$\begin{aligned}
e_{k+1}^\circ &= e_k^\circ + L_k(e_k^\circ - e_k^\dagger)\,, \\
&= \begin{pmatrix} I+L_k \\ -L_k \end{pmatrix} E_k\,. \qquad (159)
\end{aligned}$$

Hence, one gets

$$\begin{aligned}
&\Sigma_{e_{k+1}^\circ} \\
=\;& \begin{pmatrix} I+L_k \\ -L_k \end{pmatrix} \Sigma_{E_k} \begin{pmatrix} I+L_k \\ -L_k \end{pmatrix}' \\
=\;& \begin{pmatrix} I+L_k \\ -L_k \end{pmatrix} \begin{pmatrix} \Sigma_{e_k^\circ} & \Sigma_{e_k^\circ, e_k^\dagger} \\ \Sigma_{e_k^\dagger, e_k^\circ} & \Sigma_{e_k^\dagger} \end{pmatrix} \begin{pmatrix} I+L_k \\ -L_k \end{pmatrix}' \\
=\;& \begin{pmatrix} I+L_k \\ -L_k \end{pmatrix} \begin{pmatrix} \Sigma_{e_k^\circ}(I+L_k)' - \Sigma_{e_k^\circ, e_k^\dagger} L_k' \\ \Sigma_{e_k^\dagger, e_k^\circ}(I+L_k)' - \Sigma_{e_k^\dagger} L_k' \end{pmatrix} \\
=\;& (I+L_k)\Sigma_{e_k^\circ}(I+L_k)' - (I+L_k)\Sigma_{e_k^\circ, e_k^\dagger} L_k' - L_k \Sigma_{e_k^\dagger, e_k^\circ}(I+L_k)' + L_k \Sigma_{e_k^\dagger} L_k'\,.
\end{aligned}$$

Then, by using equations (156) and (157), we get

$$\begin{aligned}
\text{Tr}\left\{ \Sigma_{e_{k+1}^\circ} \right\} \;\leq\;& \text{Tr}\left\{ (I+L_k)\Sigma_{e_k^\circ}(I+L_k)' \right\} \\
&+\text{Tr}\left\{ L_k \Sigma_{e_k^\dagger} L_k' \right\} \\
&+2\sqrt{\text{Tr}\left\{ (I+L_k)\Sigma_{e_k^\circ}(I+L_k)' \right\} \text{Tr}\left\{ L_k \Sigma_{e_k^\dagger} L_k' \right\}}\,, \qquad (160)
\end{aligned}$$

which proves (ii).

(iii) Let us now consider the infinite-horizon case. Then, $L_k$ is replaced by $L$, and (160) becomes

$$\begin{aligned}
\text{Tr}\left\{ \Sigma_{e_{k+1}^\circ} \right\} \;\leq\;& \text{Tr}\left\{ (I+L)\Sigma_{e_k^\circ}(I+L)' \right\} \\
&+\text{Tr}\left\{ L \Sigma_{e_k^\dagger} L' \right\} \\
&+2\sqrt{\text{Tr}\left\{ (I+L)\Sigma_{e_k}^\circ(I+L)' \right\} \text{Tr}\left\{ L \Sigma_{e_k^\dagger} L' \right\}}\,, \qquad (161)
\end{aligned}$$

where the matrix $(I+L)$ has all its eigenvalues inside the unit circle by [7, Section 4.1, Proposition 4.1], and $\Sigma_{e_k^\dagger}$ tends to the 0 matrix when $k$ tends to $+\infty$, due to (32), (34), and (151). We first show that the non-negative sequence

$$\left\{\mathrm{Tr}\left\{\Sigma_{e_k^\circ}\right\}, k=0,1,\ldots\right\} \tag{162}$$

is bounded. Indeed, let $M_0$ be a symmetric and positive-semidefinite $d\times d$ matrix such that

$$\mathrm{Tr}\left\{M_0\right\} \geq \mathrm{Tr}\left\{\Sigma_{e_0}^\circ\right\}\ . \tag{163}$$

We first note that the matrix $L$ is symmetric by its definition (39), since the symmetric matrices $(\overline{K}+I)^{-1}$ and $\overline{K}$ commute, being associated with the same basis of eigenvectors. Then, by [20, Theorem 1], one has

$$\mathrm{Tr}\left\{(I+L)\Sigma_{e_0^\circ}(I+L)'\right\} \leq (|\lambda|_{\max}(I+L)|)^2\,\mathrm{Tr}\left\{\Sigma_{e_0^\circ}\right\} < \mathrm{Tr}\left\{\Sigma_{e_0^\circ}\right\}\ , \tag{164}$$

since the spectral radius $|\lambda|_{\max}(I+L)<1$, and

$$\mathrm{Tr}\left\{L\Sigma_{e_0^\dagger}L'\right\} \leq (|\lambda|_{\max}(L))^2\,\mathrm{Tr}\left\{\Sigma_{e_0^\dagger}\right\}\ . \tag{165}$$

Combining the inequalities (164) and (165) with (161), we obtain

$$\begin{aligned}
\mathrm{Tr}\left\{\Sigma_{e_1^\circ}\right\} \leq\ & (|\lambda|_{\max}(I+L))^2\,\mathrm{Tr}\left\{\Sigma_{e_0^\circ}\right\} \\
& +|\lambda|_{\max}(L)^2\,\mathrm{Tr}\left\{\Sigma_{e_0^\dagger}\right\} \\
& +2\sqrt{(|\lambda|_{\max}(I+L))^2\,\mathrm{Tr}\left\{\Sigma_{e_0^\circ}\right\}(|\lambda|_{\max}(L))^2\,\mathrm{Tr}\left\{\Sigma_{e_0^\dagger}\right\}} \\
\leq\ & (|\lambda|_{\max}(I+L))^2\,\mathrm{Tr}\left\{M_0\right\} \\
& +(|\lambda|_{\max}(L))^2\,\mathrm{Tr}\left\{\Sigma_{e_0^\dagger}\right\} \\
& +2\sqrt{(|\lambda|_{\max}(I+L))^2\,\mathrm{Tr}\left\{M_0\right\}(|\lambda|_{\max}(L))^2\,\mathrm{Tr}\left\{\Sigma_{e_0^\dagger}\right\}}\,.
\end{aligned} \tag{166}$$

Moreover, if $\mathrm{Tr}\left\{M_0\right\}$ is sufficiently large, one gets[29]

$$\begin{aligned}
& (|\lambda|_{\max}(I+L))^2\,\mathrm{Tr}\left\{M_0\right\}+(|\lambda|_{\max}(L))^2\,\mathrm{Tr}\left\{\Sigma_{e_0^\dagger}\right\} \\
& +2\sqrt{(|\lambda|_{\max}(I+L))^2\,\mathrm{Tr}\left\{M_0\right\}(|\lambda|_{\max}(L))^2\,\mathrm{Tr}\left\{\Sigma_{e_0^\dagger}\right\}} \\
\leq\ & \mathrm{Tr}\left\{M_0\right\}\ .
\end{aligned} \tag{167}$$

29. The inequality (167) is of the form

$$a_1x+a_2+a_3\sqrt{x}\leq x\,,$$

where $0\leq a_1<1$, and $a_2>0$, $a_3>0$, and $x:=\mathrm{Tr}\left\{M_0\right\}>0$, and it holds for

$$x\geq\left(\frac{a_3+\sqrt{a_3^2+4a_2(1-a_1)}}{2(1-a_1)}\right)^2\,.$$

Since $M_0$ has also to satisfy (163), one finally chooses

$$x\geq\max\left\{\mathrm{Tr}\left\{\Sigma_{e_0^\circ}\right\},\left(\frac{a_3+\sqrt{a_3^2+4a_2(1-a_1)}}{2(1-a_1)}\right)^2\right\}\,.$$

Hence, one concluses

$$\mathrm{Tr}\left\{\Sigma_{e_1^\circ}\right\} \quad \leq \quad \mathrm{Tr}\left\{M_0\right\} . \tag{168}$$

By a similar reasoning, also for any $k = 2, 3, \ldots$, one obtains

$$\mathrm{Tr}\left\{\Sigma_{e_k^\circ}\right\} \quad \leq \quad \mathrm{Tr}\left\{M_0\right\} , \tag{169}$$

since the sequence $\{\mathrm{Tr}\{\Sigma_{e_k^\dagger}\}, k = 0, 1, \ldots\}$ is non-increasing (see (30) and (151)). This shows that the sequence (162) is bounded.

Now, we investigate the convergence of $\mathrm{Tr}\left\{\Sigma_{e_k^\circ}\right\}$ as $k$ tends to $+\infty$. Let

$$\alpha := (|\lambda|_{\max}(I + L))^2 < 1 ,$$

and choose any $\beta > 0$ such that $\alpha + \beta < 1$. Moreover, let $M_k$ denote any symmetric and positive-semidefinite $d \times d$ matrix such that

$$\mathrm{Tr}\left\{M_k\right\} \geq \mathrm{Tr}\left\{\Sigma_{e_k^\circ}\right\} . \tag{170}$$

Then, proceeding likewise in the proof of (166) and exploting the fact that

$$|\lambda|_{\max}(L) \leq 2 ,$$

(which follows from (46) and $|\lambda|_{\max}(I + L) \leq 1$), one obtains

$$\begin{aligned}
& \mathrm{Tr}\left\{\Sigma_{e_{k+1}^\circ}\right\} \\
\leq\ & (|\lambda|_{\max}(I + L))^2\, \mathrm{Tr}\left\{M_k\right\} \\
& +(|\lambda|_{\max}(L))^2\, \mathrm{Tr}\left\{\Sigma_{e_k^\dagger}\right\} \\
& +2\sqrt{(|\lambda|_{\max}(I + L))^2\, \mathrm{Tr}\left\{M_k\right\} (|\lambda|_{\max}(L))^2\, \mathrm{Tr}\left\{\Sigma_{e_k^\dagger}\right\}} \\
\leq\ & (|\lambda|_{\max}(I + L))^2\, \mathrm{Tr}\left\{M_k\right\} \\
& +(|\lambda|_{\max}(L))^2\, \mathrm{Tr}\left\{\Sigma_{e_k^\dagger}\right\} \\
& +4\sqrt{(|\lambda|_{\max}(I + L))^2\, \mathrm{Tr}\left\{M_k\right\}\, \mathrm{Tr}\left\{\Sigma_{e_k^\dagger}\right\}} \\
\leq\ & (\alpha + \beta)\mathrm{Tr}\left\{M_k\right\} \\
\leq\ & \mathrm{Tr}\left\{M_k\right\}
\end{aligned} \tag{171}$$

for any symmetric and positive-semidefinite $d \times d$ matrix $M_k$ such that[30]

$$\begin{aligned}
\mathrm{Tr}\{M_k\} \quad & \geq \quad \max\left\{\mathrm{Tr}\{\Sigma_{e_k^\circ}\}, \frac{4\left(\sqrt{\alpha} + \sqrt{\alpha + \beta}\right)^2\ \mathrm{Tr}\{\Sigma_k^\dagger\}}{\beta^2}\right\} \\
& \geq \quad \max\left\{\mathrm{Tr}\{\Sigma_{e_k^\circ}\}, \frac{\left(2\sqrt{\alpha} + \sqrt{4\alpha + \beta(|\lambda|_{\max}(L))^2}\right)^2\ \mathrm{Tr}\{\Sigma_k^\dagger\}}{\beta^2}\right\} .
\end{aligned} \tag{172}$$

30. Likewise in footnote 29, formula (172) is obtained by reducing (171) to a quadratic inequality, expressing its solution through the two roots of the associated quadratic equality, and using (170).

Moreover, if (172) is satisfied, one gets

$$\mathrm{Tr}\left\{\Sigma_{e^{\circ}_{k+j}}\right\} \quad \leq \quad \mathrm{Tr}\left\{M_k\right\} \tag{173}$$

also for any $j = 2, 3, \ldots$, since the sequence $\{\mathrm{Tr}\{\Sigma_{e^{\dagger}_{k+j}}\}, j = 0, 1, \ldots\}$ is non-increasing. Finally, for a given value of $\beta$, we generate a sequence of symmetric and positive-semidefinite $d \times d$ matrices $M_0, M_1 \ldots$ as follows:

- $M_0$ is chosen such that (163) and (167) are satisfied;
- for $k = 0, 1, \ldots$:

$$M_{k+1} := \begin{cases} (\alpha+\beta) M_k & \text{if (172) is satisfied;} \\ M_k & \text{otherwise}. \end{cases} \tag{174}$$

By construction, one has

$$\mathrm{Tr}\{\Sigma^{\circ}_{e_k}\} \leq \mathrm{Tr}\{M_k\} \tag{175}$$

for every $k = 0, 1, \ldots$. Moreover, since, as already shown,

$$\lim_{k \to +\infty} \mathrm{Tr}\{\Sigma^{\dagger}_k\} = 0\,,$$

the first condition in (174) is satisfied an infinite number of times, which, combined with $0 < \alpha + \beta < 1$, shows that

$$\lim_{k \to +\infty} \mathrm{Tr}\{M_k\} = 0\,,$$

hence, by (175), also

$$\lim_{k \to +\infty} \mathrm{Tr}\{\Sigma^{\circ}_{e_k}\} = 0\,,$$

which, combined with (158), proves (52).

(iv) The estimate can be derived by combining the upper bound (47) on the rate of convergence to 0 of the mean-square error associated with the KF estimate of $w$ with equation (161) and the procedure to generate the matrices $M_k$ described in equation (174). ■

**Remark 18** An estimate of the rate of convergence to 0 of the MSE of the OLL estimate may be derived by combining the upper bound (47) on the rate of convergence to 0 of the MSE of the KF estimate with some results contained in the proof of Proposition 6: particularly, equation (161) and the procedure used to generate the matrices $M_k$ described in equation (174).

**Proof of Proposition 7**

By Proposition 6 (iii) and the Cauchy-Schwarz inequality, one gets

$$\begin{aligned}
& \lim_{k\to+\infty} \mathop{\mathbb{E}}_{u_k^\circ} \left\{(u_k^\circ)'(u_k^\circ)\right\} \\
= & \lim_{k\to+\infty} \mathop{\mathbb{E}}_{\hat{w}_{k+1}^\circ, \hat{w}_k^\circ} \left\{\left(\hat{w}_{k+1}^\circ - \hat{w}_k^\circ\right)'\left(\hat{w}_{k+1}^\circ - \hat{w}_k^\circ\right)\right\} \\
= & \lim_{k\to+\infty} \mathop{\mathbb{E}}_{\hat{w}_{k+1}^\circ, \hat{w}_k^\circ, w} \left\{\left((\hat{w}_{k+1}^\circ - w) + (w - \hat{w}_k^\circ)\right)'\left((\hat{w}_{k+1}^\circ - w) + (w - \hat{w}_k^\circ)\right)\right\} \\
= & \lim_{k\to+\infty} \mathop{\mathbb{E}}_{\hat{w}_{k+1}^\circ, w} \left\{\left((\hat{w}_{k+1}^\circ - w)\right)'\left((\hat{w}_{k+1}^\circ - w)\right)\right\} \\
& + \lim_{k\to+\infty} \mathop{\mathbb{E}}_{\hat{w}_k^\circ, w} \left\{\left((\hat{w}_k^\circ - w)\right)'\left((\hat{w}_k^\circ - w)\right)\right\} \\
& +2 \lim_{k\to+\infty} \mathop{\mathbb{E}}_{\hat{w}_{k+1}^\circ, \hat{w}_k^\circ, w} \left\{\left((\hat{w}_{k+1}^\circ - w)\right)'\left((\hat{w}_k^\circ - w)\right)\right\} \\
\leq & \lim_{k\to+\infty} \mathop{\mathbb{E}}_{\hat{w}_{k+1}^\circ, w} \left\{\left((\hat{w}_{k+1}^\circ - w)\right)'\left((\hat{w}_{k+1}^\circ - w)\right)\right\} \\
& + \lim_{k\to+\infty} \mathop{\mathbb{E}}_{\hat{w}_k^\circ, w} \left\{\left((\hat{w}_k^\circ - w)\right)'\left((\hat{w}_k^\circ - w)\right)\right\} \\
& +2 \lim_{k\to+\infty} \sqrt{\mathop{\mathbb{E}}_{\hat{w}_{k+1}^\circ, w} \left\{\left((\hat{w}_{k+1}^\circ - w)\right)'\left((\hat{w}_{k+1}^\circ - w)\right)\right\}} \sqrt{\mathop{\mathbb{E}}_{\hat{w}_k^\circ, w} \left\{\left((\hat{w}_k^\circ - w)\right)'\left((\hat{w}_k^\circ - w)\right)\right\}} \\
= & \ 0 \, .
\end{aligned} \tag{176}$$

Hence, the proof is concluded likewise for the KF estimates (see formulas (57) and (58) in Subsection 6.1), exploiting again the fact that, by Proposition 6 (iii), one has $\lim_{k\to+\infty} \mathrm{MSE}_k^\circ = 0$. ■

**Proof of Proposition 8**

By the optimality for the problem with $\gamma > 0$ and the limit problem, respectively, one has (with $e_0^\circ = e_0^\dagger = e_0$)

$$\begin{aligned}
& \mathop{\mathbb{E}}_{\{e_k^\circ\}_{k=0}^N, \{u_k^\circ\}_{k=0}^{N-1}, \{Q_k\}_{k=0}^N} \left\{\sum_{k=0}^{N-1} \left[(e_k^\circ)'Q_k(e_k^\circ) + \gamma(u_k^\circ)'(u_k^\circ)\right] + (e_N^\circ)'Q_N(e_N^\circ)\right\} \\
\leq & \mathop{\mathbb{E}}_{\{e_k^\dagger\}_{k=0}^N, \{u_k^\dagger\}_{k=0}^{N-1}, \{Q_k\}_{k=0}^N} \left\{\sum_{k=0}^{N-1} \left[(e_k^\dagger)'Q_k(e_k^\dagger) + \gamma(u_k^\dagger)'(u_k^\dagger)\right] + (e_N^\dagger)'Q_N(e_N^\dagger)\right\} ,
\end{aligned} \tag{177}$$

and

$$\begin{aligned}
& \mathop{\mathbb{E}}_{\{e_k^\dagger\}_{k=0}^N, \{Q_k\}_{k=0}^N} \left\{\sum_{k=0}^{N-1} \left[(e_k^\dagger)'Q_k(e_k^\dagger)\right] + (e_N^\dagger)'Q_N(e_N^\dagger)\right\} \\
\leq & \mathop{\mathbb{E}}_{\{e_k^\circ\}_{k=0}^N, \{Q_k\}_{k=0}^N} \left\{\sum_{k=0}^{N-1} \left[(e_k^\circ)'Q_k(e_k^\circ)\right] + (e_N^\circ)'Q_N(e_N^\circ)\right\} .
\end{aligned} \tag{178}$$

Hence, combining (177) and (178), one obtains

$$\begin{aligned}
&\mathop{\mathbb{E}}_{\{u_k^\circ\}_{k=0}^{N-1}}\left\{\sum_{k=0}^{N-1}\left[\gamma(u_k^\circ)'(u_k^\circ)\right]\right\}\\
\le\ &\mathop{\mathbb{E}}_{\{u_k^\circ\}_{k=0}^{N-1}}\left\{\sum_{k=0}^{N-1}\left[\gamma(u_k^\circ)'(u_k^\circ)\right]\right\}+\mathop{\mathbb{E}}_{\{e_k^\circ\}_{k=0}^{N},\{Q_k\}_{k=0}^{N}}\left\{\sum_{k=0}^{N-1}\left[(e_k^\circ)'Q_k(e_k^\circ)\right]+(e_N^\circ)'Q_N(e_N^\circ)\right\}\\
-\ &\mathop{\mathbb{E}}_{\{e_k^\dagger\}_{k=0}^{N},\{Q_k\}_{k=0}^{N}}\left\{\sum_{k=0}^{N-1}\left[(e_k^\dagger)'Q_k(e_k^\dagger)\right]+(e_N^\dagger)'Q_N(e_N^\dagger)\right\}\\
\le\ &\mathop{\mathbb{E}}_{\{u_k^\dagger\}_{k=0}^{N-1}}\left\{\sum_{k=0}^{N-1}\left[\gamma(u_k^\dagger)'(u_k^\dagger)\right]\right\}. \qquad (179)
\end{aligned}$$

Concluding, as $\gamma > 0$, one gets the result. ■

**Proof of Proposition 9**

We recall here the equations (25), (24), and (23), which are needed to compute the KF estimate $\hat{w}_k^\dagger$:

$$\Sigma_{k+1} = \Sigma_k - \Sigma_k(C_{k+1}^{(\phi)})'((C_{k+1}^{(\phi)})\Sigma_k(C_{k+1}^{(\phi)})' + \sigma_\varepsilon^2)^{-1}(C_{k+1}^{(\phi)})\Sigma_k\,, \qquad (180)$$

$$H_{k+1} := \Sigma_{k+1}(C_{k+1}^{(\phi)})'(\sigma_\varepsilon^2)^{-1}\,, \qquad (181)$$

$$\hat{w}_{k+1}^\dagger = \hat{w}_k^\dagger + H_{k+1}(y_{k+1} - (C_{k+1}^{(\phi)})\hat{w}_k^\dagger)\,, \qquad (182)$$

where we have used the notation $C_k^{(\phi)}$ to recall that they are now defined using $\phi(x_k)$ instead than $x_k$, i.e., one has

$$C_k^{(\phi)} := (\phi(x_k))'\,.$$

We also recall the initializations

$$\hat{w}_{-1}^\dagger = 0$$

and

$$\Sigma_{-1} = \Sigma_w\,.$$

Then, for $k = -1$, using (180), (181), and (182), we obtain

$$\begin{aligned}
\Sigma_0 &= \nu I_{d_E} - \nu I_{d_E}\phi(x_0)\left(\nu(\phi(x_0))'\phi(x_0) + \sigma_\varepsilon^2\right)^{-1}(\phi(x_0))'\nu I_{d_E}\\
&= \nu I_{d_E} - \nu^2\phi(x_0)\left(\nu\mathcal{K}(x_0,x_0) + \sigma_\varepsilon^2\right)^{-1}(\phi(x_0))'\,,
\end{aligned}$$

$$\begin{aligned}
H_0 &= \left(\nu I_{d_E} - \nu^2\phi(x_0)\left(\nu\mathcal{K}(x_0,x_0) + \sigma_\varepsilon^2\right)^{-1}(\phi(x_0))'\right)\phi(x_0)(\sigma_\varepsilon^2)^{-1}\\
&= \left(\nu - \nu^2\left(\nu\mathcal{K}(x_0,x_0) + \sigma_\varepsilon^2\right)^{-1}\mathcal{K}(x_0,x_0)\right)(\sigma_\varepsilon^2)^{-1}\phi(x_0)\,,
\end{aligned}$$

$$\hat{w}_0^\dagger = \left(\nu - \nu^2 \left(\nu \mathcal{K}(x_0, x_0) + \sigma_\varepsilon^2\right)^{-1} \mathcal{K}(x_0, x_0)\right) (\sigma_\varepsilon^2)^{-1} y_0 \phi(x_0)\,,$$

and

$$\begin{aligned} (\hat{w}_0^\dagger)' \phi(x) &= \left(\nu - \nu^2 \left(\nu \mathcal{K}(x_0, x_0) + \sigma_\varepsilon^2\right)^{-1} \mathcal{K}(x_0, x_0)\right) (\sigma_\varepsilon^2)^{-1} y_0 (\phi(x_0))' \phi(x)\,, \\ &= \left(\nu - \nu^2 \left(\nu \mathcal{K}(x_0, x_0) + \sigma_\varepsilon^2\right)^{-1} \mathcal{K}(x_0, x_0)\right) (\sigma_\varepsilon^2)^{-1} y_0 \mathcal{K}(x_0, x)\,, \end{aligned}$$

where the last expression does not involve an explicit computation of $\phi(x_0)$ and $\phi(x)$. Similarly, by iterating the procedure above, one obtains, for every $k$, expressions of the form (see [33, Theorem 2])

$$\Sigma_k = \nu_k I_{d_E} - \Phi_k \Psi_k \Phi_k'\,, \tag{183}$$

where $\nu_k$ and $\Psi_k$ are, respectively, suitable scalars and suitable $(k+1) \times (k+1)$ matrices (both of which can be computed recursively), and

$$\Phi_k := \left(\phi(x_0), \ldots, \phi(x_k)\right)\,.$$

This, combined with (181) and (182), respectively, provides (63), (64), and (65). ■

**Proof of Proposition 10**

(i) Following Remark 10, the matrix (operator) $L_k$ is proportional to $I_{d_E}$, say

$$L_k = \alpha_{L_k} I_{d_E}$$

for some $\alpha_{L_k} > 0$, then the update equation (22) of the OLL estimate of $w$ becomes

$$\begin{aligned} \hat{w}_{k+1}^\circ &= \hat{w}_k^\circ + \alpha_{L_k} I_{d_E} (\hat{w}_k^\circ - \hat{w}_k^\dagger) \\ &= \hat{w}_k^\circ + \alpha_{L_k} (\hat{w}_k^\circ - \hat{w}_k^\dagger)\,. \end{aligned} \tag{184}$$

This, combined with the initialization

$$\hat{w}_0^\circ = 0$$

and with (64), shows that, for $k = 1, 2, \ldots$, one gets (69) and (70).

(ii) Applying KPCA, one can show that the eigenvectors associated with the positive eigenvalues of $\overline{Q}_{\text{emp}}^{(\phi)}$ are linear combinations of $\phi(\tilde{x}_1), \ldots, \phi(\tilde{x}_{l_U})$. Hence, following again Remark 10, one obtains

$$\begin{aligned} & L_k \\ = & \sum_{j=1}^{l_U} (j\text{-th lin. comb. of } \phi(\tilde{x}_1), \ldots, \phi(\tilde{x}_{l_U}))\, (\Lambda_{L_k})_{(j,j)}\, (\text{same} j\text{-th lin. comb. of } \phi(\tilde{x}_1), \ldots, \phi(\tilde{x}_{l_U}))'\,. \end{aligned}$$

Then, reasoning as in the proof of part (i), one obtains (71) and (72). ■